\documentclass[11pt]{amsart}
\usepackage{amssymb}
\usepackage{amsmath}
\usepackage{mathrsfs}
\usepackage{fancyhdr}
\usepackage[british]{babel}
\usepackage{geometry}
\usepackage{enumitem}
\usepackage{algpseudocode}
\usepackage{dsfont}
\usepackage{centernot}
\usepackage{xstring}
\usepackage{colortbl}
\usepackage[section]{algorithm}
\usepackage{graphicx}
\usepackage[font={footnotesize}]{caption}
\usepackage[usenames,dvipsnames,table]{xcolor}
\usepackage{tikz}
\usepackage[h]{esvect}
\usepackage[
		bookmarksopen=true,
		bookmarksopenlevel=1,
		colorlinks=true,
		linkcolor=darkblue,
        linktoc=page,
		citecolor=darkblue,
]{hyperref}
\usepackage{bbm}

\title[]{A short proof of the blow-up lemma\\for approximate decompositions}
\date{\today}
\author[S.~Ehard]{Stefan Ehard}
\address[S.~Ehard]{Institut f\"ur Optimierung und Operations Research, Universit\"at Ulm, 89081 Ulm,
Germany
}
\email{stefan.ehard@uni-ulm.de}

\author[F.~Joos]{Felix Joos}
\address[F.~Joos]{Department of Mathematics, Universit\"at Hamburg, Bundesstrasse 55,
20146 Hamburg, Germany }
\email{felix.joos@uni-hamburg.de}

\thanks{The research leading to these results was partially supported by the Deutsche Forschungsgemeinschaft (DFG, German Research Foundation) -- 339933727 (F.~Joos).
}

\newtheorem{theorem}[algorithm]{Theorem}

\newtheorem{lemma}[algorithm]{Lemma}

\newtheorem{fact}[algorithm]{Fact}

\theoremstyle{definition}

\newtheorem{defin}[algorithm]{Definition}

\newcounter{stepenv}
\newenvironment{stepenv}[1][]{\refstepcounter{stepenv}}{}

\newcounter{step}[stepenv]
\newenvironment{step}[1][]{\refstepcounter{step}\par\medskip\noindent%
        \textit{Step~\thestep. #1} \itshape\rmfamily}{\medskip}
				
\newcounter{substep}[step]
\renewcommand{\thesubstep}{\thestep.\arabic{substep}}
\newenvironment{substep}[1][]{\refstepcounter{substep}\par\medskip\noindent%
        \emph{Step~\thesubstep. #1} \itshape\rmfamily}{\medskip}

\newcounter{claim}[stepenv]
\newenvironment{claim}[1][]{\refstepcounter{claim}\par\medskip\noindent%
        \textit{Claim~\theclaim. #1} \itshape\rmfamily}{\medskip}

\numberwithin{equation}{section}

\definecolor{darkblue}{rgb}{0,0,0.5}

\def\noproof{{\unskip\nobreak\hfill\penalty50\hskip2em\hbox{}\nobreak\hfill%
       $\square$\parfillskip=0pt\finalhyphendemerits=0\par}\goodbreak}
\def\endproof{\noproof\bigskip}

\def\noclaimproof{{\unskip\nobreak\hfill\penalty50\hskip2em\hbox{}\nobreak\hfill%
       $-$\parfillskip=0pt\finalhyphendemerits=0\par}\goodbreak}

\def\endclaimproof{\noclaimproof\medskip}

\newdimen\margin
\def\textno#1&#2\par{
   \margin=\hsize
   \advance\margin by -4\parindent
          \setbox1=\hbox{\sl#1}
   \ifdim\wd1 < \margin
      $$\box1\eqno#2$$
   \else
      \bigbreak
      \hbox to \hsize{\indent$\vcenter{\advance\hsize by -3\parindent
      \it\noindent#1}\hfil#2$}
      \bigbreak
   \fi}

\def\proof{\removelastskip\penalty55\medskip\noindent\begin{stepenv}\end{stepenv}{\bf Proof. }} % in each main proof, claim and step counter set back

\def\lateproof#1{\removelastskip\penalty55\medskip\noindent\begin{stepenv}\end{stepenv}{\bf Proof of #1. }} % in each main proof, claim and step counter set back

\DeclareMathOperator{\dg}{deg}

\def\claimproof{\removelastskip\penalty55\medskip\noindent{\em Proof of claim: }}

\begin{document}

\newcommand{\new}[1]{\textcolor{red}{#1}}
\def\COMMENT#1{}
\def\TASK#1{}

\newcommand{\todo}[1]{\begin{center}\textbf{to do:} #1 \end{center}}

\def\eps{{\varepsilon}}
\def\heps{{\hat{\varepsilon}}}

\newcommand{\ex}{\mathbb{E}}
\newcommand{\pr}{\mathbb{P}}
\newcommand{\cB}{\mathcal{B}}
\newcommand{\cA}{\mathcal{A}}
\newcommand{\cE}{\mathcal{E}}
\newcommand{\cS}{\mathcal{S}}
\newcommand{\cF}{\mathcal{F}}
\newcommand{\cG}{\mathcal{G}}
\newcommand{\bL}{\mathbb{L}}
\newcommand{\bF}{\mathbb{F}}
\newcommand{\bZ}{\mathbb{Z}}
\newcommand{\cH}{\mathcal{H}}
\newcommand{\cC}{\mathcal{C}}
\newcommand{\cM}{\mathcal{M}}
\newcommand{\bN}{\mathbb{N}}
\newcommand{\bR}{\mathbb{R}}
\def\O{\mathcal{O}}
\newcommand{\cP}{\mathcal{P}}
\newcommand{\cQ}{\mathcal{Q}}
\newcommand{\cR}{\mathcal{R}}
\newcommand{\cJ}{\mathcal{J}}
\newcommand{\cL}{\mathcal{L}}
\newcommand{\cK}{\mathcal{K}}
\newcommand{\cD}{\mathcal{D}}
\newcommand{\cI}{\mathcal{I}}
\newcommand{\cN}{\mathcal{N}}
\newcommand{\cV}{\mathcal{V}}
\newcommand{\cT}{\mathcal{T}}
\newcommand{\cU}{\mathcal{U}}
\newcommand{\cX}{\mathcal{X}}
\newcommand{\cZ}{\mathcal{Z}}
\newcommand{\cW}{\mathcal{W}}
\newcommand{\1}{{\bf 1}_{n\not\equiv \delta}}
\newcommand{\eul}{{\rm e}}
\newcommand{\Erd}{Erd\H{o}s}
\newcommand{\cupdot}{\mathbin{\mathaccent\cdot\cup}}
\newcommand{\whp}{whp }
\newcommand{\bX}{\mathcal{X}}
\newcommand{\bV}{\mathcal{V}}
\newcommand{\hcX}{\widehat{\mathcal{X}}}
\newcommand{\hcV}{\widehat{\mathcal{V}}}
\newcommand{\hpsi}{\widehat{\psi}}
\newcommand{\hA}{\widehat{A}}
\newcommand{\tA}{\widetilde{A}}
\newcommand{\tB}{\widetilde{B}}
\newcommand{\tH}{\widetilde{H}}
\newcommand{\tG}{\widetilde{G}}
\newcommand{\hB}{\widehat{B}}
\newcommand{\tE}{\widetilde{E}}
\newcommand{\hX}{\widehat{X}}
\newcommand{\hV}{\widehat{V}}
\newcommand{\tX}{\widetilde{X}}
\newcommand{\tV}{\widetilde{V}}
\newcommand{\oX}{\overline{X}}
\newcommand{\cbI}{\overline{\mathcal{I}^\alpha}}
\newcommand{\hAj}{\widehat{A}^\sigma_j}
\newcommand{\hVM}{V^M}
\newcommand{\bd}{\mathbf{d}}
\newcommand{\sX}{\mathscr{X}}
\newcommand{\sA}{\mathscr{A}}
\newcommand{\sB}{\mathscr{B}}
\newcommand{\sP}{\mathscr{P}}
\newcommand{\stA}{{\mathscr{\widetilde{A}}}}

\newcommand{\doublesquig}{%
  \mathrel{%
    \vcenter{\offinterlineskip
      \ialign{##\cr$\rightsquigarrow$\cr\noalign{\kern-1.5pt}$\rightsquigarrow$\cr}%
    }%
  }%
}

\newcommand{\defn}{\emph}

\newcommand\restrict[1]{\raisebox{-.5ex}{$|$}_{#1}}

\newcommand{\prob}[1]{\mathrm{\mathbb{P}}\left[#1\right]}
\newcommand{\expn}[1]{\mathrm{\mathbb{E}}\left[#1\right]}
\def\gnp{G_{n,p}}
\def\G{\mathcal{G}}
\def\lflr{\left\lfloor}
\def\rflr{\right\rfloor}
\def\lcl{\left\lceil}
\def\rcl{\right\rceil}

\newcommand{\qbinom}[2]{\binom{#1}{#2}_{\!q}}
\newcommand{\binomdim}[2]{\binom{#1}{#2}_{\!\dim}}

\newcommand{\grass}{\mathrm{Gr}}

\newcommand{\brackets}[1]{\left(#1\right)}
\def\sm{\setminus}
\newcommand{\Set}[1]{\{#1\}}
\newcommand{\set}[2]{\{#1\,:\;#2\}}
\newcommand{\krq}[2]{K^{(#1)}_{#2}}
\newcommand{\ind}[1]{$\mathbf{S}(#1)$}
\newcommand{\indC}[1]{$\mathbf{C}(#1)$}
\newcommand{\indcov}[1]{$(\#)_{#1}$}
\def\In{\subseteq}
\newcommand{\dom}{dom}
\newcommand{\IND}{\mathbbm{1}}
\newcommand{\norm}[1]{\|#1\|}
\newcommand{\normv}[1]{\|#1\|_v}
\newcommand{\normV}[2]{\|#1\|_{#2}}

\begin{abstract} 
	\noindent
	Kim, K\"uhn, Osthus and Tyomkyn (Trans.~Amer.~Math.~Soc. 371 (2019), 4655--4742)
greatly extended the well-known blow-up lemma of Koml\'os, S\'ark\"ozy and Szemer\'edi by 
 proving a `blow-up lemma for approximate decompositions' which states that multipartite quasirandom graphs can be almost decomposed into any collection of bounded degree graphs with the same multipartite structure and slightly fewer edges.
	This result has already been used by Joos, Kim, K\"uhn and Osthus to prove the tree packing conjecture due to Gy\'arf\'as and Lehel from 1976 and Ringel's conjecture from 1963 for bounded degree trees as well as
	 implicitly in the recent resolution of the Oberwolfach problem (asked by Ringel in 1967) by Glock, Joos, Kim, K\"uhn and Osthus.
	
	Here we present a new and significantly shorter proof of the blow-up lemma for approximate decompositions.
	In fact, we prove a more general theorem that yields packings with stronger quasirandom properties
	so that it can be combined with Keevash's results on designs to obtain results of the following form.
	For all $\eps>0$, $r\in \bN$ and all large $n$ (such that $r$ divides $n-1$),
	there is a decomposition of $K_n$ into any collection of $r$-regular graphs $H_1,\ldots,H_{(n-1)/r}$ on $n$ vertices
	provided that $H_1,\ldots,H_{\eps n}$ contain each at least $\eps n$ vertices in components of size at most $\eps^{-1}$.
\end{abstract}

\maketitle
\vspace{-0.4cm}
\section{Introduction}
The theme of decomposing `large' objects into `smaller' objects
or finding a maximal number of specified `small' objects in a `larger' object is among the key topics in mathematics.
In discrete mathematics, it appears in Euler's question from 1782 for which~$n$ there exist pairs of orthogonal Latin squares\footnote{This is equivalent to a $K_4$-decomposition of the $4$-partite graph $K_{n,n,n,n}$.} of order $n$,
%(this is equivalent to a $K_4$-decomposition of a balanced $4$-partite graph),
in Steiner's questions for Steiner systems from the 1850s which cumulated in the `existence of designs' question,
in Walecki's theorem on decompositions of complete graphs into edge-disjoint Hamilton cycles from the 1890s,
and in Kirkman's famous `school girl problem'.

These questions and results set off an entire branch of combinatorics and design theory.
Several decades later in the 1970s, Wilson~\cite{wilson:72a,wilson:72b,wilson:75} famously proved that (the edge set of) the complete graph on $n$ vertices can be decomposed into any fixed graph provided necessary divisibility conditions are satisfied and $n$ is large, thereby solving the `existence of designs' question for graphs.
In~2014, Keevash verified the `existence of designs' for hypergraphs~\cite{keevash:14}.
This has been reproved and {generalised} by Glock, K\"uhn, Lo and Osthus in~\cite{GKLO:16,GKLO:17}.  Keevash extended his results to a more general framework in~\cite{keevash:18b}.

In contrast to these questions and results where we aim to decompose a large graph into graphs of fixed size,
one can also ask for decompositions into larger pieces, for example into graphs with the same number of vertices as the host graph.
A prime example is the Oberwolfach problem where Ringel asked in 1967 whether one can decompose (the edge set of) $K_{2n+1}$ into $n$ copies of any $2$-regular graph on $2n+1$ vertices.
This problem received considerable attention and Glock, Joos, Kim, K\"uhn and Osthus~\cite{GJKKO:18} solved it for large~$n$.
Possibly equally well-known is Ringel's conjecture from 1963 stating that $K_{2n+1}$ can be decomposed into any tree with $n$ edges, as well as the tree packing conjecture due to Gy\'arf\'as and Lehel from 1976 stating that $K_n$ can be decomposed into any collection of trees $T_1,\ldots,T_n$, where $T_i$ has $i$ vertices.
Ringel's conjecture has been solved by Montgomery, Pokrovskiy and Sudakov~\cite{MPS:20,MPS:18} and both conjectures have been solved for bounded degree trees by Joos, Kim, K\"uhn and Osthus~\cite{JKKO:ta}; Allen, B\"ottcher, Clemens and Taraz~\cite{ABCT:19} have solved these conjectures for trees with many leaves and maximum degree $o(n/\log n)$ (in fact, they proved a more general result on degenerate graphs) -- all mentioned results apply only when $n$ is sufficiently large. We refer the reader to~\cite{BHPT:16,FLM:17,KKOT:19,MRS:16} for earlier results regarding these conjectures and to~\cite{ABHP:17,CKKO:19,FS:ta,KO:13} for further developments in the field.

On a very high level, numerous decomposition results combine approximate decomposition results with certain absorbing techniques.
This includes~\cite{ABCT:19,GKLO:16,GKLO:17,JKKO:ta,keevash:14,keevash:18b}.
For many questions in extremal combinatorics, 
the blow-up lemma due to Koml\'os, S\'ark\"ozy and Szemer\'edi~\cite{KSS:97} in combination with Szemer\'edi's regularity lemma has demonstrated its power and usefulness.
Having this in mind and in need of a powerful approximate decomposition result, 
Kim, K\"uhn, Osthus and Tyomkyn~\cite{KKOT:19} proved a far-reaching generalisation of the blow-up lemma, a `blow-up lemma for approximate decompositions', which can also be combined with the regularity lemma to obtain almost decompositions of graphs into bounded degree graphs.

The blow-up lemma for approximate decompositions has already exhibited its versatility. 
It has been applied in~\cite{JKKO:ta,KKL:18,KLMY:17}
and in~\cite{CKKO:19} for a `bandwidth theorem for approximate decompositions', which in turn is one of the key ingredients for the resolution of the Oberwolfach problem in~\cite{GJKKO:18}.
However, its very complex and long proof is an obstacle for further generalisations.
One main aim of this paper is to overcome this by presenting a new and significantly shorter proof.

Our approach makes it possible to include some more features
including an easier handling of exceptional vertices, which results in a substantially easier applicability of the theorem,
and the approximate decompositions share stronger quasirandom properties.
{To be more precise, the first yields shorter proofs of the main results in~\cite{CKKO:19} and~\cite{JKKO:ta} as certain technically involved preprocessing steps are no longer needed;}
the latter permits to combine our main result with Keevash's recent results on designs~\cite{keevash:18b}.
We dedicate a section at the end of the paper to demonstrate this and obtain new results on decomposing quasirandom graphs into regular spanning graphs.

\subsection{The blow-up lemma for approximate decompositions}

In this section, we first introduce some terminology and then state the blow-up lemma for approximate decompositions.
We say that a collection/multiset of graphs $\cH=\{H_1,\ldots, H_s\}$ \defn{packs} into a graph $G$
if there is a function $\phi: \bigcup_{H\in \cH}V(H)\to V(G)$ such that $\phi|_{V(H)}$ is injective and $\phi$ injectively maps edges onto edges. In such a case, we call $\phi$ a \emph{packing of $\cH$ into $G$}.
Our general aim is to pack a collection $\cH$ of multipartite graphs in a host graph $G$ having the same multipartite structure which is captured by a so-called `reduced graph' $R$.
To this end, let
$(H,G,R,\cX,\cV)$ be a \defn{blow-up instance} if
\begin{itemize}
	\item $H,G,R$ are graphs where $V(R)=[r]$ for some $r\geq 2$;
	\item $\cX=(X_i)_{i\in[r]}$ is a vertex partition of $H$ into independent sets,
	$\cV=(V_i)_{i\in[r]}$ is a vertex partition of $G$ such that $|V_i|=|X_i|$ for all $i\in [r]$; 
	\item $H[X_{i},X_j]$ is empty whenever {$ij\in \binom{[r]}{2}\sm E(R)$.}
\end{itemize}
%$(H,G,R,(X_i)_{i\in[r]},(V_i)_{i\in[r]})$ be an \defn{blow-up instance} if
%\begin{itemize}
%	\item $H$ and $G$ are graphs, $(X_i)_{i\in[r]}$ is a vertex partition of $H$ into independent sets,
%	$(V_i)_{i\in[r]}$ is a vertex partition of $G$, and $|V_i|=|X_i|$ for all $i\in [r]$;
%	\item $R$ is a graph with vertex set $[r]$ and $H[X_{i},X_j]$ is empty whenever $ij\notin E(R)$.
%\end{itemize}
We also refer to $\sB=(\cH,G,R,\cX,\cV)$ as a \defn{blow-up instance}
if $\cH$ is a collection of graphs and $\cX$ is a collection of vertex partitions $(X_i^H)_{i\in [r],H\in \cH}$ 
so that $(H,G,R,(X_i^H)_{i\in[r]},\cV)$ is a blow-up instance for every $H\in \cH$.

The quasirandom notion mostly used in this paper coincides with the notion used in Szemeredi's regularity lemma.
For a bipartite graph $G$ with vertex partition $(V_1,V_2)$,
we define the \defn{density} of $W_1,W_2$ with $W_i\In V_i$ by $d_G(W_1,W_2):=e_G(W_1,W_2)/|V_1||V_2|$.
We say $G$ is \defn{$(\eps,d)$-regular} if $d_G(W_1,W_2)=d\pm\eps$ for all $W_i\In V_i$ with $|W_i|\geq \eps |V_i|$,
and $G$ is \defn{$(\eps,d)$-super-regular} if in addition $|N_G(v)\cap V_{3-i}|= (d\pm \eps)|V_{3-i}|$ for each $i\in [2]$ and $v\in V_i$.
The blow-up instance $\sB$ is \defn{$(\eps,d)$-super-regular} if $G[V_i,V_j]$ is $(\eps,d)$-super-regular for all $ij\in E(R)$,
and $\sB$ is \defn{$\Delta$-bounded} if $\Delta(R),\Delta(H) \leq \Delta$ for each $H\in \cH$.
Now we are ready to state the blow-up lemma for approximate decompositions.

\begin{theorem}[Kim, K\"uhn, Osthus, Tyomkyn~\cite{KKOT:19}]\label{thm:main_extended}
	%For all $\alpha, d_0\in (0,1]$ and $\Delta,r\geq 2$, 
	For all $\alpha \in (0,1]$ and $r\geq 2$, 
	there exist $\eps=\eps(\alpha) > 0$ and $n_0 = n_0(\alpha, r)$ such that the following holds for all $n \geq n_0$ and $d\geq \alpha$. 
	Suppose $(\cH,G,R,\cX,\cV)$ is an $(\eps,d)$-super-regular and $\alpha^{-1}$-bounded blow-up instance
	such that $|V_i|=n$ for all $i\in [r]$, $|\cH|\leq \alpha^{-1}n$, and $\sum_{H\in \cH} e_H(X_{i}^H,X_j^H) \leq (1-\alpha)dn^2$ for all $ij \in E(R)$. %and $H^i[X_j^i,X_k^i]$ is empty if $jk \notin E(R)$. 
	Then there is a packing $\phi$ of $\cH$ into $G$ such that $\phi(X_i^H)= V_i$ for all $i\in [r]$ and $H\in \cH$.
\end{theorem}

We remark that there are more general versions of Theorem~\ref{thm:main_extended} in~\cite{KKOT:19}, but omit the more technical statements here.
Instead we state our main result and the interested reader can easily check that it generalises\footnote{Observe that we do not allow different densities between the cluster pairs in $G$. However, this technical complication could very easily be implemented by adding at numerous places extra indices. As this feature has never been used so far in applications, we omitted it for the sake of a clearer presentation.} the more technical versions in~\cite{KKOT:19}.

\subsection{Main result}

Most blow-up lemmas exhibit their power if they are applied in conjunction with Szemer\'edi's regularity lemma.
This, however, comes with the expense of a small set of vertices over which we have no control.
Consequently, in such a setting, when embedding a graph $H$ into $G$, it is often the case that some vertices of $H$ are already embedded and the blow-up lemma is applied only to some nice part of~$G$. 
To deal with such scenarios we consider extended blow-up instances.
We say $(H,G,R,\cX,\cV,\phi_0)$ is an \defn{extended blow-up instance}
if
%\begin{itemize}
%	\item $H$ and $G$ are graphs, $(X_i)_{i\in[r]_0}$ is a vertex partition of $H$ into independent sets,
%	$(V_i)_{i\in[r]_0}$ is a vertex partition of $G$ and $|V_i|=|X_i|$ for all $i\in [r]$;
%	\item $R$ is a graph with vertex set $[r]$ and $H[X_{i},X_j]$ is empty whenever $ij\notin E(R)$;
%	\item $\phi_0$ is an embedding of $X_0$ into $V_0$.
%\end{itemize}
\begin{itemize}
	\item $H,G,R$ are graphs where $V(R)=[r]$ for some $r\geq 2$;
	\item $\cX=(X_i)_{i\in[r]_0}$ is a vertex partition of $H$ into independent sets,
	$\cV=(V_i)_{i\in[r]_0}$ is a vertex partition of $G$ such that $|V_i|=|X_i|$ for all $i\in [r]_0$; 
	\item $H[X_{i},X_j]$ is empty whenever {$ij\in \binom{[r]}{2}\sm E(R)$;}
	\item $\phi_0$ is an {injective} embedding of $X_0$ into $V_0$.
\end{itemize}

This definition also extends as above to the case when $H$ is replaced by a collection of graphs~$\cH$ in the obvious way as before.
%We refer to $(\cH,G,R,\cX,(V_i)_{i\in[r]_0},\phi_0)$ also as an \defn{extended blow-up instance}
%if $\cH$ is a collection of graphs and $\cX$ is a collection of vertex partitions $(X_i^H)_{i\in [r]_0,H\in \cH}$ 
%so that $(H,G,R,(X_i^H)_{i\in[r]_0},(V_i)_{i\in[r]_0},\phi_0|_{X_0^H})$ is an extended blow-up instance for every $H\in \cH$.
An extended blow-up instance is \defn{$(\eps,d)$-super-regular} if $G[V_i,V_j]$ is $(\eps,d)$-super-regular for all $ij\in E(R)$.

Let $\sB=(\cH,G,R,\cX,\cV,\phi_0)$ be an extended blow-up instance. We say $\sB$ is \defn{$(\eps,\alpha)$-linked} if 
\begin{itemize}
	\item at most $\eps|X^H_i|$ vertices in $X^H_i$ have a neighbour in $X_0^H$ for all $i\in[r]$, $H\in \cH$;
	\item $|V_i\cap\bigcap_{x_0\in X_0^H\cap N_H(x)} N_G(\phi_0(x_0))| \geq \alpha |V_i|$ for all $x\in X_i^H, i\in [r], H\in \cH$;
	\item $|\phi_0^{-1}(v)|\leq \eps |\cH|$ for all $v\in V_0$;
\item$\sum_{H\in \cH}|N_H(x_0^H)\cap N_H(x_0'^H) \cap X_i^H|\leq \eps|V_i|^{1/2}$ for all $i\in[r]$ and distinct $v_0, v_0' \in V_0$ 
where $x_0^H=\phi_0^{-1}(v_0)\cap X_0^H$ and $x_0'^H=\phi_0^{-1}(v_0')\cap X_0^H$ for $H\in\cH$.
\end{itemize}
One feature of our result is that one can replace `blow-up instance' in Theorem~\ref{thm:main_extended} by `extended blow-up instance that is $(\eps,\alpha)$-linked'.
We remark that the above conditions are easily met in applications known to us and are similar to conditions found elsewhere for this purpose.

Next, we define two types of {structures} for $\sB$ and our main result yields a packing that behaves as we would expect it from an idealised typical random packing with respect to these {structures}.
We say $(W,Y_1,\ldots,Y_k)$ is an \defn{$\ell$-set tester} for $\sB$ 
if $k\leq \ell$ and there exist $i\in [r]$ and distinct $H_1,\ldots,H_k\in \cH$ such that $W\In V_i$ and $Y_j\In X_i^{H_j}$ for all $j\in [k]$.
We say $(v,\omega)$ is an \defn{$\ell$-vertex tester} for $\sB$ if $v\in V_i$ and $\omega:\bigcup_{H\in \cH}X_i^H \to [0,\ell]$ for some $i\in [r]$. 
For a weight function $\omega$ on a finite set $X$, we define $\omega(X'):=\sum_{x\in X'}\omega(x)$ for any $X'\In X$.
The following theorem is our main result.

\begin{theorem}\label{thm:main_new}
For all $\alpha \in (0,1]$ and $r\geq 2$, 
there exist $\eps=\eps(\alpha) > 0$ and $n_0 = n_0(\alpha, r)$ such that the following holds for all $n \geq n_0$ and $d\geq \alpha$. 
Suppose $(\cH,G,R,\cX,\cV,\phi_0)$ is an  $(\eps,d)$-super-regular, $\alpha^{-1}$-bounded and $(\eps,\alpha)$-linked extended blow-up instance,
$|V_i|=(1\pm \eps)n$ for all $i\in [r]$, 
$|\cH|\leq \alpha^{-1}n$, and $\sum_{H\in \cH} e_H(X_{i}^H,X_j^H) \leq (1-\alpha)dn^2$ for all $ij \in E(R)$.
Suppose $\cW_{set},\cW_{ver}$ are sets of  $\alpha^{-1}$-set testers and $\alpha^{-1}$-vertex testers of size at most $n^{\log n}$, respectively.
Then there is a packing $\phi$ of $\cH$ into $G$ which extends $\phi_0$ 
such that 
\begin{enumerate}[label=\rm (\roman*)]
	\item\label{item:partition} $\phi(X_i^H)= V_i$ for all $i\in [r]_0$ and $H\in \cH$;
	\item\label{item:set testers} $|W\cap \bigcap_{j\in [\ell]}\phi(Y_j)|= |W||Y_1|\cdots |Y_\ell|/n^\ell \pm \alpha n$ for all $(W,Y_1,\ldots,Y_\ell)\in \cW_{set}$;
	\item\label{item: vertex tester} $\omega(\bigcup_{H\in \cH}X_i^H\cap \phi^{-1}(v))=\omega(\bigcup_{H\in \cH}X_i^H)/n\pm \alpha n$ for all $(v,\omega)\in \cW_{ver}$.
\end{enumerate}
\end{theorem}

\subsection{Applications}

The multipartite framework can be used to obtain results for the non-partite setting.
The next theorem applies to graphs $G$ that are $(\eps,d)$-quasirandom;
that is, if $n$ is the order of $G$,
then $|N_G(u)|=(d\pm \eps)n$ and $|N_G(u)\cap N_G(v)|=(d^2\pm \eps)n$ for all distinct $u,v\in V(G)$.
Given~$G$ and a collection of graphs $\cH$ on at most $n$ vertices,
we say $(W,Y_1,\ldots,Y_k)$ is an \defn{$\ell$-set tester}
if $k\leq \ell$ and there exist distinct $H_1,\ldots,H_k\in \cH$ such that $W\In V(G)$ and $Y_i\In V(H_i)$ for all $i\in [k]$.
We say $(v,\omega)$ is an \defn{$\ell$-vertex tester} if $v\in V(G)$ and $\omega:\bigcup_{H\in \cH} V(H) \to [0,\ell]$.

\begin{theorem}\label{thm:quasirandom}
For all $\alpha>0$, there exist $\eps>0$ and $n_0\in\mathbb{N}$ such that the following holds for all $n\geq n_0$ and $d\geq \alpha$.
Suppose $G$ is an $(\eps,d)$-quasirandom graph on $n$ vertices and~$\cH$ is a collection of graphs on at most $n$ vertices with $|\cH|\leq \alpha^{-1}n$ and $\sum_{H\in \cH}e(H)\leq (1-\alpha)e(G)$ as well as $\Delta(H)\leq \alpha^{-1}$ for all $H\in \cH$.
Suppose $\cW_{set},\cW_{ver}$ are sets of  $\alpha^{-1}$-set testers and $\alpha^{-1}$-vertex testers of size at most~$n^{\log n}$, respectively.
Then there is a packing $\phi$ of $\cH$ into  $G$
such that 
\begin{itemize}
	\item $|W\cap \bigcap_{i\in [\ell]}\phi(Y_i)|= |W||Y_1|\cdots |Y_\ell|/n^\ell \pm \alpha n$ for all $(W,Y_1,\ldots,Y_\ell)\in \cW_{set}$;
	\item $\omega(\bigcup_{H\in \cH} V(H))\cap \phi^{-1}(v))=\omega(\bigcup_{H\in \cH} V(H))/n\pm \alpha n$ for all $(v,\omega)\in \cW_{ver}$.
\end{itemize}
\end{theorem}

In many scenarios when one applies approximate decomposition results, as for example Theorem~\ref{thm:quasirandom}, it is important that the graph $G-\phi(\cH)$ has `small' maximum degree.
Here, this can be easily achieved by {utilising} vertex testers $(v,\omega)$ where $\omega$ assigns to all $x\in\bigcup_{H\in \cH}V(H)$ the degree of~$x$.
{We remark that set and vertex testers in our main result are very flexible and capture many desirable properties.} For example, Theorem~\ref{thm:quasirandom} implies the approximate decomposition result in~\cite{ABCT:19} when restricted to graphs of bounded maximum degree.

In this paper, we give one example how to apply Theorem~\ref{thm:quasirandom}. 
By exploiting set testers, we can combine it with Keevash's results on hypergraph decompositions to decompose pseudorandom graphs into regular spanning graphs
as long as only a few graphs contain a few vertices in components of bounded size.
This is stronger as some results in~\cite{GJKKO:18} where a few graphs with almost all vertices in components of bounded size are required.\footnote{The results in~\cite{GJKKO:18} consider only $2$-regular graphs. However, their proof for the part where they consider collections of graphs $\cH$ that contain a few graphs with almost all vertices in components of bounded size carries over verbatim to $r$-regular graphs for any $r$ if $n$ is large in terms of $r$.}
For this result, we need a stronger notation of pseudorandomness as used by Keevash in~\cite{keevash:18b}.
We say a graph $G$ on $n$ vertices is \defn{$(\eps,s,d)$-typical} if $|\bigcap_{u\in U}N_G(u)|=(1\pm \eps)d^{|U|}n$ for all sets $U\In V(G)$ with $|U|\leq s$.

\begin{theorem}\label{thm:application}
	For all $\alpha>0$, there exist $\eps>0$ and $s,n_0\in \bN$ such that the following holds for all $n\geq n_0$ and $d\geq \alpha$.
	Suppose $G$ is a regular $(\eps,s,d)$-typical graph on $n$ vertices and $\cH$ is a collection of regular graphs on $n$ vertices with $\sum_{H\in \cH}e(H)=e(G)$ as well as $\Delta(H)\leq \alpha^{-1}$ for all $H\in\cH$.
	Suppose there are at least $\alpha n$ graphs $H\in \cH$ such that at least $\alpha n$ vertices in $H$ belong to components of order at most~$\alpha^{-1}$.
	Then there is a decomposition of the edge set of $G$ into $\cH$.
\end{theorem}

Theorem~\ref{thm:application} makes progress on a conjecture by Glock, Joos, Kim, K\"uhn and Osthus who conjecture in~\cite{GJKKO:18} that $K_n$ can be decomposed into any collection $\cH$ of regular bounded degree graphs with $\sum_{H\in \cH}e(H)=\binom{n}{2}$.

\section{Proof sketch}\label{sec:proof sketch}

Before we explain our approach,
we briefly sketch the approach of Kim, K\"uhn, Osthus and Tyomkyn in~\cite{KKOT:19}.
Their first step is to stack several graphs $H\in \cH$ together to a new graph $\widetilde{H}$
such that $\widetilde{H}[X_i^{\widetilde{H}},X_j^{\widetilde{H}}]$ is essentially regular for all $ij\in E(R)$.\footnote{In fact, their main theorem only applies to collections of graphs that are essentially regular and this stacking had to be performed again in~\cite{CKKO:19} and~\cite{JKKO:ta} which made the application in both cases technically involved.}
Let $\widetilde{\cH}$ be the collection of these graphs~$\widetilde{H}$.
They prove that such graphs $\widetilde{H}$ can be embedded into~$G$ by a probabilistic algorithm in a very uniform way.
For some $\gamma \ll \alpha$, they apply this algorithm to $\gamma n$ graphs in $\widetilde{\cH}$ in turn.
Observe that this may cause edge overlaps in $G$.
Nevertheless, after embedding $\gamma n$ graphs, they remove all `used' edges from $G$ and repeat.
At the end, they eliminate all edge overlaps by unembedding several vertices and complete the packing by {utilising} a thin edge slice put aside at the beginning.

Our approach is somewhat perpendicular to their approach.
We proceed cluster by cluster and find a function $\phi_i$ which maps almost all vertices in $\bigcup_{H\in \cH}X_i^H$ into $V_i$ 
and which is consistent with our partial packing so far.
{Our `Approximate Packing Lemma', stated in Section~\ref{sec:approx packing}, performs one such step using} an auxiliary hypergraph where we aim to find a large matching which is pseudorandom with respect to certain weight functions.
At the end, we complete the packing by also using a thin edge slice similar to~\cite{KKOT:19}.
At the beginning, we partition the clusters of our blow-up instance into many smaller clusters with the only purpose to ensure that $H[X_i^H,X_j^H]$ is a matching (see~Section~\ref{sec:partition}).
This preprocessing is comparably simple and first used in~\cite{RR:99}.

Both the approach in~\cite{KKOT:19} and ours draw on ideas from an alternative proof of the blow-up lemma by R\"odl and Ruci\'nski~\cite{RR:99}.
In spirit, our approach is again closer to the procedure in~\cite{RR:99} as they also embed the clusters of $H$ in turn.
Many generalisations of the original blow-up lemma build on this alternative proof.
We hope that our alternative proof of the blow-up lemma for approximate decompositions paths the way for further developments in the field.

Some ideas in this paper are taken from our recent article with Glock~\cite{EGJ:19b} on rainbow embeddings in graphs.
As a crucial part of our proof, we apply the main result from another paper with Glock~\cite{EGJ:19a}, where we prove the existence of quasirandom hypergraph matchings in hypergraphs with small codegree.
The idea that hypergraph matchings can be used to embed (almost) spanning graphs has been brought to our attention by~\cite{KKKO:18}.

\section{Preliminaries}
\subsection{Notation}
\COMMENT{Adapted from rainbow paper.}
For $k\in\mathbb{N}$, we write $[k]_0:=[k]\cup\{0\}=\{0,1,\ldots,k\}$, where $[0]=\emptyset$.
For a finite set $S$ and $k\in \mathbb{N}$,
we write $\binom{S}{k}$ for the set of all subsets of~$S$ of size~$k$
and~$2^S$ for the powerset of~$S$. 
For a set $\Set{i,j}$, we sometimes simply write~$ij$.
For $a,b,c\in \mathbb{R}$,
we write $a=b\pm c$ whenever $a\in [b-c,b+c]$.
For $a,b,c\in (0,1]$,
we sometimes write $a\ll b \ll c$ in our statements meaning that there are increasing functions $f,g:(0,1]\to (0,1]$
such that whenever $a\leq f(b)$ and $b \leq g(c)$,
then the subsequent result holds.
For $a\in (0,1]$ and $\mathbf{b}\in(0,1]^k$, we write $a\ll\mathbf{b}$ whenever $a\ll b_i$ for all $b_i\in\mathbf{b}$, $i\in[k]$.
For the sake of a clearer presentation, we avoid roundings whenever it does not affect the argument. 
%In particular, small parameters in our hierarchies are considered to be the inverse of a natural number.

For a graph $G$, let $V(G)$ and $E(G)$ denote the vertex set and edge set, respectively.
We say $u\in V(G)$ is a $G$-neighbour of $v\in V(G)$ if $uv\in E(G)$ and let $N_G(u)$ be the set of all $G$-neighbours of $u$
as well as $\dg_G(u):=|N_G(u)|$. As usual, $\Delta(G)$ denotes the maximum degree of $G$.
For $u,v\in V(G)$, let $N_G(u,v):=N_G(u)\cap N_G(v)$ denote the common neighbourhood of $u$ and~$v$, and let $N_G[u]:=N_G(u)\cup\{u\}$. 
For a set $S$, let $N_G(S):=\bigcup_{v\in S\cap V(G)}N_G(v)$.
{We frequently treat collections of graphs as the graph obtained by taking the disjoint union of all members; for example, for a collection $\cH$ of graphs, we define $N_\cH(v):=\bigcup_{H\in \cH}N_H(v)$ for all $v\in V(H)$ with $H\in\cH$.}
For disjoint subsets $A,B\In V(G)$, 
let $G[A,B]$ denote the bipartite subgraph of $G$ between $A$ and $B$
and $G[A]$ the subgraph of $G$ induced by $A$.
For convenience, let $G[A,A]:=G[A]$.
Let $e(G)$ denote the number of edges of $G$ and let $e_G(A,B)$ denote the number of edges of $G[A,B]$.
For $k\in\mathbb{N}$, let $G^k$ denote the $k$-th power of $G$, that is, the graph obtained from~$G$ by adding all edges between vertices whose distance in $G$ is at most $k$. 
{For graphs $G,H$, we write $G-H$ to denote the graph with vertex set $V(G)$ and edge set $E(G)\setminus E(H)$.}

\subsection{Probabilistic tools and graph regularity}
To verify the existence of subgraphs with certain properties we frequently consider random subgraphs and use McDiarmid's inequality to verify that specific random variables are highly concentrated around their mean.

\begin{theorem}[McDiarmid's inequality, see~\cite{mcdiarmid:89}\COMMENT{Lemma~1.2}] \label{thm:McDiarmid}
Suppose $X_1,\dots,X_m$ are independent Bernoulli random variables and suppose $b_1,\dots,b_m\in [0,B]$.
Suppose $X$ is a real-valued random variable determined by $X_1,\dots,X_m$ such that changing the outcome of $X_i$ changes $X$ by at most $b_i$ for all $i\in [m]$.
Then, for all $t>0$, we have $$\prob{|X-\expn{X}|\ge t} \le 2 \exp\left({-\frac{2t^2}{B\sum_{i=1}^m b_i}}\right).$$
\end{theorem}

We will also need the following two standard results concerning the robustness of $\eps$-regular graphs.
\begin{fact} \label{fact:regularity}
Suppose $G$ is an $(\eps,d)$-regular bipartite graph with vertex partition $(A,B)$ and $Y\In B$ with $|Y|\ge \eps|B|$. Then all but at most $2\eps|A|$ vertices of $A$ have $(d\pm \eps)|Y|$ neighbours in~$Y$.  
\end{fact}

%We will also often use the fact that super-regularity is robust with respect to small vertex and edge deletions.
\begin{fact}\label{fact:regularity robust}
Suppose $1/n \ll \eps \ll d$.
Suppose $G$ is an $(\eps,d)$-super-regular bipartite graph with vertex partition $(A,B)$, where $\eps^{1/6} n\le |A|,|B|\le n$. If $\Delta(H)\le \eps n$ and $X\In A\cup B$ with $|X|\le \eps n$, then $G[A\sm X,B\sm X]-H$ is $(\eps^{1/3},d)$-super-regular.\COMMENT{Proof: The number of edges between two vertex sets $Z_1$ and $Z_2$ is at least 
$(d-\eps)|Z_1||Z_2|- (\eps n + \eps n)(|Z_1|+|Z_2|)\geq |Z_1||Z_2|(d- \eps - 4\eps^{1/2})\geq |Z_1||Z_2|(d-\eps^{1/3})$ whenever $|Z_i|\geq \eps^{1/2}n= \eps^{1/3} \cdot \eps^{1/6}n$.}
\end{fact}

We will also use the next result from~\cite{DLR:95}.
(In~\cite{DLR:95} it is proved in the case when $|A|=|B|$ with $16\eps^{1/5}$ instead of $\eps^{1/6}$.
The version stated below can be easily derived from this.)

\begin{theorem}\label{thm: almost quasirandom}
Suppose $1/n\ll\eps \ll \gamma,d$.
Suppose $G$ is a bipartite graph with vertex partition~$(A,B)$ such that $|A|=n$, $\gamma n \leq |B|\leq \gamma^{-1}n$ and at least $(1-5\eps)n^2/2$ pairs $u,v\in A$ satisfy $\dg_G(u),\dg_G(v)\geq (d- \eps)|B|$ and $|N_G(u,v)|\leq (d+ \eps)^2|B|$.
Then $G$ is $\eps^{1/6}$-regular. 
%In particular, if $G$ is $(\eps,p)$-quasi-random, it is $(\eps^{1/6},p)$-super-regular.
\end{theorem}

At the end of our packing algorithm we apply the following version of the blow-up lemma due to Koml\'os, Sark\"ozy, and Szemer\'edi.

\begin{theorem}[Koml\'os, Sark\"ozy, and Szemer\'edi~\cite{KSS:97}]\label{thm:KSS}
Suppose $1/n\ll \eps \ll 1/\Delta,d$ and $1/n\ll 1/r$.
%Suppose $(H,G,R, (X_i)_{i\in [r]},(V_i)_{i\in [r]})$ is an $(\eps,d')$-super-regular and $\Delta$-bounded blow-up instance with $d'\geq d$,
Suppose $(H,G,R, (X_i)_{i\in [r]},(V_i)_{i\in [r]})$ is an $(\eps,d')$-super-regular and $\Delta$-bounded blow-up instance,
with $d'\geq d$ as well as $|V_i|=(1\pm\eps)n$ for all $i \in [r]$
and $(A_i)_{i\in[r]}$ is a collection of graphs
such that~$A_i$ is bipartite with vertex partition $(X_i,V_i)$ and $(\eps,d_i)$-super-regular for some $d_i\geq d$.
Then there is a packing $\phi$ of $H$ into $G$ such that $\phi(x)\in N_{A_i}(x)$ for all $x\in X_i$ and $i\in [r]$.
\end{theorem}

\subsection{Pseudorandom hypergraph matchings}
A key ingredient in the proof of our `Approximate Packing Lemma' in Section~\ref{sec:approx packing} is the main result from~\cite{EGJ:19a} on pseudorandom hypergraph matchings.

For this we need some more notation.
Given a set $X$ and an integer $\ell\in \bN$, 
an \defn{$\ell$-tuple weight function on~$X$} is a function $\omega\colon \binom{X}{\ell}\to \bR_{\ge 0}$.
For a subset $X'\In X$, we then define $\omega(X'):=\sum_{S\in \binom{X'}{\ell}}\omega(S)$. 
%Moreover, if $\cX\In \binom{X}{\ell}$, we write $\omega(\cX)$ for $\sum_{S\in \cX}\omega(S)$ as for usual weight functions. 
For $k\in[\ell]_0$ and a tuple $T\in \binom{X}{k}$, define 
$\omega(T):=\sum_{S\supseteq T}\omega(S)$ and let $\normV{\omega}{k}:=\max_{T\in \binom{X}{k}}\omega(T)$. 
Suppose $\cH$ is an $r$-uniform hypergraph and $\omega$ is an $\ell$-tuple weight function on~$E(\cH)$. Clearly, if $\cM$ is a matching, then a tuple of edges which do not form a matching will never contribute to~$\omega(\cM)$. We thus say that $\omega$ is \defn{clean} if $\omega(\cE)=0$ whenever $\cE\in\binom{E(\cH)}{\ell}$ is not a matching.

\begin{theorem}[\cite{EGJ:19a}]\label{thm:hypermatching}
Suppose $1/\Delta \ll \delta,1/r,1/L$, $r\geq 2$, and let $\eps:=\delta/50L^2r^2$.
Let $\cH$ be an $r$-uniform hypergraph with $\Delta(\cH)\leq \Delta$ and $\Delta^c(\cH)\le \Delta^{1-\delta}$ as well as $e(\cH)\leq \exp(\Delta^{\eps^2})$. 
Suppose that for each $\ell\in[L]$, we are given a set $\cW_\ell$ of clean $\ell$-tuple weight functions on $E(\cH)$ of size at most $\exp(\Delta^{\eps^2})$, such that $\omega(E(\cH))\ge \normV{\omega}{k}\Delta^{k+\delta}$ for all $\omega\in \cW_\ell$ and $k\in [\ell]$.
Then there exists a matching~$\cM$ in~$\cH$ such that $\omega(\cM)=(1\pm \Delta^{-\eps}) \omega(E(\cH))/\Delta^\ell$ for all $\ell\in [L]$ and $\omega \in \cW_\ell$.
\end{theorem}

\subsection{Refining partitions}\label{sec:partition}

Here, we provide a useful result to refine the vertex partition of a blow-up instance such that every $H\in\cH$ only induces a matching between its refined partition classes.
While in~\cite{RR:99} this procedure was easily obtained by applying the classical Hajnal--Szemer\'edi Theorem, 
we perform a random procedure to obtain more control on the mass distribution of a weight function with respect to the refined partition. 

\begin{lemma}\label{lem:split}
	Suppose $1/n\ll \beta \ll \alpha$ and $1/n\ll 1/r$.
	Suppose $\cH$ is a collection of at most $\alpha^{-1}n$ graphs, $(X_i^H)_{i\in[r]}$ is a vertex partition of $H$, and $\Delta(H)\leq \alpha^{-1}$ for every $H \in\cH$.
	Suppose $n/2\leq |X_i^H|=|X_i^{H'}| \leq 2n$ for all $H,H'\in \cH$ and $i\in [r]$.
	Suppose $\cW$ is a set of weight functions $\omega\colon \bigcup_{H\in \cH, i\in[r]} X_i^H\to [0,\alpha^{-1}]$ with $|\cW|\leq \eul^{\sqrt{n}}$.
	Then for all $H\in \cH$ and $i\in [r]$, there exists a partition $(X_{i,j}^H)_{j\in [\beta^{-1}]}$ of $X^H_i$ such that for all $H\in \cH,\omega\in \cW, i,i'\in [r],j,j'\in[\beta^{-1}]$ where $i\neq i'$ or $j\neq j'$, we have that
	\begin{enumerate}[label=\rm (\roman*)]
		\item\label{item:split1} $X^H_{i,j}$ is independent in $H^2$;% for all $H\in \cH,i,i'\in[r],j,j'\in[\beta^{-1}]$;
		\item\label{item:split2} $|X_{i,1}^H|\leq \ldots \leq |X_{i,\beta^{-1}}^H|\leq |X_{i,1}^H|+1$; %for all $H\in\cH$;
		\item\label{item:split4} $\omega(X_{i,j}^H)=\beta \omega(X_i^H) \pm \beta^{3/2}n$.
		\item\label{item:split3} $\sum_{H\in\cH}e_H(X_{i,j}^H,X_{i',j'}^H)= \beta^2\sum_{H\in \cH}e_H(X_{i}^H,X_{i'}^H) \pm n^{5/3}$;\COMMENT{Note that we defined $H[X_i,X_i]=H[X_i]$.}
	\end{enumerate}
\end{lemma}
\COMMENT{Note that we do not require $(X_i^H)_{i\in[r]}$ to be a vertex partition of $H$ into independent sets, which allows to apply the theorem for $r=1$.}
\proof
Our strategy is as follows.
We first consider every $H\in \cH$ in turn and construct a partition $(Y_{i,j})_{j\in[\beta^{-1}]}$ that essentially satisfies~\ref{item:split1}--\ref{item:split4} with $Y_{i,j}$ playing the role of $X_{i,j}^H$. 
Then we perform a vertex swapping procedure to resolve some conflicts in $Y_{i,j}$ and obtain $Z_{i,j}^H$. 
In the end, we permute the ordering of $(Z_{i,j}^H)_{j\in[\beta^{-1}]}$ for each $H\in \cH,i\in [r]$ to also ensure~\ref{item:split3}.

To simplify notation, we assume from now on that $|X_i^H|$ is divisible by $\beta^{-1}$ for all $i\in [r]$
and at the end we explain how very minor modifications yield the general case.

Fix some $H\in \cH$ and write $X_i$ for $X_i^H$.
We claim that there exist partitions $(Y_{i,j})_{j\in [\beta^{-1}]}$ of $X_i$ for each $i\in [r]$ such that for all $i\in [r],j\in [\beta^{-1}]$
\begin{enumerate}[label=(\alph*)]
	\item\label{item:correct split1} $|Y_{i,j}|=\beta |X_i|\pm\beta^2n$;

	\item\label{item:correct split2} at most $2\alpha^{-2}\beta^2 n$ pairs of vertices in $Y_{i,j}$ are adjacent in $H^2$;\COMMENT{There are at most $2n\alpha^{-2}/2$ many $P_3$'s in $H$ that start and end in $X_i$, and at most $n\alpha^{-1}$ edges in $X_i$. Each endpoint of such a $P_3$ or $P_2$ lies in $Y_{i,j}$ with probability $\beta$.}
	\item\label{item:correct split3} $\omega(Y_{i,j})=\beta \omega(X_i) \pm  \beta^2 n$ for all $\omega\in \cW$.
\end{enumerate}
Indeed, the existence of such partitions can be seen by assigning every vertex in $X_i$ uniformly at random to some $Y_{i,j}$ for $j\in [\beta^{-1}]$. 
Together with a union bound and Theorem~\ref{thm:McDiarmid} we conclude that~\ref{item:correct split1}--\ref{item:correct split3} hold simultaneously with positive probability.

Next, we slightly modify these partitions to obtain a new collection of partitions. 
These modifications can be performed for each $i\in [r]$ independently.
Hence, we fix some $i\in [r]$.
%We say a vertex $x'\in X_i$ is a $2$-neighbour of $x\in X_i$ if $x,x'$ have a common neighbour in $H$.
For all $j\in[\beta^{-1}]$, let $W_{j}\In Y_{i,j}$ be such that $|Y_{i,j} \setminus W_j|=\beta|X_i|- \beta^{5/3}n$ and $W_j$ contains all vertices in $Y_{i,j}$ that contain an $H^2$-neighbour in $Y_{i,j}$ (the sets $W_j$ clearly exist by~\ref{item:correct split1},\ref{item:correct split2}).
Let $\{w_1,\ldots,w_s\}=W:=\bigcup_{j\in [\beta^{-1}]}W_j$ and observe that $s=|X_i|-\sum_{[\beta^{-1}]}|Y_{i,j}\sm W_j|=\beta^{2/3}n$.

Now arbitrarily assign labels in $[\beta^{-1}]$ to the vertices in $W$ such that each label is used exactly~$\beta^{5/3}n$ times.
Let  $Z_j(0):=Y_{i,j}\setminus W_j$ for all $j\in[\beta^{-1}]$.
To obtain the desired partitions we perform the following swap procedure for every $t\in [s]$ in turn.
Say $w_t\in W_j$ and $w_t$ receives label $j'$.
We select $j''\in [\beta^{-1}]\setminus \{j,j'\}$ such that $w_t$ has no $H^2$-neighbour in $Z_{j''}(t-1)$ and such that $Z_{j''}(t-1)$ contains a vertex $w$ that has no $H^2$-neighbour in $Z_{j'}(t-1)$.
In such a case we say that $j''$ is \defn{selected} in step $t$.
Then we define $Z_{j''}(t):= (Z_{j''}(t-1)\cup \{w_t\})\setminus \{w\}$, $Z_{j'}(t):=Z_{j'}(t)\cup \{w\}$ and $Z_k(t):=Z_k(t-1)$ for all $k\in [\beta^{-1}]\setminus \{j',j''\}$.
Note that $\beta |X_i| - \beta^{5/3}n \leq |Z_k(t)|\leq \beta |X_i|$ for all $t\in [s],k\in [\beta^{-1}]$.
Observe also that we have always at least $\beta^{-1}/2$ choices to select $j''$ in step $t$.\COMMENT{There are at most $2\beta n\alpha^{-2}$ many $P_3$'s in $H$ that start in $Y_{i,j'}$, and at most $2\beta n\alpha^{-1}$ many $P_2$'s in $H$ that start in $Y_{i,j'}$. 
Hence, at most $6\alpha^{-2}$ entries $j''\in[\beta^{-1}]\sm\{j,j'\}$ are such that every vertex in $Y_{i,j''}$ has an $H^2$-neighbour in $Y_{i,j'}$, which implies that there are at least ${(3\beta^{-1})}/{4}$ good entries $j''\in[\beta^{-1}]\sm\{j,j'\}$ such that there exists a $w\in Y_{i,j''}$ which does not have an $H^2$-neighbour in $Y_{i,j'}$. 
Since $w_t$ has at most $\alpha^{-2}$ neighbours in $H^2$, at most $\alpha^{-2}$ of these good entries $j''\in[\beta^{-1}]\sm\{j,j'\}$ are such that $w_t$ has an $H^2$-neighbour in $Y_{i,j''}$. Hence, there are at least $\beta^{-1}/2$ choices to select $j''$ in step $t$.}
As $s= \beta^{2/3}n$, we can ensure that each $j''\in [\beta^{-1}]$ is selected, say, at most $10\beta^{5/3}n$ times.
We write $Z_{i,j}^H:=Z_j(s)$ and then it is straightforward to see that for all $j\in [\beta^{-1}]$ we have that
\begin{enumerate}[label=(\alph*$'$)]
	\item\label{item:correct splitI1} $|Z_{i,j}^H|= \beta|X_i^H|$;\COMMENT{Each label $j\in[\beta^{-1}]$ is used exactly $\beta^{5/3}n$ times. In each such step $t$ where $j$ is used, the cardinality of the current set $Z_{j}(t-1)$ is increased by exactly one and remains the same during the other steps. Since $|Z_j(0)|=|Y_{i,j}\sm W_j|=\beta|X_i|-\beta^{5/3}n$, this implies that $|Z_{i,j}^H|=\beta|X_i|$.}
	\item\label{item:correct splitI2} $Z_{i,j}^H$ is independent in $H^2$;\COMMENT{This follows immediately from the construction of the swaps because $W$ contains all vertices that contain a $H^2$-neighbour and each $w\in W$ is swapped to a set such that $w$ does not have a $H^2$-neighbour in that set.}
	\item\label{item:correct splitI3} $\omega(Z_{i,j}^H)=\beta \omega(X_i^H) \pm  \beta^{3/2} n/2$ for all $\omega\in \cW$.\COMMENT{Since $|W_j|=\beta^{5/3}n$, we have that $\omega(Z_j(0))=\omega(Y_{i,j})\pm\alpha^{-1}\beta^{5/3}n$. 
Since we perform at most $10\beta^{5/3}n$ swaps in $Y_{i,j}$ and we add exactly $\beta^{5/3}n$ vertices to $Y_{i,j}$, we have that $\omega(Z_{i,j})=\omega(Z_j(0))\pm11\alpha^{-1}\beta^{5/3}n$ and thus,
$$\omega(Z_{i,j})=\omega(Y_{i,j})\pm12\alpha^{-1}\beta^{5/3}n
\stackrel{\ref{item:correct split3}}{=}\beta \omega(X_i)\pm \beta^{3/2}n/2.
$$
}
\end{enumerate}
As $i\in [r]$ and $H\in \cH$ are chosen arbitrarily, the statements~\ref{item:correct splitI1}--\ref{item:correct splitI3} hold for all $i\in [r],H\in \cH$.

In the remainder of the proof we show how to find permutations $\{\pi_{i}^H\}_{H\in \cH,i\in[r]}$ 
such that~\ref{item:split3} holds for  $X_{i,j}^H:=Z_{i,\pi_i^H(j)}^H$.
The answer is simple, random permutations yield~\ref{item:split3} with probability, say, at least $1/2$.
To see this, fix $i,i'\in [r],j,j'\in [\beta^{-1}]$ where $i\neq i'$ or $j\neq j'$, and then Theorem~\ref{thm:McDiarmid} implies that $\sum_{H\in\cH}e_H(X_{i,j}^H,X_{i',j'}^H)= \beta^2\sum_{H\in \cH}e_H(X_{i}^H,X_{i'}^H) \pm n^{5/3}/2$ with probability at least $1-1/n$.\COMMENT{
...
}
A~union bound over all $i,i',j,j'$ completes the proof.

In the beginning we made the assumption that $\beta^{-1}$ divides $X_i^H$.
To avoid this assumption, we simply remove a set $\widetilde{X}_i^H$ of size at most $\beta^{-1}-1$ from $X_i^H$ such that $\beta^{-1}$ divides $X_i^H\setminus \widetilde{X}_i^H$ and perform the entire procedure with $X_i^H\setminus \widetilde{X}_i^H$ instead of $X_i^H$.
That is, for all $H\in \cH$ and $i\in [r]$, we obtain a partition $(\tX_{i,j}^H)_{j\in [\beta^{-1}]}$ of $X^H_i\sm \tX^H_i$ satisfying~\ref{item:split1}--\ref{item:split3}, where $|\tX_{i,j}^H|=|\tX_{i,j'}^H|$ for all $i\in[r],j,j'\in[\beta^{-1}]$, and~\ref{item:split3} and \ref{item:split4} hold with error terms `$\pm n^{5/3}/2$' and `$\pm\beta^{3/2}n/2$', respectively.
At the very end we add the vertices in $\widetilde{X}_i^H$ to the partition $(\tX_{i,j}^H)_{j\in [\beta^{-1}]}$ while preserving~\ref{item:split1} and~\ref{item:split2}.
We may do so by performing a swap argument as before.
Observe that our error bounds give us enough room to spare.\COMMENT{Fix $H\in\cH, i\in[r]$, and assign labels in $[\beta^{-1}]$ to $\tX^H_i$ such that each label is used at most once.
Let $\tX_j(0):=\tX^H_{i,j}$ and we proceed for $t\in|\tX_i^H|$ in turn.
Assume $x_t\in\tX_i^H$ receives label $j\in[\beta^{-1}]$.
Select $j'\in[\beta^{-1}]\sm\{j\}$ such that $x_t$ has no $H^2$-neighbour in $\tX_{j'}(t-1)$ and $\tX_{j'}(t-1)$ contains a vertex $x$ such that $x$ has no $H^2$-neighbour in $X_{i,j}$. As before, we have always $\beta^{-1}/2$ choices to select $j'$ in step $t$. 
Set $\tX_j(t):=\tX_j(t-1)\cup\{x\}$ and $\tX_{j'}(t):=(\tX_{j'}(t-1)\cup\{x_t\})\sm \{x\}$. 
Since each label is chosen at most once,~\ref{item:split1} holds. 
By the construction of the swap,~\ref{item:split2} holds. \ref{item:split3} and \ref{item:split4} hold since our error terms `$\pm n^{5/3}/2$' and `$\pm\beta^{3/2}n/2$' had still enough room to spare and $|\tX_i^H|\leq\beta^{-1}-1$.
}
\endproof

\section{Approximate packings}\label{sec:approx packing}
The goal of this section is to provide an `Approximate Packing Lemma' (Lemma~\ref{lem:packing lemma}).
Given a blow-up instance $(\cH,G,R,\cX,\cV)$, it allows us to embed almost all vertices of $\bigcup_{H\in \cH}X_i^H$ into $V_i$, while maintaining crucial properties for future embedding rounds of other clusters.
To describe this setup we define a \emph{packing instance} and collect some more notation.

\subsection{Packing instances}
Given a graph $G$ and a set $\cE$, we call $\psi\colon E(G)\to2^\cE$ an \emph{edge set labelling} of $G$. 
A label $\alpha\in\cE$ \emph{appears} on an edge $e$ if $\alpha\in\psi(e)$. 
We define the \emph{maximum degree $\Delta_\psi(G)$ of $\psi$} as the maximum number of edges of $G$ on which any fixed label appears.
We define the \emph{maximum codegree $\Delta^c_\psi(G)$ of $\psi$} as the maximum number of edges of $G$ on which any two fixed labels appear together.

\medskip
Let $r\in \bN_0$.
We say $(\cH,G,R,\cA,\psi)$ is a \defn{packing-instance of size $r$} if
\begin{itemize}
%[label={\rm (P\arabic*)}]
	\item\label{item:P1} $\cH$ is a collection of graphs, and $G$ and $R$ are graphs, where $V(R)=[r]_0$;

%where $(X^H_i)_{i\in[r]_0}$ is a partition of $V(H)$ into independent sets for every $H\in\cH$;
%\item $G$ and $R$ are graphs
	\item $\cA=\bigcup_{H\in\cH,i\in[r]_0}A^H_i$ is a union of balanced bipartite graphs $A^H_i$ with vertex partition~$(X^H_i,V_i)$;
	\item $(X^H_i)_{i\in[r]_0}$ is a partition of $H$ into independent sets for every $H\in\cH$, and $(V_i)_{i\in[r]_0}$ is a partition of $V(G)$;
\item $R=R_A\cup R_B$ is the union of two edge-disjoint graphs with $N_R(0)=[r]$;
	\item for all $H\in\cH$, the graph $H[X^H_i,X^H_j]$ is a matching if $ij\in E(R)$ and empty otherwise;
	\item\label{item:P5} $\psi \colon E(\cA)\to 2^\cE$ is an edge set labelling such that $\Delta_\psi(A^H_i)\leq 1$ for all $H\in\cH,i\in[r]_0$, and for every label $\alpha\in\cE$ with $\alpha\in\psi(xv)\cap\psi(x'v')$ and $v,v'\in V(G)$, we have $v=v'$.\COMMENT{This encodes the condition that the edge set containing a particular label is a star.
This property is important when establishing~\ref{emb lem:boundedness}, cf. Claim~\ref{claim1}.}  
\end{itemize}
In such a case, we write for simplicity $\sX_i:=\bigcup_{H\in\cH}X^H_i$ and $\sA_i:=\bigcup_{H\in\cH}A^H_i$ for each $i\in[r]_0$, and whenever we write $xv\in E(\sA_i)$, we tacitly assume that $x\in\sX_i, v\in V_i$.
{The only reason why $R$ is the disjoint union of two graphs lies in the nature of our approach;
while $R_A$ represents parts of~$R$ as in the statement of our main result (Lemma~\ref{lem:main matching}, which is very similar to Theorem~\ref{thm:main_new}),
the edges of~$R_B$ represent copies of edge slices of $G$ that in the end will be used to complete the approximate packing.
We use copies here to obtain a unified setup for the Approximate Packing Lemma, alternatively, we could have used parallel edges in the reduced graph.
}

{The aim of this section is to map almost all vertices of $\sX_0$ into $V_0$ by defining a function $\sigma\colon \sX_0^\sigma\to V_0$ in~$\sA_0$ 
(that is, $x\sigma(x)\in E(\sA_0)$)}
where $\sX_0^\sigma\subseteq \sX_0$.
{(Hence, we refer to subgraphs of $\sA_i$ as candidacy graphs.)}
For convenience, we identify such a function~$\sigma$ with its \emph{corresponding edge set}~$M$ defined as $M=M(\sigma):=\{xv\colon x\in \sX_0^\sigma, v\in V_0, \sigma(x)=v \}$.
We say 
\begin{align}\label{eq:conflict free}
\begin{minipage}[c]{0.8\textwidth}\em
$\sigma\colon \sX_0^\sigma\to V_0$ is a \emph{conflict-free packing} if $\sigma|_{\sX_0^\sigma\cap X_0^H}$ is injective for all $H\in\cH$ and $\psi(e)\cap\psi(f)=\emptyset$ for all distinct $e,f\in M(\sigma)$.
\end{minipage}\ignorespacesafterend
\end{align}
The set $\psi(xv)$ will encode the set of edges of $G$ that are used for the embedding when mapping~$x$ to~$v$.
The property that $\psi(e)\cap\psi(f)=\emptyset$ for all distinct $e,f\in M(\sigma)$ will guarantee that in the proof of our main result (Lemma~\ref{lem:main matching}) every edge in~$G$ is used at most once.

Given a conflict-free packing $\sigma\colon \sX_0^\sigma\to V_0$ in~$\sA_0$, we update the remaining candidacy graphs and their edge set labelling according to the following two definitions.

\begin{defin}[Updated candidacy graphs]\label{def:updated candidacy}
	For a conflict-free packing $\sigma\colon \sX_0^\sigma\to V_0$ in~$\sA_0$ and all $H\in\cH,i\in[r]$, let $A^H_i[\sigma]$ be the \defn{updated candidacy graph (with respect to~$\sigma$)} which is defined by the spanning subgraph of $A^H_i$ that contains precisely those edges $xv\in E(A^H_i)$ for which the following holds:
	if $x$ has an $H$-neighbour $x_0\in \sX_0^\sigma$ (which would be unique), then $\sigma(x_0)v\in E(G[V_0,V_i])$.
\end{defin}

\begin{defin}[Updated labelling]\label{def:updated labelling}
	For a conflict-free packing $\sigma\colon \sX_0^\sigma\to V_0$ in $\sA_0$, 
	let $\psi[\sigma]$ be the \defn{updated edge set labelling  (with respect to $\sigma$)} defined as follows: 
	for all $H\in\cH, i\in[r]$ and $xv\in E(A^H_i[\sigma])$, if $x$ has an $H$-neighbour $x_0\in \sX_0^\sigma$, then set $\psi[\sigma](xv):=\psi(xv)\cup \{\sigma(x_0),v\}$, and otherwise set $\psi[\sigma](xv):=\psi(xv)$.\COMMENT{Note that $\{\sigma(x_0),v\}\in E(G)$ because $xv\in E(A^H_i[\sigma])$.}
\end{defin}

{In order to be able to analyse our packing process in Section~\ref{sec:mainpart}, we carefully maintain quasirandom properties of the candidacy graphs throughout the procedure.
To this end,
we refer to a packing instance $(\cH,G, R,\cA,\psi)$ of size $r$ as an \emph{$(\eps,\bd)$-packing-instance}, where $\bd=(d_A,d_B,d_0,\ldots,d_r)$,
if} 
\begin{enumerate}[label={\rm (P\arabic*)}]
	\item\label{item:P G superreg} $G[V_i,V_j]$ is $(\eps,d_Z)$-super-regular for all $ij\in E(R_Z), Z\in\{A,B\}$;
	\item\label{item:P A superreg} $A^H_i$ is $(\eps,d_i)$-super-regular for all $H\in\cH,i\in[r]_0$;
	\item\label{item:P triple intersection} $e_{H}(N_{A^H_i}(v_i),N_{A^H_j}(v_j))=
(d_id_j\pm\eps)e_{H}(X^H_i,X^H_j)$ for all $H\in\cH, ij\in E(R_A), v_iv_j\in E(G[V_i,V_j])$;\COMMENT{$R_A$ will be the reduced graph corresponding to edges between candidacy graphs $A^H_i$, whereas $R_B$ corresponds to edges between candidacy graphs $B^H_i$ used for the absorption; (for $B^H_i$ we do not care about the triple intersection).}
	\item\label{item:P boundedness} $\Delta_\psi(\sA_i)\leq (1+\eps)d_i|V_i|$ for all $i\in[r]_0$.
\end{enumerate}
{Property~\ref{item:P boundedness} ensures that no edge is a potential candidate for too many graphs in $\cH$ and~\ref{item:P triple intersection} enables us to maintain this property for future embedding rounds (see Lemma~\ref{lem:packing lemma}\ref{emb lem:boundedness} below).}
Let $\sP=(\cH,G, R,\cA,\psi)$ be an $(\eps,\bd)$-packing-instance of size~$r$.
Similarly as for a blow-up instance, we say $(W,Y_1,\ldots,Y_k)$ is an \emph{$\ell$-set tester} for $\sP$ if $k\leq \ell$ and there exist distinct $H_1,\ldots,H_k\in\cH$ such that $W\subseteq V_0$ and $Y_j\subseteq X_0^{H_j}$ for all $j\in[k]$.
For $i\in N_{R_A}[0]$ and $v\in V_i$, we say $\omega\colon E(\sA_i)\to[0,\ell]$ is an \emph{$\ell$-edge tester with centre $v$} for $\sP$ if $\omega(x'v')=0$ for all $x'v'\in E(\sA_i)$ with $v'\in V_i, v'\neq v$.\COMMENT{That is, the positive weight of $\omega$ induces a star in $\sA_i$.}
{We say $\omega\colon E(\sA_0)\to[0,\ell]$ is an \emph{$\ell$-edge tester with centres in $\sX_0$} if there exist vertices $\{x_H\}_{H\in\cH}$ with $x_H\in X_0^H$ for each $H\in\cH$ such that $\omega(x'v')=0$ for all $x'v'\in E(\sA_0)$ with $x'\notin \{x_H\}_{H\in\cH}$.}
Further, let $\dim(\omega)$ be the \emph{dimension of $\omega$} defined as
\begin{align}\label{eq:def dim(w)}
\dim(\omega)=\begin{cases}
1 &\text{ if $\omega(E(\sA_i))=\omega(E(A_i^H))$ for some $H\in\cH$,}
\\ 2 & \text{ otherwise.}
\end{cases}
\end{align}
\COMMENT{$\dim(\omega)=1$ if the star lives only on one $H\in\cH$.}

Moreover, for every $H\in\cH$, let $H_+$ be an auxiliary supergraph of $H$ that is obtained by adding a maximal number of edges between $X^H_0$ and $X^H_i$ for every $i\in[r]$ subject to $H_+[X^H_0,X^H_i]$ being a matching.
We call $\cH_+:=\bigcup_{H\in\cH}H_+$ an \emph{enlarged graph of~$\cH$}.
We say that $\sP$ is \emph{nice (with respect to $\cH_+$)} if 
\begin{enumerate}[label={\rm (N\arabic*)}]
	\item\label{A good} $|N_{A^{H}_i}(x_i)\cap N_{G}(v_j)|=(d_id_Z\pm\eps)|V_i|$ for all  $x_jv_j\in E(A^H_j)$ whenever $\{x_i\}=N_{H_+}(x_j)\cap X^H_i$, $H\in\cH,ij\in E(R_Z), Z\in\{A,B\}$;
	\item\label{G good} $|N_G(v_i,v_j)\cap V_0|=(d_A^2\pm\eps)|V_0|$ for all $ij\in E(R_A-\{0\})$ and $v_iv_j\in E(G[V_i,V_j])$.
\end{enumerate}
Using standard regularity methods (see Facts~\ref{fact:regularity} and~\ref{fact:regularity robust}), it is straightforward to verify the following:
\begin{align}\label{eq:nice emb instance}
\begin{minipage}[c]{0.85\textwidth}\em
For every $(\eps,\bd)$-packing-instance $(\cH,G, R,\cA,\psi)$ of size $r$ and every enlarged graph $\cH_+$ of $\cH$, there exist spanning subgraphs $G'\subseteq G$ and $\cA'\subseteq\cA$
such that $(\cH,G', R,\cA',\psi)$ is a nice $(\eps',\bd)$-packing-instance of size $r$ with respect to $\cH_+$ for some $\eps'$ with $\eps\ll\eps'\ll 1/r$. 
\end{minipage}\ignorespacesafterend
\end{align}
\COMMENT{Let $\eps\ll\heps$. For every $H\in\cH$, let $A_0^{H,bad}$ consist of all edges $x_0v_0\in E(A_0^H)$ for which there is some $j\in N_{R_Z}(0), Z\in\{A,B\}$, $\{x_j\}=N_{H_+}(x_0)\cap X_j^H$ and
\begin{align}\label{eq:bad A_0}
|N_{A_j^H}(x_j)\cap N_{G}(v_0)|\neq (d_jd_Z\pm 3\eps)|V_j|.
\end{align}
For all $H\in\cH,j\in N_{R_Z}(0), Z\in\{A,B\}$, let $A_j^{H,bad}$ consist of all edges $x_jv_j\in E(A_j^H)$ for which there exists $\{x_0\}=N_{H_+}(x_j)\cap X_0^H$ and
\begin{align}\label{eq:bad A_j}
|N_{A_0^H}(x_0)\cap N_{G}(v_j)|\neq (d_0d_Z\pm 3\eps)|V_0|.
\end{align}
For all $jk\in E(R_A[[r]])$, let $G^{bad}$ be the spanning subgraph of $G$ such that an edge $v_jv_k\in E(G[V_j,V_k])$ belongs to $G^{bad}[V_j,V_k]$ if
\begin{align}\label{eq:bad G}
|N_G(v_j,v_k)\cap V_0|\neq(d_A^2\pm\eps)|V_0|.
\end{align}
Using Fact~\ref{fact:regularity}, it is easy to see that $\Delta(A_0^{H,bad})\leq 2r\eps |V_0|$, $\Delta(A_j^{H,bad})\leq 2\eps|V_j|$, and $\Delta(G^{bad})\leq 2r\eps n$.
Now, for all $H\in\cH,j\in[r]_0$, let
\begin{align*}
G':=G-E(G^{bad}), && A^{H,good}_j:=A^H_j-E(A^{H,bad}_j), && \sA^{good}_j:=\textstyle\bigcup_{H\in\cH}A^{H,good}_j.
\end{align*}
Using Fact~\ref{fact:regularity robust}, we conclude that $G'[V_k,V_\ell]$ is $(\eps^{1/4},d_A)$-super-regular for all $k\ell\in E(R_A[[r]])$, and $A^{H,good}_j$ is $(\eps^{1/4},d_j)$-super-regular for all $H\in\cH,j\in[r]_0$.
Crucially, we now have the following properties:
\begin{align}
\begin{split}\label{eq:good A_0}
&|N_{A^{H,good}_j}(x_j)\cap N_{G'}(v_0)|=(d_jd_Z\pm\heps)|V_j|,
\\& \text{for all  $x_0v_0\in E(A^{H,good}_0)$ whenever $H\in\cH,j\in N_{R_Z}(0), Z\in\{A,B\}$ and $\{x_j\}=N_{H_+}(x_0)\cap X_j^H$;} 
\end{split}
\\[1em]
\begin{split}\label{eq:good A_j}
&|N_{A^{H,good}_0}(x_0)\cap N_{G'}(v_j)|=(d_0d_Z\pm\heps)|V_0|,
\\& \text{for all  $x_jv_j\in E(A^{H,good}_0),H\in\cH,j\in N_{R_Z}(0), Z\in\{A,B\}$, whenever $\{x_0\}=N_{H_+}(x_j)\cap X_0^H$}; 
\end{split}
\\[1em]
\begin{split}\label{eq:good G}
&|N_{G'}(v_j,v_k)\cap V_0|=(d_A^2\pm\heps)|V_0|
\\&\text{for all  $v_jv_k\in E(G'[V_j,V,k])$ and $jk\in E(R_A[[r]])$.}
\end{split}
\\[1em]
\begin{split}\label{eq:good triple intersection}
&e_H(N_{A_j^{H,good}}(v_j),N_{A_k^{H,good}}(v_k))=(d_jd_k\pm\heps)e_{H}(X_j^H,X_k^H)
\\&\text{for all $H\in\cH, jk\in E(R_A)$, and $v_jv_k\in E(G'[V_j,V_k])$.}
\end{split}
\end{align}
Indeed, when replacing $A^H_j$ and $A^H_k$ with $A_j^{H,good}$ and $A_k^{H,good}$ in~\eqref{eq:good triple intersection}, we may have removed some edges that contributed to $e_H(N_{A_j^{H}}(v_j),N_{A_k^{H}}(v_k))$ but at most $\heps^2n$ edges for every edge $v_jv_k\in E(G'[V_j,V_k])$ because $\Delta(A_\ell^{H,bad})\leq 2r\eps |V_\ell| n$ for all $\ell\in[r]_0$, which implies~\eqref{eq:good triple intersection}.
Similar arguments hold for~\eqref{eq:good A_0}, \eqref{eq:good A_j}, and~\eqref{eq:good G}.  
}

\subsection{Approximate Packing Lemma}
We now state our Approximate Packing Lemma.
{Roughly speaking it states that given a packing instance, 
we can find a conflict-free packing such that the updated candidacy graphs are still super-regular, albeit with a smaller density.
Moreover, with respect to certain weight functions on the candidacy graphs,
the updated candidacy graphs behave as we would expect this by a random and independent deletion of the edges.
}

\begin{lemma}[Approximate Packing Lemma]\label{lem:packing lemma}
Let $1/n\ll\eps\ll\eps'\ll \bd,1/r,1/s$.
Suppose $(\cH,G,R,\cA,\psi)$ is an $(\eps,\bd)$-packing-instance of size $r$, $\|\psi\|\leq s$, {$|\cH|\leq sn$}, $|V_i|=(1\pm\eps)n$ for all $i\in[r]_0$, $\sum_{H\in\cH}e_{H}(X_0^H,X_i^H)\leq d_An^2$ for all $i\in N_{R_A}(0)$, and $e_{H}(X^H_i,X^H_j)\geq{\eps'^2n}$ for all $H\in\cH, ij\in E(R)$.
Suppose $\Delta_\psi^c(\sA_i)\leq \sqrt{n}$ for all $i\in N_{R_A}[0]$, and 
suppose $\cW_{set},\cW_{edge}$ are sets  of $s$-set testers and $s$-edge testers of size at most $n^{3\log n}$, respectively.

Then there is a conflict-free packing $\sigma\colon \sX_0^\sigma\to V_0$ in $\sA_0$ such that for all $H\in\cH$, we have $|\sX_0^\sigma\cap X_0^H|\geq (1-\eps')n$ and for all $i\in[r]$ there exists a spanning subgraph $A^{H,new}_i$ of the updated candidacy graph $A^H_i[\sigma]$ (where $\sA_i^{new}:=\bigcup_{H\in\cH}A^{H,new}_i$) with
\begin{enumerate}[label={\rm (\Roman*)\textsubscript{\ref{lem:packing lemma}}}]
\item\label{emb lem:1} $A^{H,new}_i$ is $(\eps',d_id_Z)$-super-regular for all $i\in N_{R_Z}(0), Z\in\{A,B\}$;
\item\label{emb lem:triple intersection}
$e_{H}(N_{A_i^{H,new}}(v_i),N_{A_j^{H,new}}(v_j))=
(d_id_j d_A^2\pm\eps')e_{H}(X^H_i,X^H_j)$ for all $ij\in E(R_A-\Set{0})$ and $v_iv_j\in E(G[V_i,V_j])$;

\item\label{emb lem:edge weight}
$\omega(E(\sA_i^{new}))=(1\pm\eps'^2)d_A\omega(E(\sA_i))\pm\eps'^2n^{\dim(\omega)}$ for all $\omega\in\cW_{edge}$ and $i\in N_{R_A}(0)$;

\item\label{emb lem:boundedness} $\Delta_{\psi[\sigma]}(\sA_i^{new})\leq (1+\eps')d_id_A|V_i|$  for all $i\in N_{R_A}(0)$; 
\item\label{emb lem:codegree} $\Delta_{\psi[\sigma]}^c(\sA_i^{new})\leq \sqrt{n}$ for all $i\in N_{R_A}(0)$;

\item\label{emb lem:set testers} $|W\cap \bigcap_{j\in [\ell]}\sigma(Y_j\cap \sX_0^\sigma)|= |W||Y_1|\cdots |Y_\ell|/n^\ell \pm \eps' n$ for all $(W,Y_1,\ldots,Y_\ell)\in \cW_{set}$;

\item\label{emb lem:vertex testers} 
$\omega(M(\sigma))=(1\pm\eps'){\omega(E(\sA_0))}/{d_0n}\pm \eps'n$ for all $\omega\in\cW_{edge}$ with centre in $V_0$ {or centres in~$\sX_0$}.

\end{enumerate}
\end{lemma}
\COMMENT{Deleted:
\\ (VI): $|\sigma^{-1}(v)|\geq (1-\eps')|\cH|$ and $|\sigma^{-1}(v)\cap \sX_0^\ast|\leq \eps'|\cH|$ for all $v\in V_0$; (follows by~\ref{emb lem:vertex testers}).}

Properties~\ref{emb lem:1}, \ref{emb lem:triple intersection} and \ref{emb lem:boundedness} ensure that \ref{item:P A superreg}--\ref{item:P boundedness} are also satisfied for the updated candidacy graphs $\sA_i^{new}$, respectively, and \ref{emb lem:codegree} ensures that the codegree of the updated labelling $\psi[\sigma]$ is still small on $\sA_i^{new}$. 
Property~\ref{emb lem:edge weight} states that the weight of the edge testers on the updated candidacy graphs $\sA_i^{new}$ is what we would expect by a random sparsification of the edges in~$\sA_i$, and~\ref{emb lem:set testers} and~\ref{emb lem:vertex testers} guarantee that $\sigma$ behaves like a random packing with respect to the set and edge testers.

We split the proof into two steps.
In Step~\ref{step:hypergraph}, we construct an auxiliary hypergraph and apply Theorem~\ref{thm:hypermatching} to obtain the required conflict-free packing $\sigma$.
By defining suitable weight functions in Step~\ref{step:weight functions}, we employ the conclusions of Theorem~\ref{thm:hypermatching} to establish~\ref{emb lem:1}--\ref{emb lem:vertex testers}.

\begin{proof}
Let $\cH_+$ be an enlarged graph of $\cH$, and for $i\in[r]$, let $\sA_i^{bad}$ and $\sA_i^{good}$ be spanning subgraphs of $\sA_i$ such that $\sA_i^{bad}$ contains precisely those edges $xv\in E(\sA_i)$ where $N_{\cH_+}(x)\cap \sX_0=\emptyset$, and $E(\sA_i^{good}):=E(\sA_i)\sm E(\sA_i^{bad})$.
{We may assume that $|\cH|=sn$ and $\sum_{H\in\cH}e_{H}(X_0^H,X_i^H)\leq (d_A+\eps'^{3/2})n^2$ for all $i\in N_{R_A}(0)$, where the last inequality will be only used in~\eqref{eq:bound N}. (Otherwise we artificially add some graphs to $\cH$ subject to the condition that still $e_{H}(X^H_i,X^H_j)\geq{\eps'^2n}$ for all $H\in\cH, ij\in E(R)$, and accordingly we add some graphs to $\cA$ satisfying~\ref{item:P G superreg}--\ref{item:P boundedness}.)}
We may also assume that
$\psi\colon E(\cA)\to2^\cE$ is such that $|\psi(e)|=s$ for all $e\in E(\sA_0)$ (otherwise we add artificial labels that we delete at the end again),\COMMENT{We only need that $|\psi(e)|=s$ in order to obtain a $(s+2)$-uniform hypergraph in Step~\ref{step:hypergraph}.} and $(\cH,G,R,\cA,\psi)$ is a nice $(\eps,\bd)$-packing-instance with respect to $\cH_+$ (otherwise we may employ~\eqref{eq:nice emb instance} and replace $\eps$ by some $\tilde{\eps}$, where $\eps\ll\tilde{\eps}\ll{\eps'}$; observe also that this does not cause problems with the weight of the edge testers in~\ref{emb lem:edge weight} and~\ref{emb lem:vertex testers}, as the operation in~\eqref{eq:nice emb instance} only deletes few edges of $\cA$ incident to every vertex).\COMMENT{Note that when we employ~\eqref{eq:nice emb instance}, and obtain $\sA_0'\subseteq \sA_0$ we also have to argue that for each $\omega\in\cW_{edge}$ we do not delete too much $\omega$-weight on $\sA_0$. 
For $v\in V_0$, we remove at most $2r\eps|V_0|$ edges of $\sA_0$ for each $H\in\cH$. 
That is, if $\dim(\omega)=1$ (weight is a star only on one $H\in\cH$) we remove at most $s2r\eps|V_0|\leq\eps^{1/2}n$ weight incident to $v$, and if $\dim(\omega)=2$ (weight appears on more than one  $H\in\cH$) we remove at most $s2r\eps|V_0||\cH|\leq\eps^{1/2}n^2$ weight incident to $v$. 
Recall that we defined $\cW_{edge}$ only as weight functions on a star, that is, if for $\omega\in\cW_{edge}$ we have $\omega(xv)>0$ then $\omega(x'v')=0$ for all $v'\neq v$.
Hence, $\omega(E(\sA_0))=\omega(E(\sA_0'))\pm\eps^{1/2}n^2$ and~\ref{emb lem:vertex testers} still holds.
\\\ref{emb lem:edge weight} also holds because the error term is `$\pm\eps'^2n$' for $\dim(\omega)=1$ and `$\pm\eps'^2n^2$' for $\dim(\omega)=2$.}

\begin{step}\label{step:hypergraph}
Constructing an auxiliary hypergraph
\end{step}

We want to use Theorem~\ref{thm:hypermatching} to find the required conflict-free packing $\sigma$ in $\sA_0$.
To this end,  let~$(V_0^H)_{H\in\cH}$ be disjoint copies of $V_0$, and for $H\in\cH$ and $e=x_0v_0\in E(A_0^{H})$, let $e^H:=x_0v_0^H$ where $v_0^H$ is the copy of $v_0$ in $V_0^H$.
Let $f_e:=e^H\cup \psi(e)$ for each $e\in E(A_0^H), H\in\cH$ and let $\cH^{aux}$ be the $(s+2)$-uniform hypergraph with vertex set $\bigcup_{H\in\cH}(X_0^H\cup V_0^H)\cup \cE$
and edge set $\{f_e\colon e\in E(\sA_0) \}$.
A~key property of the construction of $\cH^{aux}$ is a bijection between conflict-free packings $\sigma$ in $\sA_0$ and matchings $\cM$ in $\cH^{aux}$ by assigning $\sigma$ to $\cM=\{f_e\colon e\in M(\sigma) \}$. (Recall that $M=M(\sigma)$ is the edge set corresponding to $\sigma$.\COMMENT{$M :=\{xv\colon x\in \sX_0^\sigma, v\in V_0, \sigma(x)=v \}$})
%Indeed, for a matching $\cM$ in $\cH$ and $f_e,f_{e'}\in \cM$, 
%we have $f_e\cap f_{e'}=\emptyset$ and thus, $\psi(e)\cap \psi(e')=\emptyset$.
%Conversely, any conflict-free packing $\sigma\colon \sX_0^\sigma\to V_0$ in $\sA_0$ gives rise to a matching in $\cH$ by setting $\cM$ as above.
%%$\cM=\{f_e\colon e\in \{x_0v_0\colon x_0\in \sX_0^\sigma, v_0=\sigma(x_0) \} \}$.

It is easy to estimate $\Delta(\cH^{aux})$ and $\Delta^c(\cH^{aux})$ in order to apply Theorem~\ref{thm:hypermatching}. 
Since $A_0^{H}$ is $(\eps,d_0)$-super-regular for each $H\in\cH$, $|X_0^H|=|V_0|=(1\pm\eps)n$, and $\Delta_\psi(\sA_0)\leq (1+\eps)d_0|V_0|$, we conclude~that
\begin{align}\label{eq:Delta(H)}
\Delta(\cH^{aux})\leq(d_0+3\eps)n=:\Delta.
\end{align}

Note that the codegree in $\cH^{aux}$ of two vertices in $\bigcup_{H\in\cH}(X_0^H\cup V_0^H)$ is at most~1, and similarly, the codegree in $\cH^{aux}$ of a vertex in $\bigcup_{H\in\cH}(X_0^H\cup V_0^H)$ and a label in $\cE$ is at most~1 because $\Delta_\psi(A_0^H)\leq 1$ for all $H\in\cH$.
By assumption, $\Delta_\psi^c(\sA_0)\leq \sqrt{n}$.
Altogether, this implies that
\begin{align}\label{eq:Delta^c(H)}
\Delta^c(\cH^{aux})\leq \sqrt{n}\leq \Delta^{1-\eps^2}.
%\frac{n}{m^{3t+6}\log^3n}\leq \frac{(d_0+\eps)n}{m^{3t+6}\log^2n}.
\end{align}

Suppose $\cW=\bigcup_{\ell\in[s]}\cW_\ell$ is a set of size at most $n^{4\log n}$ of given weight functions $\omega\in\cW_\ell$ for $\ell\in[s]$ with $\omega\colon \binom{E(\sA_0)}{\ell}\to[0,s]$.
Note that every weight function $\omega\colon\binom{E(\sA_0)}{\ell}\to[0,s]$ naturally corresponds to a weight function $\omega_{\cH^{aux}}\colon\binom{E(\cH^{aux})}{\ell}\to[0,s]$ by defining $\omega_{\cH^{aux}}(\{f_{e_1},\ldots,f_{e_\ell}\}):=\omega(\{e_1,\ldots,e_\ell\})$.
We will explicitly specify~$\cW$ in Step~\ref{step:weight functions} and it is simple to check that for each $\omega\in\cW$ the corresponding weight function $\omega_{\cH^{aux}}$ will be clean.
Our main idea is to find a hypergraph matching in $\cH^{aux}$ that behaves like a typical random matching with respect to $\{\omega_{\cH^{aux}}\colon \omega\in\cW \}$ in order to establish~\ref{emb lem:1}--\ref{emb lem:vertex testers}.

Suppose $\ell\in[s]$ and $\omega\in\cW_\ell$. 
If $\omega(E(\sA_0))\geq n^{1+\eps/2}$ or $\ell\geq 2$,\COMMENT{We need the $\max$-bound in~\eqref{eq:w bound codegree} only for $\ell=1$.} define $\tilde{\omega}:=\omega$.
Otherwise, choose $\tilde{\omega}\colon E(\sA_0)\to[0,s]$ such that $\omega\leq\tilde{\omega}$ and $\tilde{\omega}(E(\sA_0))= n^{1+\eps/2}$.
By~\eqref{eq:Delta(H)} and~\eqref{eq:Delta^c(H)}, we can apply Theorem~\ref{thm:hypermatching} 
(with $(d_0+3\eps)n,\eps^2,s+2,s,\{\tilde{\omega}_{\cH^{aux}}\colon{\omega}\in\cW_\ell \}$ playing the roles of $\Delta,\delta,r,L,\cW_\ell$) to obtain a matching $\cM$ in $\cH^{aux}$ that corresponds to a conflict-free packing $\sigma\colon \sX_0^\sigma\to V_0$ in $\sA_0$ with its corresponding edge set $M=M(\sigma)$ that satisfies the following properties (where $\heps:=\eps^{1/2}$):
\begin{align}
\label{eq:w packing size}
\omega(M)&=(1\pm\heps)\frac{\omega(E(\sA_0))}{(d_0n)^\ell} \text{ for } \omega\in\cW_\ell, \ell\in[s]
\text{ where } \omega(E(\sA_0))\geq \normV{\omega}{k}\Delta^{k+\eps^2} \text{ for all } k\in[\ell];\\
\label{eq:w bound codegree}
\omega(M)&\leq\max\Big\{(1+\heps)\frac{\omega(E(\sA_0))}{d_0n}, n^\eps\Big\} \text{ for all }\omega\in \cW_1.
\end{align}
\COMMENT{$$(1\pm\Delta^{-{\eps^2}/{50s^{2}(\Delta+2)^2 }})\frac{\omega(E(\sA_0))}{((d_0+3\eps)n)^\ell}=
(1\pm\Delta^{-{\eps^2}/{50s^{2}(\Delta+2)^2 }})\frac{\omega(E(\sA_0))}{(1\pm\eps^{2/3})(d_0n)^\ell}
=(1\pm\eps^{1/2})\frac{\omega(E(\sA_0))}{(d_0n)^\ell}
$$}
\COMMENT{For~\eqref{eq:w packing size} and $\omega\in\cW_1$, notice that $\tilde{\omega}(E(\sA_0))\geq n^{1+\eps/2}\geq s((d_0+3\eps)n)^{1+\eps^2}\geq \normV{\tilde{\omega}}{1}\Delta^{1+\eps^2}$.
Note further that $\Delta^{-{\eps^2}/{50s^{2}(t+2)^2 }}\leq \eps$.
\\For~\eqref{eq:w bound codegree} and $\omega\in\cW_1$ with $\omega(E(\sA_0))<n^{1+\eps/2}$, we have $\omega(M)\leq\tilde{\omega}(M)\leq (1+\heps)\frac{n^{1+\eps/2}}{d_0n}\leq n^\eps$.}

One way to exploit~\eqref{eq:w packing size} is to control the number of edges in $M$ between sufficiently large sets of vertices.
To this end, for subsets $S\subseteq X_0^H$ and $T\subseteq V_0$ for some $H\in\cH$ with $|S|,|T|\geq 2\eps n$, we define a weight function $\omega_{S,T}\colon E(A_0^{H})\to\{0,1\}$ with
\begin{align}\label{eq:def w_S,T}
\omega_{S,T}(e):=\mathbbm{1}_{\{e\in E(A_0^{H}[S,T])\}}.
\end{align}
That is, $\omega_{S,T}(M)$ counts the number of edges in $A_0^{H}$ between $S$ and $T$ that lie in $M$.
Since $A_0^H$ is $(\eps,d_0)$-super-regular we have that $e_{A_0^{H}}(S,T)=(d_0\pm\eps)|S||T|\geq\eps^3n^2$ which implies together with~\eqref{eq:w packing size} that whenever $\omega_{S,T}\in\cW$, then $\sigma$ is chosen such that 
\begin{align}\label{eq:size w_S,T}
|\sigma(S\cap \sX_0^\sigma)\cap T|=\omega_{S,T}(M)=(1\pm 2\heps)\frac{|S||T|}{n}.
\end{align}

\begin{step}
\label{step:weight functions}
Employing weight functions to conclude~\ref{emb lem:1}--\ref{emb lem:vertex testers}
\end{step}

By Step~\ref{step:hypergraph}, we may assume that~\eqref{eq:w packing size} and~\eqref{eq:w bound codegree} hold for a set of weight functions~$\cW$ that we will define during this step.
We will show that for this choice of~$\cW$ the conflict-free packing $\sigma\colon \sX_0^\sigma\to V_0$ as obtained in Step~\ref{step:hypergraph} satisfies~\ref{emb lem:1}--\ref{emb lem:vertex testers}.
Similarly as in Definition~\ref{def:updated candidacy} (here, $\cH$ is replaced by~$\cH_+$), we define subgraphs $A^{H,\ast}_i$ of $A^{H}_i$ as follows.
\begin{align}\label{eq:def A_i^H,*}
\begin{minipage}[c]{0.85\textwidth}\em
For all $H\in\cH$, $i\in[r]$, let $A^{H,\ast}_i$ be the spanning subgraph of $A^{H}_i$ containing precisely those edges $xv\in E(A^{H}_i)$ for which the following holds: if $\{x_0\}=N_{\cH_+}(x)\cap \sX_0^\sigma$, then $\sigma(x_0)v\in E(G[V_0,V_i])$.
%Let $\sA_i^\ast:=\bigcup_{H\in\cH}A^{H,\ast}_i$.
\end{minipage}\ignorespacesafterend
\end{align}

\COMMENT{Note that $N_{A^{H,\ast}_j}(x)=N_{A^{H}_j}(x)$ if $x$ has no $\cH_+$-neighbour in $\sX_0^\sigma$}
Observe that $A^{H,\ast}_i$ is a spanning subgraph of the updated candidacy graph $A^H_i[\sigma]$ as in Definition~\ref{def:updated candidacy}.
By taking a suitable subgraph of $A^{H,\ast}_i$, we will in the end obtain the required candidacy graph $A^{H,new}_i$.

First, we show that $|\sX_0^\sigma\cap X_0^H|\geq (1-3\heps)n$ for each $H\in\cH$. 
Adding $\omega_{X_0^H,V_0}$ as defined in~\eqref{eq:def w_S,T} for every $H\in\cH$ to $\cW$ and using~\eqref{eq:size w_S,T} yields
\begin{align}\label{eq:packing size}
|\sX_0^\sigma\cap X_0^H|=\omega_{X_0^H,V_0}(M)\geq (1-3\heps)n.
\end{align}
\medskip

\begin{substep}
Checking~\ref{emb lem:1}
\end{substep}

For all $H\in\cH$ and $i\in N_{R_Z}(0), Z\in\{A,B\}$ we proceed as follows.
Let $Y^{H}_i:=N_{H_+}(\sX_0^\sigma)\cap X^H_i$.
We first show that $A^{H,\ast}_i[Y^{H}_i,V_i]$ is $(\heps^{1/18},d_id_Z)$-super-regular (see~\eqref{eq:tilde A superregular}). 
We do so by showing that every vertex in $Y^{H}_i\cup V_i$ has the appropriate degree, and that the common neighbourhood of most pairs of vertices in~$V_i$ have the correct size such that we can employ Theorem~\ref{thm: almost quasirandom} to guarantee the super-regularity of $A^{H,\ast}_i[Y^{H}_i,V_i]$.

Note that $|Y^{H}_i|\geq |\sX_0^\sigma\cap X_0^H|-2\eps n\geq (1-4\heps)n$ by~\eqref{eq:packing size}. 
For every vertex $x\in Y^{H}_i$ with $\{x_0\}=N_{H_+}(x)\cap \sX_0^\sigma$, we have $\dg_{A^{H,\ast}_i}(x)=|N_{A^{H}_i}(x)\cap N_{G}(\sigma(x_0))|$. 
Since our packing-instance is nice, \ref{A good} implies that 
%\begin{align*}
$\dg_{A^{H,\ast}_i}(x)=(d_id_Z\pm\eps)|V_i|$.
%\end{align*}
For $v\in V_i$, let $N_v:=N_{A^{H}_i}(v)$. Observe that
\begin{align}\label{eq:degree v A^i_j def}
\dg_{A^{H,\ast}_i[Y^{H}_i,V_i]}(v)=|\sigma(N_{H_+}(N_v)\cap \sX_0^\sigma)\cap N_G(v)|,
\end{align}
and $|N_{H_+}(N_v)\cap X_0^H|=|N_v|\pm2\eps n=(d_i\pm5\eps)n$, and $|N_G(v)\cap V_0|=(d_Z\pm3\eps)n$.
Adding for every vertex $v\in V_i$, the weight function $\omega_{S,T}$ as defined in~\eqref{eq:def w_S,T} for $S:=N_{H_+}(N_v)\cap X_0^H$ and $T:=N_G(v)\cap V_0$ to $\cW$, we obtain by~\eqref{eq:size w_S,T} and~\eqref{eq:degree v A^i_j def} that
\begin{align}\label{eq:degree v A^i_j}
\dg_{A^{H,\ast}_i[Y^{H}_i,V_i]}(v)
=
(1\pm2\heps)|N_{H_+}(N_v)\cap X_0^H||N_G(v)\cap V_0|n^{-1}
=(d_id_Z\pm{\heps}^{1/2})|Y^{H}_i|.
\end{align}

Next, we use Theorem~\ref{thm: almost quasirandom} to show that $A^{H,\ast}_i[Y^{H}_i,V_i]$ is $(\heps^{1/18},d_id_Z)$-super-regular. 
We call a pair of vertices $u,v\in V_i$ \emph{good} if $|N_{A^{H}_i}(u,v)|=(d_i\pm\eps)^2|X^H_i|$, and $|N_G(u,v)\cap V_0|=(d_Z\pm\eps)^2|V_0|$.
By the $\eps$-regularity of $A^{H}_i$ and $G[V_0,V_i]$, there are at most $2\eps |V_i|^2$ pairs $u,v\in V_i$ which are not good.

For all good pairs $u,v\in V_i$, let $S_{u,v}:=N_{H_+}(N_{A^{H}_i}(u,v))\cap X_0^H$ and $T_{u,v}:=N_G(u,v)\cap V_0$.
We add the weight function $\omega_{S_{u,v},T_{u,v}}$ as defined in~\eqref{eq:def w_S,T} to $\cW$. 
Observe that $|S_{u,v}|=|N_{A^{H}_i}(u,v)|\pm2\eps n=(d_i\pm \eps^{1/2})^2n$\COMMENT{$|N_{A^{H}_i}(u,v)|\pm2\eps n=(d_i\pm\eps)^2(1\pm\eps)n\pm2\eps n=(d_i\pm \eps^{1/2})^2n$.} and $|T_{u,v}|=(d_Z\pm\eps^{1/2})^2n.$
By~\eqref{eq:size w_S,T}, we obtain for all good pairs $u,v\in V_j$ that
\begin{align}
|N_{A^{H,\ast}_i[Y^{H}_i,V_i]}(u,v)|
&=|\sigma(S_{u,v}\cap \sX_0^\sigma)\cap T_{u,v}|=(1\pm2\heps)|S_{u,v}||T_{u,v}|n^{-1}\nonumber\\
&\leq (d_id_Z+\heps^{1/3})^2|Y^{H}_i|.
\label{eq:size good pairs u,v}
\end{align}
\COMMENT{$n\leq (1+5\heps)|Y^{H}_i|$}Now, by~\eqref{eq:degree v A^i_j} and~\eqref{eq:size good pairs u,v}, we can apply Theorem~\ref{thm: almost quasirandom}, and obtain that
\begin{align}\label{eq:tilde A superregular}
\text{$A^{H,\ast}_i[Y^{H}_i,V_i]$ is $(\heps^{1/18},d_id_Z)$-super-regular.}
\end{align}

In order to complete the proof of~\ref{emb lem:1}, we show that we can find a spanning subgraph $A^{H,new}_i$ of~$A^{H,\ast}_i$ that is $(\eps',d_id_Z)$-super-regular.
Let 
\begin{align}\label{eq:def A^new}
E(A^{H,new}_i[Y^{H}_i,V_i]):=E(A^{H,\ast}_i[Y^{H}_i,V_i]).
\end{align}
For every vertex $x\in X^H_i\sm Y^{H}_i$, we have that $\dg_{A^{H,*}_i}(x)=(d_i\pm\eps)|V_i|$ because $A^{H}_i$ is $(\eps,d_i)$-super-regular. 
Suppose $\cW^{bad}$ is a collection of at most $n^{4\log n}$ weight functions $\omega^{bad}\colon E(\sA_i^{bad})\to[0,s]$;\COMMENT{Recall that $\sA_i^{bad}$ is the spanning subgraph of $\sA_i$ containing precisely those edges $xv\in E(\sA_i)$ where $N_{\cH_+}(x)\cap \sX_0=\emptyset$.} we will specify $\cW^{bad}$ explicitly when we establish~\ref{emb lem:edge weight}.
We claim that we can delete for every vertex $x\in X^H_i\sm Y^{H}_i$ some incident edges in $A^{H,*}_i$ and obtain a subgraph $A_i^{H,new}$ such that 
\begin{align}
\label{eq:degree x A^new}
\dg_{A^{H,new}_i}(x)&=(d_id_Z\pm2\eps)|V_i| \text{ for every $x\in X^H_i\sm Y^{H}_i$;}
\\ \label{eq:bound random sparsification}
\omega^{bad}(E(\sA_i^{new}))&=(1\pm\eps)d_Z\omega^{bad}(E(\sA_i^{bad}))\pm\eps n
\text{ for every $\omega^{bad}\in\cW^{bad}$.}
\end{align}
This can be easily seen by a probabilistic argument: 
For all $H\in\cH$ and $x\in X^H_i\sm Y^{H}_i$, we keep each edge incident to $x$ in $A^{H}_i$ independently at random with probability~$d_Z$. 
Then, McDiarmid's inequality (Theorem~\ref{thm:McDiarmid}) together with a union bound yields that~\eqref{eq:degree x A^new} and~\eqref{eq:bound random sparsification} hold simultaneously with probability at least, say, $1/2$.\COMMENT{
\eqref{eq:degree x A^new} is obvious.
For~\eqref{eq:bound random sparsification}, let $(X_e)_{e\in E(\sA_i^{bad})}$ be independent Bernoulli random variables with $p=d_Z$ indicating whether $e$ is still present in $\sA_i^{new}$ and let $X:=\sum_{e\in E(\sA_i^{bad})}\omega^{bad}(e)X_e$. 
Clearly, $\expn{X}=d_Z\omega^{bad}(E(\sA_i^{bad}))$.
Changing the outcome of $X_e$ changes $X$ by at most $\omega^{bad}(e)=:b_e$. 
Hence, $\sum_{e\in E(\sA_i^{bad})}b_e=\omega^{bad}(E(\sA_i^{bad}))$. 
Thus, McDiarmid's inequality implies that
\begin{align*}
\prob{|X-\expn{X}|\geq \eps d_Z\omega^{bad}(E(\sA_i^{bad}))+\eps n}
&\leq 2\exp\left(-\frac{2\left(\eps d_Z\omega^{bad}(E(\sA_i^{bad}))+\eps n\right)^2}{s\omega^{bad}(E(\sA_i^{bad}))}\right)\leq \eul^{-\sqrt{n}}.
\end{align*}
A union bound implies that the claimed statements~\eqref{eq:degree x A^new} and~\eqref{eq:bound random sparsification} hold with probability at least $1/2$.
} 

Since $|X^H_i|=(1\pm\eps)n$, we have that $|X^H_i\sm Y^{H}_i|\leq 4\heps n$ by~\eqref{eq:packing size}.
Hence,~\eqref{eq:tilde A superregular} implies together with~\eqref{eq:degree x A^new} that $A^{H,new}_i$ is $(\eps',d_id_Z)$-super-regular, which establishes~\ref{emb lem:1}.

\begin{substep}
Checking~\ref{emb lem:triple intersection}
\end{substep}

For all $H\in\cH,ij\in E(R_A-\{0\})$, and $v_iv_j\in E(G[V_i,V_j])$ we proceed as follows.
Let $\tE:=E(H[N_{A^H_i}(v_i),N_{A^H_j}(v_j)])$ and
\begin{align*}
S&:=\{\{{x_i'},{x_j'} \}\subseteq X_0^H \colon x_ix_i',x_jx_j'\in E(H_+), x_ix_j\in\tE \},
\\S_1&:=\{S'\in S\colon |S'|=1  \}, \quad \text{and}\quad S_2:=\{S'\in S\colon |S'|=2  \},
\\E_1&:=\{xv\in E(A_0^H)\colon x\in S_1,v\in N_G(v_i,v_j) \},
%\\S_2&:=\{\{x_{x_j},x_{x_k}\}\in S_{v_j,v_k}\colon x_{x_j}\neq x_{x_k}  \},
\\E_2&:=\{\{xv,x'v' \}\in\textstyle\binom{E(A_0^H)}{2}\colon 
\{x,x' \}\in S_2, v\in N_G(v_i), v'\in N_G(v_j), v\neq v'
\}.
\end{align*}
By assumption (see~\ref{item:P triple intersection}), we have that $|\tE|=(d_id_j\pm\eps)e_{H}(X^H_i,X^H_j)$.
Since $e_{H}(X^H_i,X^H_j)\geq \eps'^2n$, we conclude that
\begin{align}\label{eq:size S}
|S|
=|S_1|+|S_2|
=(d_id_j\pm\eps)e_{H}(X^H_i,X^H_j)\pm 4\eps n
=(d_id_j\pm\heps)e_{H}(X^H_i,X^H_j).
\end{align}
Note that the term of `$\pm 4\eps n$' in~\eqref{eq:size S} accounts for possible vertices $x_i\in N_{A_i^H}(v_i)$ and $x_j\in N_{A^H_j}(v_j)$ that do not have an $H_+$-neighbour in $X_0^H$.\COMMENT{Recall that $|X_\ell^H|=(1\pm\eps)n$ for all $\ell\in[r]_0$.}

We define the following weight functions $\omega_1\colon E(A_0^H)\to\{0,1\}$ and $\omega_2\colon \binom{E(A_0^H)}{2}\to\{0,1\}$ by setting $\omega_1(e):=\mathbbm{1}_{\{e\in E_1 \}}$ and $\omega_2(\{e_1,e_2\}):=\mathbbm{1}_{\{\{e_1,e_2\}\in E_2 \}}$ and add them to $\cW$.\COMMENT{Note that $e_1\cap e_2=\emptyset$ because $v'\in N_G(v_j)\sm\{v\}$ and thus $\omega_2$ corresponds to a clean weight function.}
By the definition of $A_i^{H,new}$ (recall~\eqref{eq:def A_i^H,*} and~\eqref{eq:def A^new}), we crucially observe that
\begin{align}\label{eq:weight triple intersection}
e_{H}(N_{A_i^{H,new}}(v_i),N_{A_j^{H,new}}(v_j))=\omega_1(M)+\omega_2(M)\pm 5\heps n.
\end{align}
Note that the term of `$\pm 5\heps n$' in~\eqref{eq:weight triple intersection} accounts for possible vertices $x_i\in N_{A_i^H}(v_i)$ and $x_j\in N_{A^H_j}(v_j)$ that do not have an $H_+$-neighbour in $X_0^H$ (at most $4\eps n$), and possible vertices in $S$ that are left unembedded (at most $4\heps n$ by~\eqref{eq:packing size}).

Let us for the moment assume that $|S_1|,|S_2|\geq{ \eps'^5 n}$ (otherwise the claimed estimations in~\eqref{eq:size omega =} and~\eqref{eq:size omega neq} below are trivially true). 
Since $A_0^H$ is $(\eps,d_0)$-super-regular and $|N_G(v_i,v_j)\cap V_0|=(d_A^2\pm3\eps)n$ by~\ref{G good}, we obtain that
\begin{align}
\label{eq:weight E_1}
\omega_1(E(A_0^H))=|E_1|=(d_0\pm\eps)|S_1||N_G(v_i,v_j)\cap V_0|=(d_0d_A^2\pm \heps)|S_1|n.
\end{align}
By Fact~\ref{fact:regularity}, all but at most $6\eps n$ elements $\{x_i',x_j'\}\in S_2$ are such that $x_k'$ has $(d_0\pm\eps)|N_G(v_k)\cap V_0|$ neighbours in $N_G(v_k)\cap V_0$ for both $k\in\{i,j\}$ because $A_0^H$ is $(\eps,d_0)$-super-regular and  $G[V_0,V_k]$ is $(\eps,d_A)$-super-regular. 
Each of these $6\eps n$ exceptional elements contributes at most $|N_G(v_i)\cap V_0||N_G(v_j)\cap V_0|\leq 3n^2$ to $\omega_2(E(A_0^H))$.
This implies that
\begin{align}
\omega_2(E(A_0^H))=|E_2|&=(d_0\pm2\eps)^2(|S_2|\pm6\eps n)|N_G(v_i)\cap V_0||N_G(v_j)\cap V_0|\pm18\eps n^3
\nonumber\\&=(d_0^2d_A^2\pm \heps)|S_2|n^2.
\label{eq:weight E^neq}
\end{align}
For all $e_1\in E(A_0^H)$, the number of edges $e_2$ for which $\{e_1,e_2\}\in E_2$ is at most $2n$,\COMMENT{For fixed $e_1=xv$, also $x'\in e_2=x'v'$ is already fixed. At most $2n$ choices for $v'$.} implying
$\normV{\omega_2}{1}\Delta^{1+\eps^2}\leq 2n\Delta^{1+\eps^2}\leq\omega_2(E(A_0^H))$,
and clearly, $\normV{\omega_2}{2}\Delta^{2+\eps^2}\leq \Delta^{2+\eps^2}\leq \omega_2(E(A_0^H))$ and $\normV{\omega_1}{1}\Delta^{1+\eps^2}\leq \Delta^{1+\eps^2}\leq\omega_1(E(A_0^H))$. (Recall that $\Delta=(d_0+3\eps)n$.)
Hence, by~\eqref{eq:w packing size}, we conclude that
\begin{align}
\label{eq:size omega =}
\omega_1(M)
&=(1\pm\heps)\frac{\omega_1(E(A_0^H))}{d_0n}
\stackrel{\eqref{eq:weight E_1}}{=}
d_A^2|S_1|\pm{\eps'^2} e_{H}(X^H_i,X^H_j),
\\\label{eq:size omega neq}
\omega_2(M)
&=(1\pm\heps)\frac{\omega_2(E(A_0^H))}{(d_0n)^2}
\stackrel{\eqref{eq:weight E^neq}}{=}
d_A^2|S_2|\pm{\eps'^2} e_{H}(X^H_i,X^H_j).
\end{align}
\COMMENT{For~\eqref{eq:size omega =} (\eqref{eq:size omega neq} follows similarly):
\begin{align*}
(1\pm\heps)\frac{(d_0d_A^2\pm\heps)|S_1|n}{d_0n}
=(1\pm\heps)\left(1\pm\frac{\heps}{d_A^2d_0}\right)d_A^2|S_1|
=d_A^2|S_1|\pm \heps^{2/3}|S_1|
=d_A^2|S_1|\pm \heps^{1/2}e_H(X_i^H,X_j^H).
\end{align*}
}Clearly, the final equalities in~\eqref{eq:size omega =} and~\eqref{eq:size omega neq} are also true if {$|S_1|,|S_2|<\eps'^5 n$ because $e_{H}(X^H_i,X^H_j)\geq \eps'^2n$}.
Now, together with~\eqref{eq:size S} and~\eqref{eq:weight triple intersection} this implies that
\begin{align*}
e_{H}(N_{A_i^{H,new}}(v_i),N_{A_j^{H,new}}(v_j))
=d_A^2|S|\pm 3{\eps'^2} e_{H}(X^H_i,X^H_j)
=
(d_id_jd_A^2\pm \eps')e_{H}(X^H_i,X^H_j),
\end{align*}
which establishes~\ref{emb lem:triple intersection}.

\begin{substep}\label{step:III}
Checking~\ref{emb lem:edge weight}
\end{substep}

We will even show that~\ref{emb lem:edge weight} holds for all $\omega\in\cW_{edge}\cup\cW_{edge}'$ with $\omega\colon E(\sA_i)\to[0,s]$ and centre $v\in V_i, i\in N_{R_A}(0)$, where $\cW_{edge}'$ is a set of edge testers that we will explicitly specify in Step~\ref{step:IV} when establishing~\ref{emb lem:boundedness}.
For all $\omega\in\cW_{edge}\cup\cW_{edge}'$ with centre $v\in V_i$, $i\in N_{R_A}(0)$ we define a weight function $\omega_0\colon E(\sA_0)\to [0,s]$ by
\begin{align*}
\omega_0(x_0v_0):=
\begin{cases}
\omega(x_iv) &\text{if $ \{x_i\}= N_{\cH_+}(x_0)\cap \sX_i$, $x_iv\in E(\sA_i^{good})$ and $v_0v\in E(G)$,}
\\ 0 &\text{otherwise,}
\end{cases}
\end{align*}
\COMMENT{We could also write $x_iv\in E(\sA_i)$ instead of $x_iv\in E(\sA_i^{good})$. However, I think it is more intuitive to afterwards understand~\eqref{eq:bound updated weight linear} as we split $\omega(E(\sA_i))$ into the contributions in $\sA_i^{bad}$ and $\sA_i^{good}$.}and we add $\omega_0$ to $\cW$. 
(Recall that $\sA_i^{good}$ is the spanning subgraph of $\sA_i$ containing precisely those edges $x_iv_i\in E(\sA_i)$, where $N_{\cH_+}(x_i)\cap \sX_0\neq\emptyset$.)

For every edge $x_iv\in E(\sA_i^{good})$ with $\{x_0\}= N_{\cH_+}(x_i)\cap\sX_0$, property~\ref{A good} yields that
\begin{align*}
|N_{\sA_0}(x_0)\cap N_G(v)|=(d_0d_A\pm3\eps)n.
\end{align*}
Hence, every edge $x_iv\in E(\sA_i^{good})$ contributes weight $\omega(x_iv)\cdot(d_0d_A\pm3\eps) n$ to $\omega_0(E(\sA_0))$, and we obtain
\begin{align*}
\omega_0(E(\sA_0))&= \omega(E(\sA_i^{good})) (d_0d_A\pm3\eps)n.
\end{align*}
By the definition of $\sA_i^{new}$ (recall~\eqref{eq:def A_i^H,*} and~\eqref{eq:def A^new}), if\COMMENT{Note that we can not simply write `if and only if' as $x_iv$ might also be present in $E(\sA_i^{new})$ if $x_iv\in\Gamma$ as defined below.} $\sigma(x_0)\in N_G(v)$ for $\{x_0\}= N_{\cH_+}(x_i)\cap\sX_0$, then the edge $x_iv\in E(\sA_i^{good})$ is in $E(\sA_i^{new})$. 
Hence, if $x_0v_0\in M(\sigma)=M$, then this contributes weight $\omega_0(x_0v_0)$ to $\omega(E(\sA_i^{new}))$.
If $\omega(E(\sA_i^{good}))\geq\eps n$, then $\omega_0(E(\sA_0))\geq n^{1+\eps}\geq s\Delta^{1+\eps^2}\geq \normV{\omega_0}{1}\Delta^{1+\eps^2}$, and thus~\eqref{eq:w packing size} implies that
\begin{align}\label{eq:weight omega_0}
\omega_0(M)
=(1\pm\heps)\frac{\omega_0(E(\sA_0))}{d_0n}
=(1\pm2\heps)d_A\omega(E(\sA_i^{good}))\pm\heps n.
\end{align}
If $\omega(E(\sA_i^{good}))<\eps n$, then~\eqref{eq:w bound codegree} implies that
\begin{align*}
\omega_0(M)\leq \max\Big\{(1+\heps)\frac{\omega_0(E(\sA_0))}{d_0n},n^\eps \Big\}\leq \heps n,
\end{align*}
and hence,~\eqref{eq:weight omega_0} also holds in this case.

We now make a key observation:
\begin{align}\label{eq:weight A^new edges}
\omega(E(\sA_i^{new}))=\omega_0(M)+\omega(\Lambda)\pm\omega(\Gamma),
\end{align}
for $\Gamma:=\{x_iv\in E(\sA_i^{good})\colon N_{\cH_+}(x_i)\cap\sX_0^\sigma=\emptyset \}$\COMMENT{Observe that for $x_iv_i\in \Gamma$, we have $N_{\cH_+}(x_i)\cap\sX_0\neq\emptyset$. Note further that due to the random sparsification not all edges in $\Gamma$ survive in $\sA_i^{new}$, which is the reason for `$\pm$' in~\eqref{eq:weight A^new edges}.} and $\Lambda:=E(\sA_i^{bad})\cap E(\sA_i^{new})$.\COMMENT{That is, $\Lambda$ contains the edges of $\sA_i^{bad}$ on which the label $e$ appears and which are still present in $\sA_i^{new}$ after the random sparsification from which we obtain $\sA_i^{new}$.}
Next, we want to control $\omega(\Gamma)$ and $\omega(\Lambda)$.

In order to bound $\omega(\Gamma)$, we define a weight function $\omega_\Gamma\colon E(\sA_0)\to[0,s]$ by
\begin{align*}
\omega_\Gamma(x_0v_0):=
\begin{cases}
\omega(x_iv) &\text{if $\{x_i\}=N_{\cH_+}(x_0)\cap \sX_i$, $x_iv\in E(\sA_i^{good})$,}
\\0 &\text{otherwise,}
\end{cases}
\end{align*}
and we add $\omega_\Gamma$ to $\cW$. 
Observe that $\omega_\Gamma(M)$ accounts for the $\omega$-weight of edges $x_iv\in E(\sA_i^{good})$ such that $N_{\cH_+}(x_i)\cap\sX_0\in\sX_0^\sigma$ and thus $x_iv\notin \Gamma$. 
Hence $\omega(\Gamma)=\omega(E(\sA_i^{good}))-\omega_\Gamma(M)$.\COMMENT{Note that $\omega(\Gamma)$ incorporates the weight of edges in $\sA_i^{good}$ whose neighbours in $\sX_0$ are left unembedded by $\sigma$. That is ``unembedded weight''$=$``total weight''$-$``embedded weight''.}
For every vertex $x_0\in\sX_0$, we  have $\dg_{\sA_0}(x_0)=(d_0\pm3\eps)n$.
Hence, every edge $x_iv\in E(\sA_i^{good})$ contributes weight $\omega(x_iv)\cdot (d_0\pm3\eps)n$ to $\omega_\Gamma(E(\sA_0))$, and we obtain
\begin{align*}
\omega_\Gamma(E(\sA_0))=\omega(E(\sA_i^{good}))(d_0\pm3\eps)n.
\end{align*}
If $\omega(E(\sA_i^{good}))\geq\eps n$, then $\omega_\Gamma(E(\sA_0))\geq n^{1+\eps}\geq s\Delta^{1+\eps^2}\geq \normV{\omega_\Gamma}{1}\Delta^{1+\eps^2}$, and 
thus~\eqref{eq:w packing size} implies that
\begin{align}\label{eq:weight omega_Gamma}
\omega_\Gamma(M)
=(1\pm\heps)\frac{\omega_\Gamma(E(\sA_0))}{d_0n}
=(1\pm2\heps)\omega(E(\sA_i^{good}))\pm\heps n.
\end{align}
Again, if $\omega(E(\sA_i^{good}))<\eps n$, then~\eqref{eq:w bound codegree} implies that~\eqref{eq:weight omega_Gamma} also holds in this case.

Hence, we conclude that 
\begin{align}\label{eq:weight Gamma}
\omega(\Gamma)
=\omega(E(\sA_i^{good}))-\omega_\Gamma(M)
\stackrel{\eqref{eq:weight omega_Gamma}}{\leq} 2\heps \omega(E(\sA_i^{good}))+\heps n.
\end{align}

In order to bound $\omega(\Lambda)$, we use~\eqref{eq:bound random sparsification} and add $\omega|_{E(\sA_i^{bad})}$ to $\cW^{bad}$.
Then~\eqref{eq:bound random sparsification} implies that 
\begin{align}\label{eq:weight Lambda}
\omega(\Lambda)=(1\pm\eps)d_A\omega(E(\sA_i^{bad}))\pm\eps n.
\end{align}
Finally, equations~\eqref{eq:weight omega_0},~\eqref{eq:weight A^new edges},~\eqref{eq:weight Gamma} and~\eqref{eq:weight Lambda} yield that
\begin{align}
\omega(E(\sA_i^{new}))\nonumber
&= (1\pm2\heps)d_A\omega(E(\sA_i^{good}))
+(1\pm\eps)d_A\omega(E(\sA_i^{bad}))
\pm2\heps \omega(E(\sA_i^{good}))
\pm3\heps n
\\&=(1\pm\eps'^2)d_A\omega(E(\sA_i))\pm\eps'^2n.
\label{eq:bound updated weight linear}
\end{align} 
This establishes~\ref{emb lem:edge weight} for all $\omega\in\cW_{edge}\cup\cW_{edge}'$.

\begin{substep}\label{step:IV}
Checking~\ref{emb lem:boundedness}
\end{substep}

We show that for the updated edge set labelling $\psi[\sigma]$, we have $\Delta_{\psi[\sigma]}(\sA_i^{new})\leq (1+\eps')d_id_A|V_i|$ for every $i\in N_{R_A}(0)$.
Recall that we defined $\psi[\sigma]$ in Definition~\ref{def:updated labelling} such that for $xv\in E(A^{H,new}_i)$,\COMMENT{In fact, Definition~\ref{def:updated labelling} is stated for $xv\in E(A_i^H[\sigma])$. Observe that $A^{H,new}_i\subseteq A^{H,\ast}_i$ is a subgraph of $A_i^H[\sigma]$.} we have $\psi[\sigma](xv)=\psi(xv)\cup \{\sigma(x_0)v \}$, if $x$ has an $H$-neighbour $x_0\in\sX_0^\sigma$, and otherwise $\psi[\sigma](xv)=\psi(xv)$.
We split the proof of~\ref{emb lem:boundedness} into two claims, where Claim~\ref{claim1} bounds the number of edges on which an `old' label of $\psi$ appears on the updated candidacy graph, and Claim~\ref{claim2} bounds the number of edges on which a `new' label that we additionally added to $\psi[\sigma]$ appears in the updated candidacy graph.
Let $\psi_i\colon E(\sA_i)\to 2^{\cE_i}$ be the (old) edge set labelling~$\psi$ restricted to~$\sA_i$ and we may assume that $|\cE_i|\leq n^4$.\COMMENT{By assumption, $\|\psi\|\leq s$. Since $|E(\sA_i)|\leq2sn^3$, we may assume that $|\cE_i|\leq n^4$.}

\begin{claim}\label{claim1}
We can add at most $n^5$ weight functions to $\cW_{edge}'$ to ensure that $\Delta_{\psi_i}(\sA_i^{new})\leq (1+\eps')d_id_A|V_i|$ for every $i\in N_{R_A}(0)$.
\end{claim}
\claimproof
For all $i\in N_{R_A}(0)$ and $e\in \cE_i$, let $\omega_e\colon E(\sA_i)\to\{0,1\}$ be such that $\omega_e(x_iv_i):=\mathbbm{1}_{\{e\in\psi_i(x_iv_i) \}}$ and we add $\omega_e$ to $\cW_{edge}'$.\COMMENT{Since $\psi_i$ is defined as a star, also $\omega_e$ is a star weight function.}
By assumption (see~\ref{item:P boundedness}), we have $\Delta_\psi(\sA_i)\leq (1+\eps)d_i|V_i|$, which implies that $\omega_e(E(\sA_i))\leq (1+\eps)d_i|V_i|$.
Since~\eqref{eq:bound updated weight linear} in Step~\ref{step:III} is also valid for $\omega_e\in\cW_{edge}'$, we conclude that $e$ appears on at most 
$$(1+\eps'^2)d_A(1+\eps)d_i|V_i|+\eps'^2n\leq (1+\eps')d_id_A|V_i|$$
edges of $\sA_i^{new}$, which completes the proof of Claim~\ref{claim1}.
\endclaimproof

\begin{claim}\label{claim2}
We can add at most $n^3$ weight functions to $\cW$ to ensure that each $e\in E(G[V_0,V_i])$ appears on at most $(1+\eps')d_id_A|V_i|$ edges of $\sA_i^{new}$ for every $i\in N_{R_A}(0)$.
\end{claim}
\claimproof
For all $i\in N_{R_A}(0)$ and $e=v_0v_i\in E(G[V_0,V_i])$, we proceed as follows.
Let $N:= N_{\sA_0}(v_0)\cap N_\cH(N_{\sA_i}(v_i))$.
We define a weight function $\omega_e\colon E(\sA_0)\to \{0,1\}$ by $\omega_e(xv):=\mathbbm{1}_{\{v=v_0 \text{ and }x\in N\}}$
for every $xv\in E(\sA_0)$, and we add $\omega_e$ to~$\cW$.
Then, $e$ appears on $\omega_e(M)$ edges of~$\sA_i^{new}$.
Observe that
\begin{align}
\omega_e(E(\sA_0))
=|N|
&\stackrel{\hphantom{\ref{item:P triple intersection}}}{=}
\sum_{H\in\cH}e_{H}(N_{A_0^H}(v_0),N_{A^H_i}(v_i))\nonumber
\\&\stackrel{\ref{item:P triple intersection}}{=}
\sum_{H\in\cH}(d_0d_i\pm\eps)e_{H}(X_0^H,X^H_i)
\leq (d_0d_id_A+{2\eps'^{3/2}})n^2,\label{eq:bound N}
\end{align}
where the last inequality holds because $\sum_{H\in \cH}e_{H}(X_0^H,X^H_i)\leq {(d_A+\eps'^{3/2})n^2}$, by assumption.
With~\eqref{eq:w bound codegree}, we obtain that 
\begin{align*}
\omega_e(M)\leq
\max\Big\{(1+\heps)\frac{\omega_e(E(\sA_0))}{d_0n},n^\eps \Big\}
\stackrel{\eqref{eq:bound N}}{\leq}(1+\eps')d_id_A|V_i|,
\end{align*}
which completes the proof of Claim~\ref{claim2}.
\endclaimproof

\begin{substep}
Checking~\ref{emb lem:codegree}
\end{substep}

Recall that $\psi_i\colon E(\sA_i)\to 2^{\cE_i}$ denotes the edge set labelling~$\psi$ restricted to~$\sA_i$.
For each $i\in[r]$, $e=v_0v_i\in E(G[V_0,V_i])$, and $f\in\cE_i$, we show that $\{e,f\}$ appears on at most $\sqrt{n}$ edges of $\sA^{new}_i$.
This will imply~\ref{emb lem:codegree} because any set $\{e',f'\}\in\binom{\cE_i}{2}$ appears also on at most $\Delta^c_{\psi}(\sA_i)\leq\sqrt{n}$ edges of $\sA^{new}_i$, and no two edges of $E(G[V_0,V_i])$ appear together as a label on an edge of~$\sA^{new}_i$.
Let $$\sX_0^f:=N_{\cH_+}(\{x_i\in\sX_i\colon x_iv_i\in E(\sA_i),f\in\psi_i(x_iv_i)\})\cap \sX_0.$$
We define a weight function $\omega_{e,f}\colon E(\sA_0)\to\{0,1\}$ by $\omega_{e,f}(xv):=\mathbbm{1}_{\{v=v_0 \text{ and } x\in\sX_0^f\}}$ for every $xv\in E(\sA_0)$ and add $\omega_{e,f}$ to $\cW$.
Since $\Delta_{\psi}(\sA_i)\leq (1+\eps)d_i|V_i|$ by~\ref{item:P boundedness}, we obtain that $\omega_{e,f}(E(\sA_0))\leq 2n.
$
Note that $\{e,f\}$ appears on at most $\omega_{e,f}(M)$ edges of $\sA^{new}_i$. 
Now,~\eqref{eq:w bound codegree} implies that $\omega_{e,f}(M)\leq n^\eps$, which establishes~\ref{emb lem:codegree}.

\begin{substep}
Checking~\ref{emb lem:set testers}
\end{substep}

For each $(W,Y_1,\ldots,Y_\ell)\in\cW_{set}$ with $H_1,\ldots,H_\ell\in\cH$ such that $Y_j\subseteq X_0^{H_j}$, we define 
\begin{align*}
E_{(W,Y_1,\ldots,Y_\ell)}:=\left\{\textstyle\bigcup_{j\in[\ell]}\{xy_j \}\colon xy_j\in E(A_0^{H_j}[W,Y_j])\text{ for all $j\in[\ell]$} \right\}\subseteq \binom{E(\sA_0)}{\ell},
\end{align*}
and a weight function $\omega_{(W,Y_1,\ldots,Y_\ell)}\colon\binom{E(\sA_0)}{\ell}\to\{0,1\}$ by $\omega_{(W,Y_1,\ldots,Y_\ell)}(\{e_1,\ldots,e_\ell\}):=\mathbbm{1}_{\{ \{e_1,\ldots,e_\ell\}\in E_{(W,Y_1,\ldots,Y_\ell)} \}}$ and we add $\omega_{(W,Y_1,\ldots,Y_\ell)}$ to $\cW$.\COMMENT{Note that $\omega_{(W,Y_1,\ldots,Y_\ell)}$ corresponds to a clean weight function.}
Observe that 
\begin{align}\label{eq:weight set testers}
\omega_{(W,Y_1,\ldots,Y_\ell)}(M)=\left|W\cap \textstyle\bigcap_{j\in [\ell]}\sigma(Y_j\cap \sX_0^\sigma)\right|.
\end{align}
In view of the statement, we may assume that $|W|,|Y_j|\geq{\eps'^2 n}$ for all $j\in[\ell]$.\COMMENT{Otherwise $|W\cap \bigcap_{j\in [\ell]}\sigma(Y_j\cap \sX_0^\sigma)|\leq \eps'^2 n$ and $\eps'^2 n= |W||Y_1|\cdots |Y_\ell|/n^\ell \pm \eps' n$.}
Since $\ell\leq s$ and $A_0^H$ is $(\eps,d_0)$-super-regular for every $H\in\cH$, we obtain with Fact~\ref{fact:regularity} that are at most $\eps^{1/2}n$ vertices in~$W$ that do not have $(d_0\pm\eps)|Y_j|$ many neighbours in $Y_j$ for every $j\in[\ell]$.
Hence we obtain that
\begin{align}\label{eq:weight W_set}
\omega_{(W,Y_1,\ldots,Y_\ell)}(E(\sA_0))=|E_{(W,Y_1,\ldots,Y_\ell)}|=(d_0^\ell\pm{\eps'^2})|W||Y_1|\cdots|Y_\ell|.
\end{align}
\COMMENT{For each of the $|W|\pm\eps^{1/2}n$ good vertices $x$ in $W$, we have $(d_0\pm\eps)^\ell|Y_1|\cdots|Y_\ell|$ choices to select an $\ell$-tuple from $E_{(W,Y_1,\ldots,Y_\ell)}$ containing $x$.
Each of the $\eps^{1/2}n$ bad vertices in $W$ contributes at most $(2n)^\ell$ to $E_{(W,Y_1,\ldots,Y_\ell)}$.
That is,
$$|E_{(W,Y_1,\ldots,Y_\ell)}|
=(|W|\pm\eps^{1/2}n)(d_0\pm\eps)^\ell|Y_1|\cdots|Y_\ell|\pm 2^\ell\eps^{1/2}n^{\ell+1}
=(d_0^\ell\pm\eps'^2)|W||Y_1|\cdots|Y_\ell|,$$
where the last equality holds because $|W|,|Y_j|\geq\eps' n$ for all $j\in[\ell]$.
}For $k\in[\ell]$, any set of $k$ edges $\{e_1,\ldots,e_{k}\}$ is contained in at most $(2n)^{\ell-k}$ $\ell$-tuples in $E_{(W,Y_1,\ldots,Y_\ell)}$,\COMMENT{Let $\{e_1,\ldots,e_{k}\}$ be $k$ edges all containing one specific $x\in W$.
This also fixes the endpoint $x$ for the edges $\{e_{k+1},\ldots,e_\ell\}$. For each of these $\ell-k$ edges, there are at most $2n$ choices for its other endpoint.} which implies that
\begin{align*}
\normV{\omega_{(W,Y_1,\ldots,Y_\ell)}}{k}\Delta^{k+\eps^2}\leq (2n)^{\ell-k}\Delta^{k+\eps^2}
\stackrel{\eqref{eq:weight W_set}}{\leq} \omega_{(W,Y_1,\ldots,Y_\ell)}(E(\sA_0)).
\end{align*}
Hence, by~\eqref{eq:w packing size}, we conclude that 
\begin{align*}
\omega_{(W,Y_1,\ldots,Y_\ell)}(M)
=(1\pm\heps)\frac{\omega_{(W,Y_1,\ldots,Y_\ell)}(E(\sA_0))}{(d_0n)^\ell}
\stackrel{\eqref{eq:weight W_set}}{=}\frac{|W||Y_1|\cdots|Y_\ell|}{n^\ell}\pm\eps'n,
\end{align*}
which establishes~\ref{emb lem:set testers} by~\eqref{eq:weight set testers}.

\begin{substep}
Checking~\ref{emb lem:vertex testers}
\end{substep}

We add $\cW_{edge}$ to $\cW$ and fix some $\omega\in\cW_{edge}$.
If $\omega(E(\sA_0))\leq n^{1+\eps/2}$, then we obtain by~\eqref{eq:w bound codegree} that $\omega(M)\leq n^{\eps}$ and thus, $\omega(M)=(1\pm\eps'){\omega(E(\sA_0))}/{d_0n}\pm \eps'n$.
If $\omega(E(\sA_0))\geq n^{1+\eps/2}$, then we obtain by~\eqref{eq:w packing size} that
$\omega(M)=(1\pm\heps)\omega(E(\sA_0))/d_0n$.
This establishes~\ref{emb lem:vertex testers} and completes the proof of Lemma~\ref{lem:packing lemma}.
\end{proof}

\section{Proof of the main result}\label{sec:mainpart}

The following lemma is very similar to Theorem~\ref{thm:main_new}.
We only require additionally that all graphs in $\cH$ only span a matching between two clusters that is either empty or not too small.
This reduction has already been used in~\cite{RR:99} (and in several other extensions of the blow-up lemma) and it is also not complicated in our framework.

\begin{lemma}\label{lem:main matching}
Let $1/n\ll\eps\ll\alpha,d$ and $1/n\ll 1/r$. 
Suppose $(\cH,G,R,\cX,\cV,\phi^\ast)$ is an  $(\eps,d)$-super-regular, $\alpha^{-1}$-bounded and $(\eps,\alpha)$-linked extended blow-up instance,
$|V_i|=(1\pm \eps)n$ for all $i\in [r]$, 
and $|\cH|\leq \alpha^{-1}n$.
Suppose that $\sum_{H\in \cH} e_H(X_{i}^H,X_j^H) \leq (1-\alpha)dn^2$ for all $ij \in E(R)$
and $H[X_i^H,X_j^H]$ is a matching of size at least $\alpha^{2}n$ if $ij\in E(R)$ and empty if $ij\in\binom{[r]}{2}\sm E(R)$ for each $H\in\cH$.
Suppose $\cW_{set},\cW_{ver}$ are sets of  $\alpha^{-1}$-set testers and $\alpha^{-1}$-vertex testers of size at most $n^{2\log n}$, respectively.
Then there is a packing $\phi$ of $\cH$ in  $G$ which extends $\phi^\ast$ 
such that 
\begin{enumerate}[label={\rm (\roman*)}]
	\item\label{item:partition1} $\phi(X_i^H)=V_i$ for all $i\in [r]_0$ and $H\in \cH$;
	\item\label{item:set testers1} $|W\cap \bigcap_{j\in [\ell]}\phi(Y_j)|= |W||Y_1|\cdots |Y_\ell|/n^\ell \pm \alpha n$ for all $(W,Y_1,\ldots,Y_\ell)\in \cW_{set}$;
	\item\label{item: vertex tester1} $\omega(\bigcup_{H\in\cH}X_i^H\cap \phi^{-1}(v))=\omega(\bigcup_{H\in\cH}X_i^H)/n\pm \alpha n$ for all $(v,\omega)\in \cW_{ver}$ and $v\in V_i$.
\end{enumerate}
\end{lemma}

We first prove our main result (Theorem~\ref{thm:main_new}) assuming Lemma~\ref{lem:main matching}.

\lateproof{Theorem~\ref{thm:main_new}}
We choose a new constant $\beta$ such that $\eps\ll\beta\ll \alpha,d$.
For each $(W,Y_1,\ldots,Y_\ell)\in\cW_{set}$ with $W\subseteq V_i$, $i\in[r]$ and $k\in[\ell]$, let $\omega_{Y_{k}}\colon \bigcup_{H\in\cH}X_i^H\to\{0,1\}$ be such that $\omega_{Y_{k}}(x)=\mathbbm{1}_{\{x\in Y_{k} \}}$, and let $\cW_Y$ be the set containing all those weight functions.\COMMENT{$\cW_{Y}$ is used to control the size of the refined $Y_{k,j}$.}
{We delete from every $H\in \cH$ the set $X_0^H$ and apply Lemma~\ref{lem:split} to this collection of graphs and the set  of weight functions $\cW^\ast:=\{\omega\colon (v,\omega)\in\cW_{ver}  \}\cup \cW_Y$, which yields a refined partition of $\cH$;
to be more precise, for all $H\in\cH$ and $i\in[r]$, we obtain a partition $(X^H_{i,j})_{j\in[\beta^{-1}]}$ of $X^H_i$ satisfying~\ref{item:split1}--\ref{item:split3} of Lemma~\ref{lem:split}.
Let $\cX'$ be the collection of vertex partitions of the graphs in $\cH$ given by $(X_0^H)_{H\in\cH}$ and $(X^H_{i,j})_{H\in\cH,i\in[r],j\in[\beta^{-1}]}$.}
In particular, Lemma~\ref{lem:split}\ref{item:split4} yields that
\begin{align}\label{eq:thm split 4}
\omega(X_{i,j}^H)=\beta\omega(X_i^H)\pm\beta^{3/2}n, \text{ for all $H\in\cH,\omega\in\cW^\ast,i\in[r],j\in[\beta^{-1}]$.}
\end{align}
Let $R'$ be the graph with vertex set $[r]\times [\beta^{-1}]$ and two vertices $(i,j),(i',j')$ are joined by an edge if $ii'\in E(R)$.
Note that $\Delta(R')\leq\alpha^{-1}\beta^{-1}$ because $\Delta(R)\leq\alpha^{-1}$.

According to the refinement $\cX'$ of $\cX$, we claim that there exists a refined partition $\cV'$ of $\cV$ consisting of the collection of $V_0$ together with $(V_{i,j})_{i\in[r],j\in[\beta^{-1}]}$, where $(V_{i,j})_{j\in[\beta^{-1}]}$ is a partition of~$V_i$ for every $i\in[r]$ such that
\begin{enumerate}[label=(\alph*)]
\item\label{item:refined W} $|W\cap V_{i,j}|=\beta|W|\pm\beta^{3/2}n$ for all $(W,Y_1,\ldots,Y_\ell)\in \cW_{set}$ and $j\in[\beta^{-1}]$ with $W\subseteq V_i, \ell\in[\alpha^{-1}]$;
\item\label{item:refined blowup instance} $(\cH,G,R',\cX',\cV',\phi_0)$ is an $(\eps^{1/2},d)$-super-regular, $\beta^{-2}$-bounded and $(\eps^{1/2},\alpha/2)$-linked extended blow-up instance.
\end{enumerate}

Indeed, the existence of $\cV'$ follows by a simple probabilistic argument.
For each $i\in[r]$, let $\tau_i\colon V_i\to[\beta^{-1}]$ where $\tau_i(v)$ is chosen uniformly at random for every $v\in V_i$, all independently, and let $V_{i,j}:=\{v\in V_i\colon \tau_i(v)=j\}$ for every $j\in[\beta^{-1}]$.
%We label the vertices in $(\bigcup \cV)\setminus V_0$ independently and uniformly at random with labels in $[\beta^{-1}]$.
Chernoff's inequality and a union bound imply that~\ref{item:refined W} holds simultaneously together with the following properties with probability at least $1-\eul^{-\sqrt{n}}$:
\begin{itemize}
%\item\ref{item:refined W};
\item$G[V_{i,j},V_{i',j'}]$ is $(\eps^{1/2},d)$-super-regular for all $ii'\in E(R), j,j'\in[\beta^{-1}]$;
\item$|\bigcap_{x_0\in X_0^H\cap N_H(x)} N_G(\phi_0(x_0))\cap V_{i,j}| \geq \alpha/2 |V_{i,j}|$ for all $x\in X_{i,j}^H, i\in [r],j\in[\beta^{-1}], H\in \cH$.
\end{itemize}
Standard properties of the multinomial distribution yield that $|V_{i,j}|=|X_{i,j}^H|$ for all $i\in[r]$, $j\in[\beta^{-1}]$, $H\in\cH$ with probability at least $\Omega(n^{-r\beta^{-1}})$.
To see in~\ref{item:refined blowup instance} that the instance is $(\eps^{1/2},\alpha/2)$-linked, observe further that the number of vertices in $X_{i,j}^H$ that have a neighbour in $X_0^H$ is at most $\eps |X_i^H|\leq\eps^{1/2}|X_{i,j}^H|$ and $\sum_{H\in \cH}|N_H(\phi_0^{-1}(v_0),\phi_0^{-1}(v_0'))\cap X_{i,j}^H|\leq \eps|V_i|^{1/2}\leq \eps^{1/2}|V_{i,j}|^{1/2}$ for all $i\in[r]$, $j\in[\beta^{-1}]$ and distinct $v_0,v_0'\in V_0$.
Thus, for every $i\in[r]$, there exists a partition $(V_{i,j})_{j\in[\beta^{-1}]}$ of $V_i$ satisfying~\ref{item:refined W} and~\ref{item:refined blowup instance}.
Let $n':=\beta n$.

Next we show how to lift the vertex and set testers from the original blow-up instance to the just defined blow-up instance.
For each $(W,Y_1,\ldots,Y_\ell)\in \cW_{set}$ and distinct $H_1,\ldots,H_\ell\in\cH$ such that $W\subseteq V_i$ for some $i\in[r]$ and $Y_{k}\subseteq X_i^{H_{k}}$ for all $k\in[\ell]$, 
we define $(W_j,Y_{1,j},\ldots,Y_{\ell,j})$ by setting $W_j:=W\cap V_{i,j}$ and $Y_{k,j}:=Y_{k}\cap X_{i,j}^{H_{k}}$ for all $j\in[\beta^{-1}], k\in[\ell]$. 
By~\ref{item:refined W}, we conclude that $|W_j|=\beta |W|\pm\beta^{3/2}n$, 
and by~\eqref{eq:thm split 4}, we have that $|Y_{k,j}|=\omega_{Y_{k}}(X_{i,j}^{H_{k}})=\beta \omega_{Y_{k}}(X_i^{H_{k}}) \pm\beta^{3/2}n =\beta |Y_{k}|\pm\beta^{3/2}n$.

Let $\cW_{set}':=\{(W_j,Y_{1,j},\ldots,Y_{\ell,j})\colon 
j\in[\beta^{-1}],(W,Y_1,\ldots,Y_\ell)\in \cW_{set} \}$.
For each $(v,\omega)\in\cW_{ver}$ with $v\in V_{i,j}$, let $\omega':=\omega|_{\bigcup_{H\in \cH}X_{i,j}^H}$ and $\cW_{ver}':=\{(v,\omega')\colon (v,\omega)\in \cW_{ver} \}$.

Next, we add some edges to the graphs in $\cH$ ensuring that all matchings between two clusters are either empty or of small linear size.
To this end, we add a minimum number of edges to $H[X_{i,j}^H,X_{i',j'}^H]$ for all $\{(i,j),(i',j')\}\in E(R')$ and $H\in\cH$
such that the obtained supergraph $H'[X_{i,j}^{H'},X_{i',j'}^{H'}]$ is a matching of size at least $\beta^4n$. 
Note that $\Delta(H')\leq\Delta(R')+\alpha^{-1}\leq\beta^{-2}$.
Let $\cH'$ be the collection of graphs $H'$ obtained in this manner.
Together with Lemma~\ref{lem:split}\ref{item:split3}, we conclude for all $\{(i,j),(i',j')\}\in E(R')$ that
\begin{align*}
\sum_{H'\in \cH'}e_{H'}(X_{i,j}^{H'},X_{i',j'}^{H'})
&\leq {2\beta^4n\cdot \alpha^{-1}n}+ \sum_{H\in\cH}e_H(X_{i,j}^{H},X_{i',j'}^{H})
\\&\leq \beta^3n^2+ \beta^2\sum_{H\in \cH}e_H(X_{i}^H,X_{i'}^H) + n^{5/3}
\leq (1-\alpha/2)dn'^2.
\end{align*}
Obviously, it suffices to construct a packing of $\cH'$ into $G$ which extends $\phi_0$ and satisfies Theorem~\ref{thm:main_new}\ref{item:partition}--\ref{item: vertex tester}.
By~\ref{item:refined blowup instance} and because $\beta\ll\alpha$, also $(\cH',G,R',\cX',\cV',\phi_0)$ is an $(\eps^{1/2},d)$-super-regular, $\beta^{-2}$-bounded and $(\eps^{1/2},{\beta^2})$-linked extended blow-up instance, and we can apply Lemma~\ref{lem:main matching} to $(\cH',G,R',\cX',\cV',\phi_0)$ with set testers $\cW_{set}'$ and vertex testers $\cW_{ver}'$ as follows:\COMMENT{Note that $n'^{2\log n'}=(\beta n)^{2\log(\beta n)}
\geq(\beta n)^{3\log n/2}\geq\beta^{-1}n^{\log n}\geq |\cW_{set}'|,|\cW_{ver}'|.$}

\begin{center}
\begin{tabular}{c|c|c|c|c}
 $n'$  & $\eps^{1/2}$ & {$\beta^2$} & $d$ & $r\beta^{-1}$
\\ \hline
 $n$  & $\eps$  & $\alpha$ & $d$  & $r$ 
\vphantom{\Big(}
\end{tabular}
\end{center} 
Hence, we obtain a packing $\phi$ of $\cH'$ in $G$ which extends $\phi_0$ such that for all $i\in[r], j\in[\beta^{-1}]$
\begin{enumerate}[label=(\Roman*)]
\item\label{item:I} $\phi(X_{i,j}^H)\subseteq V_{i,j}$ for all $H\in\cH$;
\item\label{item:II} $|W_j\cap \bigcap_{k\in [\ell]}\phi(Y_{k,j})|= |W_j||Y_{1,j}|\cdots |Y_{\ell,j}|/n'^\ell \pm {\beta^2} n'$ for all $(W_j,Y_{1,j},\ldots ,Y_{\ell,j})\in \cW_{set}'$;
\item \label{item:III} $\omega'(\bigcup_{H\in\cH}X_{i,j}^H \cap\phi^{-1}(v))=\omega'(\bigcup_{H\in\cH}X_{i,j}^H)/n'\pm{\beta^2} n'$ for all $(v,\omega')\in \cW_{ver}'$ with $v\in V_{i,j}$.
\end{enumerate}
Observe that~\ref{item:I} establishes Theorem~\ref{thm:main_new}\ref{item:partition}.

For $(W,Y_1,\ldots,Y_\ell)\in\cW_{set}$, we conclude that
\begin{align*}
\left|W\cap \bigcap_{k\in [\ell]}\phi(Y_{k})\right| 
&\stackrel{\hphantom{\rm \ref{item:II},\eqref{eq:thm split 4}}}{=} \sum_{j\in[\beta^{-1}]}\left|W_j\cap \bigcap_{k\in [\ell]}\phi(Y_{k,j})\right|\\
&\stackrel{\rm \ref{item:II},\eqref{eq:thm split 4}}{=} \sum_{j\in[\beta^{-1}]}\left(\frac{\beta^{\ell+1}\left(|W||Y_1|\cdots|Y_\ell|\pm\beta^{1/3}n^{\ell+1}\right)}{(\beta n)^\ell}\pm{\beta^2} n'\right)\\
&\stackrel{\hphantom{\rm \ref{item:II},\eqref{eq:thm split 4}}}{=}|W||Y_1|\cdots|Y_\ell|/n^\ell \pm\alpha n.
\end{align*}
\COMMENT{Note that by~\ref{item:II}, we have $$|W_j||Y_{1,j}|\cdots |Y_{\ell,j}|=\beta^{\ell+1}\Big(|W||Y_1|\cdots|Y_\ell|\pm\beta^{1/3}n^{\ell+1}\Big)$$
because $|W_j|=\beta |W|\pm\beta^{3/2}n$ and $|Y_{k,j}|=\beta |Y_{k}|\pm\beta^{3/2}n$ by~\eqref{eq:thm split 4}. Hence,
\begin{align*}
\sum_{j\in[\beta^{-1}]}\left|W_j\cap \bigcap_{k\in [\ell]}\phi(Y_{k,j})\right|
&=\sum_{j\in[\beta^{-1}]}\left(\frac{\beta^{\ell+1}\left(|W||Y_1|\cdots|Y_\ell|\pm\beta^{1/3}n^{\ell+1}\right)}{(\beta n)^\ell}\pm{\beta^2} n'\right)
\\&=\frac{|W||Y_1|\cdots|Y_\ell|}{n^\ell}\pm\beta^{1/3}n\pm {\beta^2} n
\\&=|W||Y_1|\cdots|Y_\ell|/n^\ell \pm\alpha n.
\end{align*}
}Hence, Theorem~\ref{thm:main_new}\ref{item:set testers} holds.

For $(v,\omega)\in\cW_{ver}$ with $v\in V_{i,j}$ and its corresponding tuple $(v,\omega')\in \cW_{ver}'$, we conclude that 
\begin{align*}
\omega(\textstyle\bigcup_{H\in\cH}X_{i}^H\cap \phi^{-1}(v))
&\stackrel{~\rm\ref{item:I}~}{=}\omega'(\textstyle\bigcup_{H\in\cH}X_{i,j}^H\cap \phi^{-1}(v))
\stackrel{\rm \ref{item:III}}{=}
\omega'(\bigcup_{H\in\cH}X_{i,j}^H)/n'\pm{\beta^2} n'
\\
&\stackrel{\eqref{eq:thm split 4}}{=}
\frac{\beta \omega(\textstyle\bigcup_{H\in\cH} X_i^H)\pm\beta^{4/3}n^2}{\beta n}\pm{\beta^2} n'
= \omega(\textstyle\bigcup_{H\in\cH} X_i^H)/n\pm\alpha n.
\end{align*}
\COMMENT{\begin{align*}
&\frac{\omega'(\bigcup_{H\in\cH}X_{i,j}^H)}{n'}\pm\beta^2 n'
=\frac{\sum_{H\in \cH}\omega(X_{i,j}^H)}{n'}\pm\beta^2 n'
\stackrel{\eqref{eq:thm split 4}}{=}
\\
&\frac{\sum_{H\in \cH}\left(\beta\omega(X_{i}^H)\pm \beta^{3/2}n\right)}{n'}\pm\beta^2 n'
=\frac{\beta \omega(\textstyle\bigcup_{H\in\cH} X_i^H)\pm\alpha^{-1}\beta^{3/2}n^2}{\beta n}\pm\beta^2 n'
= \omega(\textstyle\bigcup_{H\in\cH} X_i^H)/n\pm\alpha n.
\end{align*}
}This yields Theorem~\ref{thm:main_new}\ref{item: vertex tester} and completes the proof.
\endproof

Theorem~\ref{thm:quasirandom} can be easily deduced from Theorem~\ref{thm:main_new} {by randomly partitioning $G$ and applying Lemma~\ref{lem:split} to $\cH$ with $r=1$.}
In particular, the proof is very similar to the proof of Theorem~\ref{thm:main_new} and therefore omitted.%
\COMMENT{
\lateproof{Theorem~\ref{thm:quasirandom}}
We choose a new constant $\beta$ with $1/n\ll\beta\ll\alpha$.
As given in the statement, let $G$ be an $(\eps,d)$-quasirandom graph on $n$ vertices and $\cH$ a collection of graphs on at most $n$ vertices with $|\cH|\leq\alpha^{-1}n$, $e(\cH)\leq(1-\alpha)e(G)$, and $\Delta(H)\leq\alpha^{-1}$ for al $H\in\cH$. 
Let $\cW_{set}$ and $\cW_{ver}$ be the given sets of set testers and vertex testers, respectively.
For each $(W,Y_1,\ldots,Y_\ell)\in\cW_{set}$ and $k\in[\ell]$, let $\omega_{Y_{k}}\colon \bigcup \cX\to\{0,1\}$ be such that $\omega_{Y_{k}}(x)=\mathbbm{1}_{\{x\in Y_{k} \}}$, and let $\cW_Y$ be the set containing all those weight functions.
By applying Lemma~\ref{lem:split} to $\cH$ with trivial vertex partition $V(H)$ for every $H\in\cH$ and the set of weight functions $\{\omega\colon (v,\omega)\in\cW_{ver}  \}\cup \cW_Y$, we obtain a partition $(X_i^H)_{i\in[\beta^{-1}]}$ of $V(H)$ for all $H\in\cH$ satisfying the conclusions~\ref{item:split1}--\ref{item:split3} of Lemma~\ref{lem:split}.
Let $R$ be the graph with vertex set $[\beta^{-1}]$ and two vertices $i,i'$ are joined by an edge if $H[X_i,X_{i'}]$ is non-empty for some $H\in\cH$. 
Clearly, $\Delta(R)\leq \beta^{-1}$.\\
Accordingly, we claim that there exists a partition $(V_i)_{i\in[\beta^{-1}]}$ of $V(G)$ such that
\begin{enumerate}[label=(\alph*)]
\item\label{item:quasi refined W} $|W\cap V_{i}|=\beta|W|\pm\beta^{3/2}n$ for all $(W,Y_1,\ldots,Y_\ell)\in \cW_{set}$ and $i\in[\beta^{-1}], \ell\in[\alpha^{-1}]$;
\item\label{item:quasi blowup instance} $(\cH,G,R,(X_i)_{i\in[\beta^{-1}],X\in\cH},(V_i)_{i\in[\beta^{-1}]})$ is an $(2\eps,d)$-super-regular, $\beta^{-2}$-bounded blow-up instance.
\end{enumerate}
Indeed, the existence of such a partition can be seen by assigning every vertex in $V(G)$ uniformly at random to some $V_i$ for $i\in[\beta^{-1}]$, all independently. Apply Chernoff's bound and standard properties of multinomial distribution ...
Let $n':=\beta n$.\\
Next, we adjust the vertex testers and set testers for the quasirandom graph to the just defined blow-up instance.
For each $(W,Y_1,\ldots,Y_\ell)\in \cW_{set}$ and distinct $H_1,\ldots,H_\ell\in\cH$ such that $Y_{k}\subseteq V(H_k)$ for all $k\in[\ell]$, 
we define $(W_i,Y_{1,i},\ldots,Y_{\ell,i})$ by setting $W_i:=W\cap V_{i}$ and $Y_{k,i}:=Y_{k}\cap X_{i}^{H_{k}}$ for all $i\in[\beta^{-1}], k\in[\ell]$. 
By~\ref{item:quasi refined W}, we conclude that $|W_i|=\beta |W|\pm\beta^{3/2}n$, 
and by Lemma~\ref{lem:split}\ref{item:split4}, we have that $|Y_{k,i}|=\omega_{Y_{k}}(X_{i}^{H_{k}})=\beta \omega_{Y_{k}}(V({H_{k}})) \pm\beta^{3/2}n =\beta |Y_{k}|\pm\beta^{3/2}n$.\\
Let $\cW_{set}':=\{(W_i,Y_{1,i},\ldots,Y_{\ell,i})\colon 
i\in[\beta^{-1}],(W,Y_1,\ldots,Y_\ell)\in \cW_{set} \}$.
For each $(v,\omega)\in\cW_{ver}$ with $v\in V_{i}$, let $\omega':=\omega|_{\bigcup_{H\in \cH}X_{i}^H}$ and $\cW_{ver}':=\{(v,\omega')\colon (v,\omega)\in \cW_{ver} \}$.\\
By Lemma~\ref{lem:split}\ref{item:split3}, we conclude that 
\begin{align*}
\sum_{H\in \cH}e_H(X_i^H,X_{i'}^H)
=\beta^2\sum_{H\in \cH}e(H)\pm n^{5/3}
\leq (1-\alpha/2)d(n')^2.
\end{align*}
Hence, we can apply Theorem~\ref{thm:main_new} to the blow-up instance $(\cH,G,R,(X_i)_{i\in[\beta^{-1}],X\in\cH},(V_i)_{i\in[\beta^{-1}]})$ with set testers $\cW_{set}'$ and vertex testers $\cW_{ver}'$ as follows:
\begin{center}
\begin{tabular}{c|c|c|c|c}
 $n'$  & $2\eps$ & $\beta$ & $d$  & $\beta^{-1}$
\\ \hline
 $n$  & $\eps$  & $\alpha$ & $d$  & $r$ 
\vphantom{\Big(}
\end{tabular}
\end{center} 
Hence, we obtain a packing $\phi$ of $\cH$ in $G$ such that
\begin{enumerate}[label=(\Roman*)]
\item\label{item:quasi I} $\phi(X_{i}^H)\subseteq V_{i}$ for all $i\in[\beta^{-1}]$ and $H\in\cH$;
\item\label{item:quasi II} $|W_i\cap \bigcap_{k\in [\ell]}\phi(Y_{k,i})|= |W_i||Y_{1,i}|\cdots |Y_{\ell,i}|/(n')^\ell \pm \beta n'$ for all $(W_i,Y_{1,i},\ldots ,Y_{\ell,i})\in \cW_{set}'$;
\item \label{item:quasi III} $\omega'(\bigcup_{H\in\cH}X_{i}^H \cap\phi^{-1}(v))=\omega'(\bigcup_{H\in\cH}X_{i}^H)/n'\pm\beta n'$ for all $(v,\omega')\in \cW_{ver}'$ with $v\in V_{i}$.
\end{enumerate}
For $(W,Y_1,\ldots,Y_\ell)\in\cW_{set}$, we conclude that
\begin{eqnarray*}
\left|W\cap \bigcap_{k\in [\ell]}\phi(Y_{k})\right| 
&=& \sum_{i\in[\beta^{-1}]}\left|W_i\cap \bigcap_{k\in [\ell]}\phi(Y_{k,i})\right|\\
&\stackrel{\rm \ref{item:quasi II}}{=}& \sum_{i\in[\beta^{-1}]}\left(\frac{\beta^{\ell+1}\left(|W||Y_1|\cdots|Y_\ell|\pm\beta^{1/3}n^{\ell+1}\right)}{(\beta n)^\ell}\pm\beta n'\right)\\
&=&|W||Y_1|\cdots|Y_\ell|/n^\ell \pm\alpha n.
\end{eqnarray*}
For $(v,\omega)\in\cW_{ver}$ with $v\in V_{i}$ and its corresponding tuple $(v,\omega')\in \cW_{ver}'$, we conclude that 
\begin{eqnarray*}
\omega(\textstyle\bigcup_{H\in\cH}V(H)\cap \phi^{-1}(v))
&\stackrel{\rm\ref{item:I}}{=}&\omega'(\textstyle\bigcup_{H\in\cH}X_{i}^H\cap \phi^{-1}(v))
\stackrel{\rm \ref{item:III}}{=}
\omega'(\bigcup_{H\in\cH}X_{i}^H)/n'\pm\beta n'
\\
&=&
\frac{\beta \omega(\textstyle\bigcup_{H\in\cH} V(H))\pm\beta^{4/3}n^2}{\beta n}\pm\beta n'
= \omega(\textstyle\bigcup_{H\in\cH} V(H))/n\pm\alpha n,
\end{eqnarray*}
where we used Lemma~\ref{lem:split}\ref{item:split4} in the third estimation. 
This completes the proof of Theorem~\ref{thm:quasirandom}.
\endproof
}
We proceed with the proof of Lemma~\ref{lem:main matching}.

\lateproof{Lemma~\ref{lem:main matching}}
We split the proof into four steps.
In Step~\ref{step:partitioning}, we partition $G$ into two edge-disjoint subgraphs $G_A$ and $G_B$.
In Step~\ref{step:candidacy}, we define `candidacy graphs' that we track for the partial packing in Step~\ref{step:induction}, where we iteratively apply Lemma~\ref{lem:packing lemma} to consider the clusters in turn. We only use the edges of $G_A$ for the partial packing in Step~\ref{step:induction} such that we can complete the packing in Step~\ref{step:completion} using the edges of $G_B$ and the ordinary blow-up lemma.

We will proceed cluster by cluster in Step~\ref{step:induction} to find a function that packs almost all vertices of $\cH$ into $G$. 
Since $r$ may be much larger than $\eps^{-1}$, we need to carefully control the growth of the error term. 
We do so, by considering a proper vertex colouring $c\colon V(R)\to[T]$ of $R^3$ where $T:=\alpha^{-3}$, and choose new constants $\eps_0,\eps_1,\ldots,\eps_T,\mu,\gamma$ such that $\eps\ll\eps_0\ll\eps_1\ll\cdots\ll\eps_T\ll\mu\ll\gamma\ll \alpha,d$. 
{To obtain the order in which we consider the clusters in turn, we simply relabel the cluster indices such that the colour values are non-decreasing; that is,
$c(1)\leq\cdots\leq c(r)$.}
Note that the sets $(c^{-1}(k))_{k\in[T]}$ are independent in $R^3$. 
For $i\in[r], t\in[r]_0$, let 
\begin{align}\label{eq:def c_i m_i}
c_i(t):=\max\{\{0\}\cup\{c(j)\colon j\in N_{R}[i]\cap[t]\}\},\quad \text{and} \quad m_i(t):=|N_R(i)\cap [t]|.
\end{align}
That is, if we think of $[t]$ as the indices of clusters that have already been embedded, then $c_i(t)$ denotes the largest colour of an already embedded cluster in the closed neighbourhood of $i$ in~$R$, and $m_i(t)$ denotes the number of neighbours of $i$ in $R$ that have already been embedded.

For $t\in[r]_0$, let
\begin{align*}
\sX_{t}&:=\textstyle\bigcup_{H\in\cH}X^H_{t}, 
& \cX_t&:=\textstyle\bigcup_{\ell\in[t]_0}\sX_{\ell}, 
& \cV_t&:=\textstyle\bigcup_{\ell\in [t]_0}V_{\ell}. 
\end{align*}

For every vertex tester $(v,\omega)\in\cW_{ver}$ with $v\in V_i$ for some $i\in[r]$, we define its \defn{corresponding function $\omega_v$} on $\{x_iv_i\colon x_i\in\sX_i,v_i\in V_i \}$ by setting $\omega_v(x_iv_i):=\omega(x_i)\mathbbm{1}_{\{v_i=v\}}$.~Let 
\begin{align}\label{eq:def W_edge}
\cW_{edge}^i:=\{\omega_v \colon (v,\omega)\in\cW_{ver}, v\in V_i \}.
\end{align}

\begin{step}\label{step:partitioning}
Partitioning the edges of $G$
\end{step}

In order to reserve an exclusive set of edges for the completion in Step~\ref{step:completion}, we partition the edges of $G$ into two subgraphs $G_A$ and $G_B$.
For each edge $e$ of $G$ independently, we add $e$ to $G_B$ with probability $\gamma$ and otherwise to $G_A$. 
Let $d_A:=(1-\gamma)d$, $d_B:=\gamma d$, $\alpha_A:=(1-\gamma)^{\alpha^{-1}}\alpha/2$ and $\alpha_B:=\gamma^{\alpha^{-1}}\alpha/2$.\COMMENT{Note that a vertex $x\in\sX_i$ might have $\alpha^{-1}$ neighbours in $\sX_0$.}
Using Chernoff's inequality, we can easily conclude that with probability at least $1-1/n$ we have for all $Z\in\{A,B\}$ that
\begin{align}\label{eq:G_Z sr}
\begin{minipage}[c]{0.8\textwidth}\em
$G_Z[V_i,V_j]$ is $(2\eps,d_Z)$-super-regular for all $ij\in E(R)$,
\end{minipage}
\\
\label{eq:linked}
\begin{minipage}[c]{0.8\textwidth}\em
$|V_i\cap\bigcap_{x_0\in X_0^H\cap N_H(x)} N_{G_Z}(\phi^\ast(x_0))|\geq \alpha_Z|V_i|$ for all $x\in X_i^H, i\in [r], H\in\cH$.
\end{minipage}
\ignorespacesafterend
\end{align}
Hence, we may assume that $G$ is partitioned into $G_A$ and $G_B$ such that~\eqref{eq:G_Z sr} and~\eqref{eq:linked} hold.

\begin{step}\label{step:candidacy}
Candidacy graphs
\end{step}

For $t\in[r]_0$, we call $\phi\colon \bigcup_{H\in\cH,i\in[t]_0}\hX^H_i\to\cV_t$ a \emph{$t$-partial packing} if $\hX^H_i\subseteq X^H_i$, $\phi|_{\hX^H_0}=\phi^\ast|_{X_0^H}$, and $\phi(\hX^H_i)\subseteq V_i$ for all $H\in\cH, i\in[t]_0$ such that $\phi$ is a packing of $(H[\hX^H_0\cup\ldots\cup \hX^H_t])_{H\in\cH}$ into $G_A[\cV_t]$. Note that $\hX^H_0=X^H_0$, and $\phi|_{\hX^H_i}$ is injective for all $H\in\cH$ and $i\in[t]_0$.
For convenience, we often write
\begin{align*}
\cX_t^{\phi} :=\textstyle\bigcup_{H\in\cH,i\in[t]_0}\hX^H_i. 
\end{align*}

Suppose $t\in[r]_0$ and $\phi_t\colon \cX_t^{\phi_t}\to \cV_t$ is a $t$-partial packing.
We introduce the notion of \emph{candidates} (with respect to $\phi_t$) for future packing rounds and track those relations in two kinds of bipartite auxiliary graphs that we call candidacy graphs:
A graph $A^H_i(\phi_t)$ with bipartition $(X^H_i,V_i)$, $i\in[r]$ that will be used to extend the $t$-partial packing $\phi_t$ to a $(t+1)$-partial packing $\phi_{t+1}$ via Lemma~\ref{lem:packing lemma} in Step~\ref{step:induction}, and a graph $B^H_i(\phi_t)$ that will be used for the completion in Step~\ref{step:completion}.
For convenience, we define $B^H_i(\phi_t)$ on a copy $(X^{H,B}_i,V_i^B)$ of $(X^H_i,V_i)$.
That is, for all $H\in\cH, i\in[r]$, let $X^{H,B}_i$ and $V_i^B$ be disjoint copies of $X^H_i$ and $V_i$, respectively. 
Let $\pi$ be the bijection that maps a vertex in $\bigcup_{H\in\cH,i\in[r]}(X^H_i\cup V_i)$ to its copy in $\bigcup_{H\in\cH,i\in[r]}(X^{H,B}_i\cup V_i^B)$.
Let $G_+$ and $H_+$ be supergraphs of $G_A$ and $H\in\cH$ with vertex partitions $(V_0,\ldots,V_r,V_1^B,\ldots,V_r^B)$ and $(X^H_0,\ldots, X^H_r,X^{H,B}_1,\ldots, X^{H,B}_r)$, respectively, and edge sets
\begin{align*}
E(G_+) &:= E(G_A)\cup \{u\pi(v)\colon uv\in E(G_B) \}, \\
E(H_+) &:= E(H)\cup \{u\pi(v)\colon uv\in E(H) \}.
\end{align*}
Let $R_B$ be the graph on $[r]\cup\{1^B,,\ldots,r^B\}$ with edge set $E(R_B) :=\{ij^B\colon ij\in E(R) \}.$
By taking copies $(X^{H,B}_i,V_i^B)$ for all $(X^H_i,V_i)$ and defining the candidacy graphs $B^H_i(\phi_t)$ on these copies, and by enlarging $G$, $H$ and $R$ accordingly to $G_+$, $H_+$ and $R\cup R_B$, we will be able to update the candidacy graphs $A^H_i(\phi_t)$ and $B^H_i(\phi_t)$ simultaneously in Step~\ref{step:induction} when we apply Lemma~\ref{lem:packing lemma} in order to extend $\phi_t$ to a $(t+1)$-partial packing $\phi_{t+1}$.

We now define $A^H_i(\phi_t)$ and $B^H_i(\phi_t)$.
Let $X^{H,A}_i:=X^H_i$ and $V_i^A:=V_i$ for all $H\in\cH$, $i\in[r]$.
For $Z\in\{A,B\}$, $H\in\cH$ and $i\in[r]$, we say that $v\in V_i^Z$ is a \emph{candidate for $x\in X_i^{H,Z}$ given~$\phi_t$}~if
\begin{align}\label{eq:candidates}
\phi_t(N_{H_+}(x)\cap \cX_t^{\phi_t})\subseteq N_{G_+}(v).
\end{align}
For all $Z\in\{A,B\}$, let $Z^H_i(\phi_t)$ be a bipartite graph with vertex partition $(X^{H,Z}_i,V_i^Z)$ and edge set\COMMENT{Note that if $Z=A$, we could restrict $i$ to $[r]\sm[t]$ since $A^H_i$ for $i\in[r]\sm[t]$ incorporate the `packing' candidacy graphs for future rounds given $\phi_t$; whereas $B^H_i$ incorporate the `absorber' candidacy graphs such that for $i\leq t$, also $B^H_i$ has to be updated when we extend $\phi_t$ to $\phi_{t+1}$.} 
\begin{align}\label{eq:def candidacy}
E(Z^H_i(\phi_t)):=\{xv\colon x\in X^{H,Z}_i, v\in V_i^Z,\text{ and $v$ is a candidate for $x$ given $\phi_t$} \}.
\end{align}
We call every spanning subgraph of $Z^H_i(\phi_t)$ a \emph{candidacy graph (with respect to $\phi_t$)}.

Furthermore, for all $H\in\cH$ and $i\in[r]$, we assign to every edge $xv\in E(A^H_i(\phi_t))$ an edge set labelling $\psi_t(xv)$ of size at most~{$\alpha^{-1}$}.
This set encodes the edges between $v$ and $\phi_t(N_{H}(x)\cap \cX_t^{\phi_t})$ in $G_A$ that are covered if we embed $x$ onto $v$; to be more precise, for all $H\in\cH$, $i\in[r]$, and every edge $xv\in E(A^H_i(\phi_t))$, we set 
\begin{align}\label{eq:def edge set labelling}
\psi_t(xv):= 
E\big(G_A\big[\phi_t(N_{H}(x)\cap \cX_t^{\phi_t}),\{v\}\big]\big). 
\end{align}
Tracking this set enables us to extend a $t$-partial packing $\phi_t$ to a $(t+1)$-partial packing $\phi_{t+1}$ by finding a conflict-free embedding (see definition in~\eqref{eq:conflict free}) in $\bigcup_{H\in\cH}A^H_{t+1}(\phi_t)$ via Lemma~\ref{lem:packing lemma}.
Since $|N_{H}(x)\cap \cX_t^{\phi_t}|\leq {\alpha^{-1}}$, we have $|\psi_t(xv)|\leq {\alpha^{-1}}$.

Before we proceed to Step~\ref{step:induction} and extend $\phi_t$ to $\phi_{t+1}$, we consider the candidacy graphs and their edge set labelling with respect to $\phi^\ast$.
\begin{claim}\label{claim exceptional}
For all $H\in\cH, i\in[r], Z\in\{A,B\}$, there exists a candidacy graph $Z^H_i\subseteq Z^H_i(\phi^\ast)$ with respect to $\phi^\ast$ (where $\sA_i:=\bigcup_{H\in \cH}A_i^H$) such that
\begin{enumerate}[label={\rm (C\ref{claim exceptional}.\arabic*)}]
\item\label{exc sr} $Z^H_i$ is $(\eps_{0},\alpha_Z)$-super-regular;
\item\label{exc boundedness} $\Delta_{\psi_0}(\sA_i)\leq \eps_0 n$;\COMMENT{In particular, $\Delta_{\psi_0}(\sA_i)\leq(1+\eps_0)\alpha_A|V_i|$}
\item\label{exc codegree} $\Delta^c_{\psi_0}(\sA_i)\leq \sqrt{n}$; 
\item\label{exc triple} $e_{H}(N_{A^H_i}(v_i),N_{A^H_j}(v_j))=(\alpha_A^2\pm\eps_{0})e_{H}(X_i^H,X_j^H)$ for all $v_iv_j\in E(G_A[V_i,V_j])$, $ij\in E(R)$;
\item \label{exc testers} $\omega_v(E(\sA_i))=\alpha_A\omega(\sX_{i})\pm\eps_0n^2$ for all $\omega_v\in\cW_{edge}^i$.
\end{enumerate}
\end{claim}
\claimproof
We fix $H\in\cH, i\in[r], ij\in E(R), Z\in\{A,B\}$, $v_iv_j\in E(G_A[V_i,V_j])$ and $\omega_v\in\cW_{edge}^i$ as defined in~\eqref{eq:def W_edge}. 
For each $k\in\{i,j\}$, let $\tX^H_k$ be the set of vertices in $X^H_k$ that have a neighbour in $X_0^H$. 
Observe that $Z^H_k(\phi^\ast)[X^H_k\sm\tX^H_k,V_k]$ is a complete bipartite graph for $k\in\{i,j\}$, and $e_{H}(N_{Z^H_i(\phi^\ast)[X^H_i\sm\tX^H_i,V_i]}(v_i),N_{Z^H_j(\phi^\ast)[X^H_j\sm\tX^H_j,V_j]}(v_j))=e_{H}(X^H_i,X^H_j)\pm 4\eps n$ because $|\tX^H_k|\leq \eps |X^H_k|\leq 2\eps n$.
By~\eqref{eq:linked} and the definition of candidates in~\eqref{eq:candidates}, we obtain that $\dg_{Z^H_i(\phi^\ast)}(x_i)\geq \alpha_Z|V_i|$ for all $x_i\in\tX^H_i$.
Note further that $\omega_v(E(\bigcup_{H\in \cH}A^H_i[X_i^H\sm\tX^H_i,V_i]))=\omega(\sX_i)\pm\eps^{1/2}n^2$.\COMMENT{$\omega(\bigcup_{H\in \cH}\tX^H_i)\leq 2\eps\alpha^{-2}n^2\leq{\eps}^{1/2}n^2$.}
Hence, there exists a subgraph $Z^H_i\subseteq Z^H_i(\phi^\ast)$ that satisfies~\ref{exc sr}, \ref{exc triple} and~\ref{exc testers}, which can be seen by keeping each edge in $Z^H_i(\phi^\ast)[X^H_i\sm\tX^H_i,V_i]$ independently at random with probability $\alpha_Z$ and by possibly removing some edges incident to $x_i\in\tX_i^H$ in  $Z^H_i(\phi^\ast)$ deterministically.\COMMENT{Use McDiarmid's inequality to obtain concentration for~\ref{exc testers}.}
In order to see~\ref{exc boundedness}, note that $|\psi_0^{-1}(v_0v_i)|\leq|N_{\cH}(\phi_0^{-1}(v_0))\cap \sX_i|
\leq {\alpha^{-1}}|\phi_0^{-1}(v_0)|\leq {\alpha^{-1}}\eps |\cH|\leq \eps_0n.$
Since the blow-up instance is $(\eps,\alpha)$-linked, we have $\sum_{H\in \cH}|N_H(\phi_0^{-1}(v_0),\phi_0^{-1}(v_0'))\cap X_i^H|\leq \eps|V_i|^{1/2}$ for all distinct $v_0v_i,v_0'v_i\in E(G[V_0,V_i])$, $i\in[r]$, which implies~\ref{exc codegree}.
This completes the proof of the claim.
\endclaimproof

\begin{step}\label{step:induction}
Induction
\end{step}

We inductively prove that the following statement~\ind{t} holds for all $t\in[r]_0$, which will provide a partial packing of $\cH$ into $G_A$. 
\begin{itemize}
\item[\ind{t}.] 
For all $H\in\cH$ and $Z\in\{A,B\}$, there exists a $t$-partial packing $\phi_t\colon \cX_t^{\phi_t}\to \cV_t$ with $|\cX_t^{\phi_t}\cap X^H_i|\geq (1-\eps_{c_i(t)})n$ for all $i\in [t]$, and there exists a candidacy graph $Z^H_i\subseteq Z^H_i(\phi_t)$ (where $\sA_i:=\bigcup_{H\in \cH}A_i^H$) such that 
\begin{enumerate}[label=(\alph*)]
\item\label{Z candidacy superregular} $Z^H_i$ is $(\eps_{c_i(t)},\alpha_Z d_Z^{m_i(t)})$-super-regular for all $i\in[r]\sm[t]$ if $Z=A$ and for all $i\in[r]$ if $Z=B$;\COMMENT{Recall that $m_i(t)=|N_R(i)\cap [t]|$ and $c_i(t)=\max\{c(j)\colon j\in N_{R^2}(i)\cap[t]\}$, $m_i(0)=c_i(0)=0$.}
\item\label{A boundedness} $\Delta_{\psi_t}(\sA_i)\leq (1+\eps_{c_i(t)})\alpha_Ad_A^{m_i(t)}|V_i|$ for all $i\in[r]\sm[t]$;
\item\label{A codegree} $\Delta^c_{\psi_t}(\sA_i)\leq\sqrt{n}$ for all $i\in[r]\sm[t]$;  
\item\label{A triple intersection} $e_{H}(N_{A^H_i}(v_i),N_{A^H_j}(v_j))=(\alpha_A^2 d_A^{m_i(t)+m_j(t)}\pm\eps_{\max\{c_i(t),c_j(t)\}})e_{H}(X^H_i,X^H_j)$ for all $H\in\cH$, $ij\in E(R-[t])$ and $v_iv_j\in E(G_A[V_i,V_j])$;

\item\label{leftover sparse} {$|\phi_t^{-1}(v)|\geq (1-\eps_{c_i(t)}^{1/2})|\cH|-\eps_{c_i(t)}n$} and $|\phi_t^{-1}(v)\cap N_{\cH}(\cX_t\sm\cX_t^{\phi_t})|\leq \eps_{c_i(t)}^{1/2}{n}$ for all $v\in V_i$, $i\in[t]$;

\item \label{updated vertex testers}$\omega_v(E(\sA_i))=\alpha_Ad_A^{m_i(t)}\omega(\sX_{i})\pm\eps_{c_i(t)}n^2$ for all $\omega_v\in\cW_{edge}^i$ and $i\in[r]\sm[t]$;

\item\label{set testers}$|W\cap \bigcap_{j\in [\ell]}\phi_t(Y_j\cap\cX_t^{\phi_t})|= |W||Y_1|\cdots |Y_\ell|/n^\ell \pm \alpha n/2$ for all $(W,Y_1,\ldots,Y_\ell)\in \cW_{set}$ with $W\subseteq V_i$, $i\in[t]$;\COMMENT{Could immediately also use error term `$\pm\eps_{c_i(t)}n$'. However, $\alpha n/2$ is sufficient for Step~\ref{step:completion}.}

\item\label{vertex testers} $\omega(\sX_i\cap \phi_t^{-1}(v))=\omega(\sX_i)/n\pm \alpha n/2$ for all $(v,\omega)\in \cW_{ver}$ with $v\in V_i$, $i\in[t]$.\COMMENT{Could also use an error term in `$\eps_{c_i(t)}$' but would have to be a bit more careful when proving~\ind{t}\ref{vertex testers}. However, $\alpha n/2$ is sufficient for Step~\ref{step:completion}.}

\end{enumerate}
\end{itemize}
{Properties~\ind{t}\ref{Z candidacy superregular}--\ref{A triple intersection} will be used particularly to establish~\ind{t+1} by applying Lemma~\ref{lem:packing lemma}.
Property~\ref{updated vertex testers} enables us to establish~\ref{vertex testers}, which together with~\ref{set testers} basically implies Lemma~\ref{lem:main matching}\ref{item:set testers1} and~\ref{item: vertex tester1} as we merely modify the $r$-partial packing $\phi_r$ for the completion in Step~\ref{step:completion} where we exploit~\ref{Z candidacy superregular} (for $Z=B$) and~\ref{leftover sparse}. 
 }

The statement \ind{0} holds for $\phi_0=\phi^\ast$ by Claim~\ref{claim exceptional}.\COMMENT{Note that \ref{leftover sparse}, \ref{set testers}, \ref{vertex testers} are vacuously true.}
Hence, we assume the truth of~\ind{t} for some $t\in[r-1]_0$ and let $\phi_t\colon \cX_t^{\phi_t}\to \cV_t$ and $A^H_i$ and $B^H_i$ be as in~\ind{t}; we set $\sA_i:=\bigcup_{H\in \cH}A_i^H$ and $\sB_i:=\bigcup_{H\in\cH}B^H_i$.
We will extend $\phi_t$ such that \ind{t+1} holds. 
Any function $\sigma\colon \sX_{t+1}^\sigma\to V_{t+1}$ with $\sX_{t+1}^\sigma\subseteq\sX_{t+1}$ extends $\phi_t$ to a function $\phi_{t+1}\colon \cX_t^{\phi_t}\cup\sX_{t+1}^\sigma\to \cV_{t+1}$ as follows:
\begin{align}\label{eq:phi_{t+1}}
\phi_{t+1}(x):=\begin{cases}
\phi_t(x) & \text{if $x\in \cX_t^{\phi_t}$,}\\
\sigma(x) & \text{if $x\in \sX_{t+1}^\sigma$.}
\end{cases}
\end{align}
We now make a key observation:
By definition of the candidacy graphs $\sA_{t+1}$ and their edge set labellings as in~\eqref{eq:def edge set labelling}, if $\sigma$ is a conflict-free packing in $\sA_{t+1}$ as defined in~\eqref{eq:conflict free}, then $\phi_{t+1}$ is a $(t+1)$-partial packing.

We aim to apply Lemma~\ref{lem:packing lemma} in order to obtain a conflict-free packing $\sigma$ in~$\sA_{t+1}$.
Let 
\begin{align*}
\cH_{t+1} &:=\bigcup_{H\in\cH}H_+\left[\textstyle\bigcup_{i\in N_R[t+1]\sm[t]}X^H_i\cup \bigcup_{i\in N_R(t+1)}X_i^{H,B}\right],
\\G_{t+1} &:=G_+\left[\textstyle\bigcup_{i\in N_R[t+1]\sm[t]}V_i\cup \bigcup_{i\in N_R(t+1)}V_i^{B}\right],
\\ \cA_{t+1}&:=\textstyle\bigcup_{i\in N_R[t+1]\sm[t]}\sA_i \cup \bigcup_{i\in N_R(t+1)}\sB_i,
\\ R_{t+1}&:=R[N_{R}[t+1]\sm[t]]\cup R_B[N_{R_B}[t+1]].
\end{align*}
\COMMENT{Note that the candidacy graphs in $\sB_{t+1}$ are not updated when we find a conflict-free packing in $\sA_{t+1}$, thus we only consider $\bigcup_{i\in N_R(t+1)}\sB_i$.
However, for the reduced graph we consider the graph $R_B[N_{R_B}[t+1]]$, which only consists of edges of the form $\{t+1,j^B\}$.}
Note that $\sP:=(\cH_{t+1},G_{t+1},R_{t+1},\cA_{t+1},\psi_t|_{E(\cA_{t+1})})$ is a packing instance of size $\dg_{R_{t+1}}(t+1)$ with $t+1$ playing the role of $0$, and we claim that $\sP$ is indeed an $(\eps_{c(t+1)-1},\bd)$-packing instance, where $\bd=\big(d_A,d_B,(\alpha_A d_A^{m_i(t)})_{i\in N_R[t+1]\sm[t]}, (\alpha_B d_B^{m_i(t)})_{i\in N_R(t+1)}\big)$.\COMMENT{This notation is intuitively clear why a precise (and ugly) definition is omitted. I did a representative poll at the Ulm institute.}
Observe that by definition of $c_i(t)$ and $m_i(t)$ in~\eqref{eq:def c_i m_i}, we have:
\begin{align}\label{eq:change c_i m_i}
\begin{minipage}[c]{0.9\textwidth}\em
If $i\in N_R(t+1)$, then $m_i(t+1)=m_i(t)+1$, and $c(t+1)=c_i(t+1)>\max\{c_i(t),c_j(t)\}$ for all $j\in N_R(i)$.
If $i\in [r]\sm N_R(t+1)$, then $m_i(t+1)=m_i(t)$.
\end{minipage}\ignorespacesafterend
\end{align}
Note that for the inequality in~\eqref{eq:change c_i m_i} we used that no pair of adjacent vertices in $R$ has two neighbours in $R$ that are coloured alike as we have chosen the vertex colouring as a colouring in $R^3$.
In particular, we infer from~\eqref{eq:change c_i m_i} that $\eps_{c(t+1)-1}=\eps_{c_i(t+1)-1}\geq \eps_{c_i(t)}$ for all $i\in N_R(t+1)$. 
Therefore, \ref{item:P G superreg} follows from~\eqref{eq:G_Z sr}, 
property~\ref{item:P A superreg} follows from~\ind{t}\ref{Z candidacy superregular}, 
property~\ref{item:P triple intersection} follows from~\ind{t}\ref{A triple intersection}  with $R[N_R[t+1]\sm[t]]$ playing the role of $R_A$,\COMMENT{Note that for $ij\in E(R[N_R[t+1]\sm[t]])$, we have $\eps_{c_i(t+1)-1}\geq \eps_{\max\{c_i(t),c_j(t)\}}$.} 
and~\ref{item:P boundedness} follows from~\ind{t}\ref{A boundedness}.

Observe further that
\begin{itemize}
\item $\psi_t$ as defined in~\eqref{eq:def edge set labelling} satisfies $\|\psi_t\|\leq \alpha^{-1}$;
\item $\sum_{H\in\cH}e_{H}(X_{t+1}^H,X^H_i)\leq (1-\alpha)dn^2\leq d_An^2$ for all $i\in N_R(t+1)\sm[t]$;
\item $\Delta^c_{\psi_t}(\sA_i)\leq\sqrt{n}$ for all $i\in N_R[t+1]\sm[t]$ by \ind{t}\ref{A codegree}.
\end{itemize}
\COMMENT{By assumption, also $e_H(X_i^H,X_j^H)\geq\alpha^2 n$ for all $H\in\cH, ij\in E(R)$.}Hence, we can apply Lemma~\ref{lem:packing lemma} to $\sP$ with\COMMENT{$R[N_R[t+1]\sm[t]]$ plays the role of $R_A$, which is stated explicitly when proving~\ind{t+1}\ref{A triple intersection}.}
\begin{center}
\begin{tabular}{c|c|c|c|c|c}
 $n$  & $\eps_{c(t+1)-1}$ & $\eps_{c(t+1)}$  & $\alpha^{-1}$ & $\dg_{R_{t+1}}(t+1)$ &  $R[N_R[t+1]\sm[t]]$
\\ \hline
 $n$  & $\eps$ & $\eps'$ & $s$ & $r$ &  $R_A$
\vphantom{\Big(}
\end{tabular}
\end{center}

\medskip\noindent
\COMMENT{We could also write that we apply it with $\eps_{c_i(t+1)-1}$ and $\eps_{c_i(t+1)}$, cf~\eqref{eq:change c_i m_i}; however index $i$ is unclear.
}and with set testers $\cW_{set}^{t+1}$ where we denote by $\cW_{set}^{t+1}\subseteq\cW_{set}$ the set of set testers $(W,Y_1,\ldots,Y_\ell)$ with $W\subseteq V_{t+1}$, and with edge testers $\cW_{edge}^{t+1}\cup \cW_{edge}^\ast$ where we will define the set $\cW_{edge}^\ast$ when proving~\ind{t+1}\ref{A triple intersection} and \ref{leftover sparse} in Steps~\ref{step:d} and~\ref{step:e}.  

Let $\sigma\colon\sX_{t+1}^\sigma\to V_{t+1}^\sigma$ be the conflict-free packing in $\sA_{t+1}$ obtained from Lemma~\ref{lem:packing lemma} with $|\sX_{t+1}^\sigma\cap X_{t+1}^H|\geq (1-\eps_{c_i(t+1)})n$ for all $H\in\cH$, which extends $\phi_t$ to $\phi_{t+1}$ as defined in~\eqref{eq:phi_{t+1}}.
Fix some $H\in\cH$.
By Definition~\ref{def:updated candidacy}, the updated candidacy graphs with respect to $\sigma$ obtained from Lemma~\ref{lem:packing lemma} are also updated candidacy graphs with respect to $\phi_{t+1}$ as defined in~\eqref{eq:def candidacy} in Step~\ref{step:candidacy}. 
Hence, the graphs $A_i^{H,new}$ in Lemma~\ref{lem:packing lemma} correspond to subgraphs $\tA^H_i\subseteq A^H_i(\phi_{t+1})$ for all $i\in N_R(t+1)\sm[t]$, and $\tB^H_i\subseteq B^H_i(\phi_{t+1})$ for all $i\in N_R(t+1)$ that satisfy~\ref{emb lem:1}--\ref{emb lem:vertex testers}.

\begin{substep}
Checking~\ind{t+1}{\rm \ref{Z candidacy superregular}}
\end{substep}

By~\ref{emb lem:1} and~\eqref{eq:change c_i m_i}, we obtain that $\tA^H_i$ is $(\eps_{c_i(t+1)},\alpha_A d_A^{m_i(t+1)})$-super-regular for all $i\in N_R(t+1)\sm[t]$, and $\tB^H_i$ is $(\eps_{c_i(t+1)},\alpha_B d_B^{m_i(t+1)})$-super-regular for all $i\in N_R(t+1)$.
Note that for each $i\in[r]\sm N_R(t+1)$, we have $m_i(t)=m_i(t+1)$ and $A^H_i(\phi_t)=A^H_i(\phi_{t+1})$ and $B^H_i(\phi_t)=B^H_i(\phi_{t+1})$.
For $i\in[r]\sm[t+1], i'\in[r]$, let
\begin{align}\label{eq:def hat A_i}
\hA^H_i:=\begin{cases}
\tA^H_i & \text{if $i\in N_R(t+1)\sm[t]$,}\\
A^H_i & \text{otherwise.}
\end{cases}
&&
\hB^H_{i'}:=\begin{cases}
\tB^H_{i'} & \text{if $i'\in N_R(t+1)$,}\\
B^H_{i'} & \text{otherwise.}
\end{cases}
\end{align}
Then the graphs $\hA^H_i$ and $\hB^H_{i'}$ are candidacy graphs satisfying \ind{t+1}\ref{Z candidacy superregular}.

\begin{substep}
Checking~\ind{t+1}{\rm \ref{A boundedness}} and \ind{t+1}{\rm\ref{A codegree}}
\end{substep}

The new edge set labelling $\psi_{t+1}$ as defined in~\eqref{eq:def edge set labelling} corresponds to the updated edge set labelling as in Definition~\ref{def:updated labelling}. 
By~\ref{emb lem:boundedness}, \ind{t}\ref{A boundedness} and~\eqref{eq:change c_i m_i}, we obtain for every $i\in[r]\sm[t+1]$ that $\Delta_{\psi_{t+1}}(\bigcup_{H\in\cH}\hA^H_i)\leq (1+\eps_{c_i(t+1)})\alpha_Ad_A^{m_i(t+1)}|V_i|$.
This establishes~\ind{t+1}\ref{A boundedness}.
Similarly, by~\ref{emb lem:codegree} with $R[N_R[t+1]\sm[t]]$ playing the role of $R_A$ and by \ind{t}\ref{A codegree}, we obtain for every $i\in[r]\sm[t+1]$ that $\Delta^c_{\psi_{t+1}}(\bigcup_{H\in\cH}\hA^H_i)\leq \sqrt{n}$, which establishes~\ind{t+1}\ref{A codegree}.

\begin{substep}\label{step:d}
Checking~\ind{t+1}{\rm \ref{A triple intersection}}
\end{substep}

In order to show~\ind{t+1}\ref{A triple intersection}, fix $H\in\cH$, $ij\in E(R-[t+1])$, and $v_iv_j\in E(G[V_i,V_j])$. 
Observe, that $|\{i,j\}\cap N_R(t+1)|\in\{0,1,2\}.$

If $|\{i,j\}\cap N_R(t+1)|=2$, then this implies together with~\eqref{eq:change c_i m_i}  that $c_i(t+1)=c_j(t+1)=c(t+1)>\max\{c_i(t),c_j(t)\}$, and $m_i(t)+1=m_i(t+1)$ as well as $m_j(t)+1=m_j(t+1)$. 
Hence we obtain by~\ref{emb lem:triple intersection} (with $R[N_R[t+1]\sm[t]]$ playing the role of $R_A$) and~\ind{t}\ref{A triple intersection} that
\begin{align}\label{eq:size hA triple intersection}
e_{H}(N_{\hA^H_i}(v_i),N_{\hA^H_j}(v_j))
=\big(\alpha_A^2 d_A^{m_i(t+1)+m_j(t+1)}\pm\eps_{\max\{c_i(t+1),c_j(t+1)\}}\big)e_{H}(X^H_i,X^H_j).
\end{align}

If $|\{i,j\}\cap N_R(t+1)|=1$, say $i\in N_R(t+1)$, then this implies together with~\eqref{eq:change c_i m_i} that
\begin{align}\label{eq:ij=1}
\begin{minipage}[c]{0.7\textwidth}\em
 $c_i(t+1)=\max\{c_i(t+1),c_j(t+1)\}=c(t+1)>\max\{c_i(t),c_j(t)\}$, and $m_i(t)+1=m_i(t+1)$, $m_j(t)=m_j(t+1)$. 
\end{minipage}\ignorespacesafterend
\end{align}
By~\eqref{eq:def hat A_i}, we have that $\hA_j^H=A_j^H$ because $j\notin N_R(t+1)$.
Let $N:=N_{A^H_i}(v_i)\cap N_{H}(N_{A^H_j}(v_j))$ and we define a weight function $\omega_N\colon E(\sA_i)\to\{0,1\}$ by $\omega_N(xv):=\mathbbm{1}_{\{v=v_i \}}\mathbbm{1}_{\{x\in N \}}$ and add $\omega_N$ to~$\cW_{edge}^\ast$.
Note that $\dim(\omega_N)=1$ (with $\dim(\omega_N)$ defined as in~\eqref{eq:def dim(w)}) and that $\omega_N(E(\sA_i))=|N|=(\alpha_A^2 d_A^{m_i(t)}d_A^{m_j(t)}\pm\eps_{\max\{c_i(t),c_j(t)\}})e_{H}(X^H_i,X^H_j)$ by~\ind{t}\ref{A triple intersection}.
This implies that
\begin{align*}
e_{H}(N_{\hA^H_i}(v_i),N_{\hA^H_j}(v_j))
&\stackrel{\hphantom{\ref{emb lem:edge weight}}}{=}
|N_{\hA^H_i}(v_i)\cap N_{H}(N_{A^H_j}(v_j))|
= |N_{\hA^H_i}(v_i)\cap N| 
= \omega_N(E(\hA^H_i))
\\&\stackrel{\ref{emb lem:edge weight}}{=}
(1\pm\eps_{c(t+1)}^2)d_A \omega_N(E(\sA_i))\pm\eps_{c(t+1)}^2n
\\&\stackrel{~\eqref{eq:ij=1}~}{=}
\big(\alpha_A^2 d_A^{m_i(t+1)+m_j(t+1)}\pm\eps_{\max\{c_i(t+1),c_j(t+1)\}}\big)e_{H}(X^H_i,X^H_j).
\end{align*}
\COMMENT{Note that $e_{H}(X^H_i,X^H_j)\geq \alpha^2n$.}

If $|\{i,j\}\cap N_R(t+1)|=0$, then this implies together with~\eqref{eq:change c_i m_i} and~\eqref{eq:def hat A_i}, that $m_i(t)=m_i(t+1)$, $m_j(t)=m_j(t+1)$, and $\hA^H_i=A^H_i$. Consequently,~\eqref{eq:size hA triple intersection} holds which establishes~\ind{t+1}\ref{A triple intersection}.

\begin{substep}\label{step:e}
Checking~\ind{t+1}{\rm \ref{leftover sparse}}
\end{substep}

In order to establish~\ind{t+1}\ref{leftover sparse}, {we first consider $v_i\in V_i$ for $i\in N_R(t+1)\cap[t]$.
We define a weight function $\omega_i^\ast\colon E(\sA_{t+1})\to\{0,1\}$ by $\omega_{v_{i}}^\ast(xv):=\mathbbm{1}_{\{ x\in \sX_{i}^\ast\}}$ for $\sX_{i}^\ast:=N_{\cH}(\phi_t^{-1}(v_i))\cap\sX_{t+1}$ and every $xv\in E(\sA_{t+1})$, and we add $\omega_{v_{i}}$ to~$\cW_{edge}^\ast$.
By~\ind{t}\ref{Z candidacy superregular}, we have
\begin{align}\label{eq:weight w_v_i}
\omega_{v_{i}}(E(\sA_{t+1}))=(\alpha_Ad_A^{m_{t+1}(t)}\pm3\eps_{c_{t+1}(t)})|\sX_i^\ast|n,
\end{align}  and by~\ind{t}\ref{leftover sparse}, we have 
\begin{align}\label{eq:St e bound}
|\phi_{t+1}^{-1}(v_i)\cap N_{\cH}(\cX_{t+1}\sm\cX_{t+1}^{\phi_{t+1}})|\leq \eps_{c_i(t)}^{1/2}{n} + |\sX_{i}^\ast|-\omega_{v_{i}}(M)
\end{align}
with $M=M(\sigma)$ being the corresponding edge set to $\sigma$.
By~\ref{emb lem:vertex testers}, we obtain that
\begin{align*}
\omega_{v_{i}}(M)&
\stackrel{\ref{emb lem:vertex testers}}{=}
(1\pm\eps_{c(t+1)})\frac{\omega_{v_{i}}(E(\sA_{t+1}))}{\alpha_Ad_A^{m_{t+1}(t)}n}\pm \eps_{c(t+1)}n
\stackrel{\eqref{eq:weight w_v_i}}{\geq}
(1-\eps_{c(t+1)}^{3/4})|\sX_i^\ast|-{\eps_{c(t+1)}n}.
\end{align*}
Together with~\eqref{eq:St e bound}, this implies that $|\phi_{t+1}^{-1}(v_i)\cap N_{\cH}(\cX_{t+1}\sm\cX_{t+1}^{\phi_{t+1}})|\leq \eps_{c(t+1)}^{1/2}{n}.$
}

{Hence, it now suffices to establish~\ind{t+1}\ref{leftover sparse}} for all $v_{t+1}\in V_{t+1}$ by~\ind{t}\ref{leftover sparse}.
We define weight functions $\omega_{v_{t+1}},\omega_{v_{t+1}}^\ast\colon E(\sA_{t+1})\to\{0,1\}$ by
$\omega_{v_{t+1}}(xv):=\mathbbm{1}_{\{v=v_{t+1} \}}$ and $\omega_{v_{t+1}}^\ast(xv):=\mathbbm{1}_{\{v=v_{t+1} \text{ and } x\in \sX_{t+1}^\ast\}}$ for $\sX_{t+1}^\ast:=N_{\cH}(\cX_t\sm\cX_t^{\phi_t})\cap\sX_{t+1} $ and every $xv\in E(\sA_{t+1})$, and we add $\omega_{v_{t+1}}$ and $\omega_{v_{t+1}}^\ast$ to~$\cW_{edge}^\ast$. Observe that~\ind{t} implies that
\begin{align}\label{eq:weight v_t+1}
\omega_{v_{t+1}}(E(\sA_{t+1}))&=(\alpha_Ad_A^{m_{t+1}(t)}\pm3\eps_{c_{t+1}(t)})|\cH|n,
\\\label{eq:weight v_t+1^*}
\omega_{v_{t+1}}^\ast(E(\sA_{t+1}))&\leq|\sX_{t+1}^\ast|\leq \eps_{c_{t+1}(t)}^{1/2} |\cH|n,
\end{align}
and we have that $|\phi_{t+1}^{-1}(v_{t+1})|=\omega_{v_{t+1}}(M)$ and $|\phi_{t+1}^{-1}(v_{t+1})\cap \sX_{t+1}^\ast|\leq \omega_{v_{t+1}}^\ast(M)$.
By~\ref{emb lem:vertex testers}, we obtain that
\begin{align*}
\omega_{v_{t+1}}(M)&
\stackrel{\ref{emb lem:vertex testers}}{=}
(1\pm\eps_{c(t+1)})\frac{\omega_{v_{t+1}}(E(\sA_{t+1}))}{\alpha_Ad_A^{m_{t+1}(t)}n}\pm \eps_{c(t+1)}n
\stackrel{\eqref{eq:weight v_t+1}}{\geq}
(1-\eps_{c(t+1)}^{1/2})|\cH|-{\eps_{c(t+1)}n};
\\
\omega_{v_{t+1}}^\ast(M)&
\stackrel{\ref{emb lem:vertex testers}}{=}
(1\pm\eps_{c(t+1)})\frac{\omega_{v_{t+1}}^\ast(E(\sA_{t+1}))}{\alpha_Ad_A^{m_{t+1}(t)}n}\pm \eps_{c(t+1)}n
\stackrel{\eqref{eq:weight v_t+1^*}}{\leq}
\eps_{c(t+1)}^{1/2}{n}.
\end{align*}
Note that $\eps_{c(t+1)}=\eps_{c_{t+1}(t+1)}$.
{Altogether, this} establishes~\ind{t+1}\ref{leftover sparse}.

\begin{substep}
Checking~\ind{t+1}{\rm \ref{updated vertex testers}--\ref{vertex testers}}
\end{substep}

In order to establish~\ind{t+1}\ref{updated vertex testers}, consider $\omega_v\in\cW_{edge}^i$ for $i\in N_R(t+1)\sm[t+1]$. 
By~\eqref{eq:change c_i m_i}, it holds that $c(t+1)=c_i(t+1)$.
With~\ref{emb lem:edge weight} we obtain that 
\begin{align*}
\omega_v(E(\textstyle\bigcup_{H\in \cH}\hA^H_i))
=(1\pm \eps_{c(t+1)}^2)d_A\omega_v(E(\sA_i))\pm \eps_{c(t+1)}^2n^2
\stackrel{\text{\ind{t}\ref{updated vertex testers}}}{=}
\alpha_A d_A^{m_i(t+1)}\omega(\sX_i)\pm\eps_{c_i(t+1)}n^2,
\end{align*}
\COMMENT{
\begin{align*}
&(1\pm \eps_{c_i(t+1)}^2)d_A\omega_v(E(\sA_i))\pm \eps_{c_i(t+1)}^2n^2
\stackrel{\text{\ind{t}\ref{updated vertex testers}}}{=}(1\pm \eps_{c_i(t+1)}^2)d_A \Big(\alpha_Ad_A^{m_i(t)}\omega(\sX_i)\pm\eps_{c_i(t)}n^2\Big)\pm \eps_{c_i(t+1)}^2n^2
\\&=\alpha_A d_A^{m_i(t+1)}\omega(\sX_i)\pm\eps_{c_i(t+1)}^2\omega(\sX_i)\pm 2\eps_{c_i(t+1)}^2n^2
=\alpha_Ad_A^{m_i(t+1)}\omega(\sX_i)\pm\eps_{c_i(t+1)}n^2,
\end{align*}
where we used that $\omega(\sX_{i})\leq\alpha^{-2}n^2$.}which together with~\ind{t}\ref{updated vertex testers} establishes~\ind{t+1}\ref{updated vertex testers}.

Next we verify~\ind{t+1}\ref{set testers}. Note that \ref{emb lem:set testers} implies that $|W\cap \bigcap_{j\in [\ell]}\sigma(Y_j\cap\sX_{t+1}^\sigma)|= |W||Y_1|\cdots |Y_\ell|/n^\ell \pm \eps_{c(t+1)}n$ for all $(W,Y_1,\ldots,Y_\ell)\in \cW_{set}^{t+1}$, which together with~\ind{t}\ref{set testers} yields~\ind{t+1}\ref{set testers}.

In order to establish~\ind{t+1}\ref{vertex testers}, let $\cW_{ver}^{t+1}\subseteq\cW_{ver}$ be the set of vertex testers $(v,\omega)$ with $v\in V_{t+1}$.
Hence, for all $(v,\omega)\in\cW_{ver}^{t+1}$ and and its corresponding edge tester $\omega_v\in\cW_{edge}^{t+1}$ as defined in~\eqref{eq:def W_edge}, property~\ref{emb lem:vertex testers} implies that 
\begin{align*}
\omega(\sX_{t+1}\cap \sigma^{-1}(v))
= \omega_v(M)
&
\stackrel{\ref{emb lem:vertex testers}}{=}
(1\pm\eps_{c(t+1)})\frac{\omega_v(E(\sA_{t+1}))}{\alpha_Ad_A^{m_{t+1}(t)}n} \pm\eps_{c(t+1)}n
\\&\stackrel{\text{~\ind{t}\ref{updated vertex testers}~}}{=}
(1\pm\eps_{c(t+1)})\frac{\alpha_Ad_A^{m_{t+1}(t)}\omega(\sX_{t+1})\pm\eps_{c_{t+1}(t)}n^2}{\alpha_Ad_A^{m_{t+1}(t)}n}\pm\eps_{c(t+1)}n
\\&\stackrel{\hphantom{\text{~\ind{t}\ref{updated vertex testers}~}}}{=}
\frac{\omega(\sX_{t+1})}{n}\pm\alpha n/2.
\end{align*}
\COMMENT{For the last equality we used that $\omega(\sX_{t+1})\leq 2\alpha^{-2}n^2$ and $\eps_{c(t+1)}\ll\alpha_A,d_A$.}Together with~\ind{t}\ref{vertex testers}, this yields~\ind{t+1}\ref{vertex testers}.

\begin{step}\label{step:completion}
Completion
\end{step}

Let $\phi_r\colon \textstyle{\bigcup_{H\in\cH,i\in[r]_0}\hX^H_i \to\cV_r}$ be an $r$-partial packing satisfying~\ind{r} with $(\eps_T,d_i)$-super-regular candidacy graphs $B^H_i\subseteq B^H_i(\phi_r)$ where $d_i:=\alpha_Bd_B^{\dg_R(i)}$ for all $i\in[r]$. 
We aim to apply iteratively the ordinary blow-up lemma in order to complete the partial packing~$\phi_r$ using the edges in~$G_B$.
Recall that $\eps_T\ll\mu\ll\gamma\ll\alpha,d$.
Our general strategy is as follows.
For every $H\in\cH$ in turn,
we choose a set $X_i\In X^H_i$ for all $i\in [r]$ of size roughly $\mu n$ by selecting every vertex uniformly at random with the appropriate probability and adding $X^H_i\setminus \hX_i^{H}$ deterministically.
Afterwards, we apply the blow-up lemma to embed $H[{X_1\cup\ldots\cup X_r}]$ into $G_B$, which together with~$\phi_r$ yields a complete embedding of~$H$ into $G_A\cup G_B$.
Before we proceed with the details of our procedure (see Claim~\ref{claim:completion_last}), we {verify in Claim~\ref{claim:completion} that we can indeed apply the blow-up lemma to a subgraph of $H\in\cH$ provided some easily verifiable conditions are satisfied. }

Recall that we have defined the candidacy graph $B^H_i$ on a copy $(X^{H,B}_i,V_i^B)$ via the bijection~$\pi$ only to conveniently apply Lemma~\ref{lem:packing lemma} in Step~\ref{step:induction}.
That is, for all $H\in\cH, i\in[r]$, we can identify~$B^H_i$ with an isomorphic bipartite graph on $(X^H_i,V_i)$ and edge set $\{xv\colon \pi(x)\pi(v)\in E(B^H_i)\}$. 
Let $\cB$ be the union over all $H\in\cH, i\in[r]$ of these graphs.
For $H\in\cH$ and a subgraph $G^\circ$ of~$G_B$,
we write $\cB_{G^\circ}[{X_i^H,V_i}]$ for the graph that arises from $\cB[{X_i^H,V_i}]$
by deleting every edge $xv$ with $x\in X_i^H, v\in V_i$ for which there exists $x'x\in E(H)$ such that $\phi_r(x')v\in E(G^\circ)$.
(We may think of $E(G^\circ)$ as the edge set in $G_B$ that we have already used in our completion step for packing some other graphs of~$\cH$ into $G_A\cup G_B$.)%
\COMMENT{We could also write that $x'\in \widehat{X}_{j}^H\cup  X_0^H$. However, it is sufficient to say that $x'\in \widehat{X}_{j}^H$ since the candidacy graphs in $\cB$ are already updated with respect to the embedding $\phi_0$ between $X_0^H$ and $V_0$. 
That is, an edge $xv\in E(\cB[\tX^H_i,\tV^H_i])$ with $x\in \tX^H_i$, $x'x\in E(H)$, $x'\in \widehat{X}_{j}^H\cup  X_0^H$ and $\phi_r(x')v\in E(G^\circ)$ does not exist.}
For future reference, we observe that
\begin{align}\label{eq:circ}
\begin{minipage}[c]{0.6\textwidth}\em
$\Delta(\cB[X^H_i,V_i]-\cB_{G^\circ}[X_i^H,V_i])\leq \alpha^{-1}\Delta(G^\circ)$ {for all $i\in[r]$.}
\end{minipage}\ignorespacesafterend
\end{align}

\begin{claim}\label{claim:completion}
Suppose $H\in\cH$, $G^\circ \In G_B$ and {$W_i\In V_i$} for all $i\in[r]$ such that the following hold:
\begin{enumerate}[label=\rm (\alph*)]
	\item\label{item:size X} 
{$V_i\sm\phi_r(\hX_i^H)\subseteq W_i\subseteq V_i$ and $|W_i|= (\mu \pm \eps_T^{1/2}) n$ for all $i\in [r]$;}
	\item\label{item:reg B} $\cB_{G^\circ}[X_i, W_i]$ is $(\mu^{1/31},d_i)$-super-regular where {$X_i:=(X_i^H\cap \phi^{-1}_r(W_i\cap \phi_r(\hX_i^{H})))\cup (X_i^{H}\setminus\hX_i^{H})$;}
	\item\label{item:reg G} $G_B[W_i,W_j]$ is $(\eps_T^{1/3},d_B)$-super-regular for all $ij\in E(R)$;
	\item\label{item:max G circ} $|N_{G^\circ}(v)\cap W_i|\leq \mu^{3/2}n$ for all $v\in V_j$ and $ij\in E(R)$.
	\end{enumerate}
Then there exists an embedding $\phi^{\tH}$ of {$\tH:=H[X_1,\ldots,X_r]$ into $\tG:=G_B[W_1,\ldots,W_r]-G^\circ$ such that $\phi':=\phi^{\tH} \cup \phi_r|_{V(H)\sm V(\tH)}$} is an embedding of $H$ into $G$
where {$\phi'(X_i^H)= V_i$ for all $i\in[r]$ and} all edges incident to a vertex in $\bigcup_{i\in[r]}X_i$ are embedded on an edge in $G_B-G^\circ$.
\end{claim}
\claimproof
Note that $|X_i|=|W_i|= (\mu\pm \eps_T^{1/2})n$ (by~\ref{item:size X} and~\ref{item:reg B})
and $\tH[X_i,X_j]$ is empty whenever $ij\notin E(R)$ and a matching otherwise.
Moreover, by~\ref{item:reg G} and~\ref{item:max G circ}, Fact~\ref{fact:regularity robust} yields that $\tG[W_i,W_j]$ is $(\mu^{1/31},d_B)$-super-regular for all $ij\in E(R)$ (with room to spare).
This shows that $(\tH, \tG,R,(X_i)_{i\in[r]},(W_i)_{i\in[r]})$ is a $(\mu^{1/31},d_B)$-super-regular blow-up instance.
We apply Theorem~\ref{thm:KSS} to this blow-up instance (with $(\cB_{G^\circ}[X_i, W_i])_{i\in[r]}$ playing the role of $(A_i)_{i\in[r]}$) and obtain an embedding $\phi^{\tH}$ of $\tH$ into $\tG$.
Recall that $\phi_r|_{V(H)\sm V(\tH)}$ is an embedding of $H-V(\tH)$ into $G_A$.
For $x\in V(\tH)$, $x'\in V(H)\sm V(\tH)$ and $xx'\in E(H)$,
we conclude that $\phi^{\tH}(x)\phi_r(x')\in E(G_B-G^\circ)$ by the definition of the candidacy graphs in~\eqref{eq:def candidacy} and the definition of $\cB_{G^\circ}$.
\endclaimproof

\begin{claim}\label{claim:completion_last}
For all $H\in\cH$ and $i\in[r]$, there exist sets $\oX_i^H\subseteq \hX_i^H$ with $|\oX_i^H|\geq (1-2\mu)n$ and a packing~$\phi$ of~$\cH$ into~$G$ {that extends $\phi^\ast$} such that $\phi|_{\oX_i^H}=\phi_r|_{\oX_i^H}$ as well as $\phi(X_i^H)=V_i$.
\end{claim}

\claimproof
We write $\cH=\{H_1,\ldots,H_{|\cH|}\}$ and let $\cH_h:=\{H_1,\ldots,H_h\}$ for all $h \in [|\cH|]_0$.
For all $v\in V_i$, $i\in[r]$, 
let $\tau^h(v):=|\{H\in \cH_h\colon v\in V_i\setminus\phi_r(\hX_i^{H})\}|$ and 
$\sigma^{h}(v):=|\{\phi_r^{-1}(v)\cap \bigcup_{H\in \cH_h, j\in N_R(i)}N_{H}(X_j^H\setminus\hX_j^H)\}|$.
We inductively prove that the following statement \indC{h} holds for all $h\in[|\cH|]_0$.

\begin{itemize}
\item[\indC{h}.]
There exists a packing $\phi^{h}$ of $\cH_h$ into $G$ {that extends $\phi^\ast$}
such that for $G^\circ_{h}:=G_B \cap \phi^h(\cH_h)$ we have
\begin{enumerate}[label=\rm (\Alph*)]
\item\label{item:A} $\deg_{G^\circ_{h}}(v)\leq \alpha^{-1}(\tau^{h}(v)+\sigma^{h}(v))+\mu^{2/3}n$ for all $v\in V(G)$;
\item\label{item:B} for all $H\in\cH_h$ and $i\in[r]$, there exist sets $\oX_i^H\subseteq \hX_i^H$ with $|\oX_i^H|\geq (1-2\mu)n$ such that $\phi^{h}|_{\oX_i^H}=\phi_r|_{\oX_i^H}$ as well as $\phi^h(X_i^H)=V_i$.
\end{enumerate}
\end{itemize}

Let $G_0^\circ$ be the edgeless graph on $V(G)$ and $\phi^0$ be the empty function;
then \indC{0} holds.
Hence, we may assume the truth of \indC{h} for some $h\in[|\cH|-1]_0$ and let $\phi^h$ and $G^\circ_h$ be as in~\indC{h}.

By~\indC{h}\ref{item:B}, there are at most $\sum_{H\in \cH_h,j\in N_R[i]}\alpha^{-1}|X_j^H\setminus \oX_j^H|\leq 5\alpha^{-3}\mu n^2$ edges of $G_h^\circ$ incident to a vertex in~$V_i$ for each $i\in[r]$.
Hence, there are at most $10\alpha^{-3}\mu^{1/3}n$ vertices in $V_i$ of degree at least $\mu^{2/3} n/2$ in $G_h^\circ$.
Let $V_i^{high}\In V_i$ be a set of size $\mu^{1/4}n$ that contains all vertices of degree at least $\mu^{2/3} n/2$ in $G_h^\circ$.
For all $i\in [r]$,
we select every vertex in $\phi_r(\hX^{H_{h+1}}_i)\setminus V_i^{high}$ independently with probability~$\mu(1-\mu^{1/4})^{-1}$
and denote by $W_i$ their union together with $V_i\setminus \phi_r(\hX^{H_{h+1}}_i)$;
we define $X_i:=(X_i^{H_{h+1}}\cap \phi^{-1}_r(W_i\cap \phi_r(\hX_i^{H_{h+1}})))\cup (X_i^{H_{h+1}}\setminus \hX_i^{H_{h+1}})$. 
{Note that~\ind{r} implies that}
\begin{align}\label{eq:size leftover}
\begin{minipage}[c]{0.6\textwidth}\em
{$\left|X_i^{H_{h+1}}\sm\hX_i^{H_{h+1}}\right|
=\left|V_i\sm\phi_r(\hX_i^{H_{h+1}})\right|\leq 2\eps_Tn$ {for all $i\in[r]$.}}
\end{minipage}\ignorespacesafterend
\end{align}

In the following we will show that the assumptions of Claim~\ref{claim:completion} are satisfied with probability at least 1/2, say, {for $H_{h+1}$ and $G^\circ_h$ playing the role of $H$ and $G^\circ$, respectively}.
In particular, there is a choice for $W_i$ such that the assumptions of Claim~\ref{claim:completion} hold.

To obtain~\ref{item:size X},
we apply Chernoff's inequality to the sum of indicator variables which indicate whether a vertex in $V_i$ is randomly selected.
{Together with~\eqref{eq:size leftover},} this shows that~\ref{item:size X} holds with probability at least $1-1/n^3$, say.

By{~\indC{h}\ref{item:A} and}~\ind{r}\ref{leftover sparse}, we obtain that $\Delta(G_h^\circ)\leq 2\mu^{2/3}n$.
We exploit~\eqref{eq:circ} and conclude that $\Delta(\cB[X^{H_{h+1}}_i, V_i]-\cB_{G_{h}^\circ}[X^{H_{h+1}}_i, V_i])\leq 2\alpha^{-1}\mu^{2/3}n$ {for all $i\in[r]$}.
Thus Fact~\ref{fact:regularity robust} implies that $\cB_{G_{h}^\circ}[X^{H_{h+1}}_i, V_i]$
is $(\mu^{1/5},d_i)$-super-regular {for all $i\in[r]$}.
For {all $i\in[r]$ and} $x\in X^{H_{h+1}}_i$,
Chernoff's inequality implies that $|N_{\cB_{G_{h}^\circ}}(x)\cap W_i|= (d_i\pm 2\mu^{1/5})|W_i|$
and similarly, for all $v\in V_i$, we have $|N_{\cB_{G_{h}^\circ}}(v)\cap X_i|= (d_i\pm 2\mu^{1/5})|X_i|$.
Moreover, for all distinct $v,v'$ with $|N_{\cB_{G_{h}^\circ}}(v,v')|=(d_i\pm 2\mu^{1/5})^2n$ (which we call \defn{good}, and there are at least $(1-2\mu^{1/5})\binom{n}{2}$ good pairs),
we also obtain $|N_{\cB_{G_{h}^\circ}}(v,v')\cap X_i|=(d_i\pm 3\mu^{1/5})^2|X_i|$, all with probability at least $1-1/n^3$.
Observe that Theorem~\ref{thm:McDiarmid} implies that there at least $(1-3\mu^{1/5})\binom{\mu n}{2}$ good pairs in $W_i$ also with probability at least $1-1/n^3$.
Therefore, we may apply Theorem~\ref{thm: almost quasirandom} and obtain that $\cB_{G_{h}^\circ}[X_i, W_i]$ is $(\mu^{1/31},d_i)$-super-regular {for all $i\in[r]$}
which yields~\ref{item:reg B} with probability at least $1-1/n^2$, say.

To obtain~\ref{item:reg G}, for each $ij\in E(R)$, we proceed as follows.
Observe first that $G_B[V_i,V_j]$ is $(2\eps,d_B)$-super-regular by~\eqref{eq:G_Z sr}.
Hence $G_B[W_i,W_j]$ is clearly $\eps_T^{1/3}$-regular as $|W_i|,|W_j|\geq \mu n/2$ by~\ref{item:size X}.
Therefore, we only need to control the degrees of the vertices in  $G_B[W_i,W_j]$ {which follows directly} by Chernoff's inequality with probability at least $1-1/n^3$ {and because of~\eqref{eq:size leftover}.}

Since $\Delta(G_h^\circ)\leq 2\mu^{2/3}n$ {and because of~\eqref{eq:size leftover}},
we conclude by Chernoff's inequality that~\ref{item:max G circ} holds with probability at least $1-1/n^3$.

Therefore, the assumptions of Claim~\ref{claim:completion} are achieved by our construction with probability at least~$1/2$. 
Fix such a choice for $W_1,\ldots,W_r$ and apply Claim~\ref{claim:completion} that returns an  embedding $\phi^{\tH_{h+1}}$
of $\tH_{h+1}:=H_{h+1}[X_1,\ldots,X_r]$ into $G_B[W_1,\ldots,W_r]-G_h^\circ$
such that $\phi':=\phi^{\tH_{h+1}}\cup \phi_r|_{{V(H_{h+1})\sm V(\tH_{h+1})}}$ is an embedding of~$H_{h+1}$ into~$G$ where {$\phi'(X_i^H)=V_i$} and all edges incident to a vertex in $\bigcup_{i\in[r]}X_i$ are embedded on an edge in $G_B-G^\circ_h$.
We define $\phi^{h+1}:=\phi^h\cup \phi'$ and obtain~\indC{h+1}\ref{item:B} with $\oX_i^{H_{h+1}}:=X_i^{H_{h+1}}\setminus X_i$.
It is straightforward to check that by our construction also~\indC{h+1}\ref{item:A} holds.
\endclaimproof

Let $\phi$ be as in Claim~\ref{claim:completion_last}.
This directly implies conclusion~\ref{item:partition} of Lemma~\ref{lem:main matching}.
Conclusion~\ref{item:set testers1} of Lemma~\ref{lem:main matching} follows from Claim~\ref{claim:completion_last} together with~\ind{r}\ref{set testers}
as we merely modified~$\phi_r$ to obtain~$\phi$.
For a similar reason, Lemma~\ref{lem:main matching}~\ref{item: vertex tester} follows from Claim~\ref{claim:completion_last} and~\ind{r}\ref{vertex testers}.
This completes the proof.
\endproof

\section{Applications}

In what follows, we provide an illustration for an application of vertex and set testers
so that the leftover is suitably well-behaved.
For a graph $H$, let $H^+$ arise from $H$ by adding a labelled vertex~$x$ and joining~$x$ to all vertices of $H$.
We call~$x$ the \emph{apex vertex} of $H^+$.

\begin{theorem}[Keevash~{\cite[cf.~Theorem~7.8]{keevash:18b}}, cf.~{\cite[Theorem~3.6 and Corollary~3.7]{GJKKO:18}}]\label{thm:keevash}
	Suppose $1/n\ll \eps \ll 1/s \ll d_0, 1/m$.
	Suppose $H$ is an $r$-regular graph on $m$ vertices.
	Let $G$ be a graph with vertex partition $(V,W)$ such that $W$ is an independent set, $d_0n\leq |W|\leq |V|=n$ and 
	$|\bigcap_{x\in V'\cup W'}N_G(x)|=(1\pm \eps)d_V^{|V'|}d_W^{|W'|}n$ for all $V'\In V,W'\In W$ with $1\leq|V'|+|W'|\leq s$ where $d_V=rd|W|/n$ and $d_W=d$ for some $d\geq d_0$.
	Suppose that $|N_G(v)\cap V|=r|N_G(v)\cap W|$ for all $v\in V$ and $m$ divides $\dg_G(w)$ for all $w\in W$.
	Then there is a decomposition of the edge set of $G$ into copies of~$H^+$ where the apex vertices are contained in $W$.
\end{theorem}

\lateproof{Theorem~\ref{thm:application}}
Choose $\eps \ll  \delta \ll 1/s \ll \beta \ll \alpha$.
Among the $\alpha n$ graphs in $\cH$ that contain at least $\alpha n$ vertices in components of size at most $\alpha^{-1}$,
there is a collection $\cH'$ of $\beta n$ $r$-regular graphs that contain each at least $\beta n$ vertices in components all isomorphic to some graph $J$
where ${|V(J)|}\leq \alpha^{-1}$ and $r\in [\alpha^{-1}]$.

For all $H\in \cH'$, let $H^-$ arise from $H$ by deleting $\beta n/|V(J)|$ components isomorphic to $J$. 
We denote by $I^H$ a set of $\beta n$ isolated vertices disjoint from {$V(H^-)$}.
Let $\widetilde{\cH}:=(\cH\setminus \cH')\cup \bigcup_{H\in \cH'}(H^-\cup I^H)$.
Let~$G_1$ be a $({2\eps}, s,d_1)$-typical subgraph of~$G$ where $d_1:=(1+\delta)(d-\beta^2r)$
such that $G-G_1$ is $({2\eps}, s,d-d_1)$-typical; that is, $e(\widetilde{\cH}) \leq (1-\delta/2)e(G_1)$.
Clearly, $G_1$ exists by considering a random subgraph and then applying Chernoff's inequality. 
Now we apply Theorem~\ref{thm:quasirandom} to obtain a packing~$\phi$ of $\widetilde{\cH}$ in~$G_1$ with  $\delta^2$ playing the role of~$\alpha$ {and sets of set and vertex testers $\cW_{set}$, $\cW_{ver}$ defined as follows}.
{For all $H_1,\ldots,H_{\ell_1}\in \cH'$ and $v_1,\ldots,v_{\ell_2}\in V(G)$ with $1\leq \ell_1$ and $\ell_1+\ell_2\leq s$, we add the set tester $(V',I^{H_1},\ldots,I^{H_{\ell_1}})$ to $\cW_{set}$ where $V':=V(G)\cap\bigcap_{i\in [\ell_2]}N_{G-G_1}(v_i)$.}
Then Theorem~\ref{thm:quasirandom} implies that (where $d_2:=d-d_1$)
\begin{align}\label{eq:typical application}
\bigg|\bigcap_{i\in [\ell_1]}\phi(I^{H_i})\cap \bigcap_{i\in [\ell_2]}N_{G-G_1}(v_i)\bigg|= (\beta^{\ell_1}d_2^{\ell_2}\pm 2\delta^2)n.
\end{align}
For each $v\in V(G)$, we define a vertex tester $(v,\omega)$ where $\omega$ assigns every vertex in $V(H)$ its degree for all $H\in \widetilde{H}$ {and add $(v,\omega)$ to~$\cW_{ver}$}. 
Then Theorem~\ref{thm:quasirandom} implies that $\Delta(G_1-\phi(\widetilde{\cH}))\leq 2\delta n$.

{Let $G_2:= G- \phi(\widetilde{\cH})$.}
Next, we add $\beta n$ vertices $W:=\{v_H\}_{H\in \cH'}$ to $G_2$ and join $v_H$ to all vertices in $\phi(I^H)$ and denote this new graph by $G_3$.
Let $d_V:=\beta^2r$ and $d_W:=\beta$.
Hence~\eqref{eq:typical application}, the typicality of $G-G_1$, and $\Delta(G_1-\phi(\widetilde{\cH}))\leq 2\delta n$ imply that for all $w_1,\ldots,w_{\ell_1}\in W$ and $v_1,\ldots,v_{\ell_2}\in V(G)$
\begin{align*}
\bigcap_{x\in \Set{w_1,\ldots,w_{\ell_1},v_1,\ldots,v_{\ell_2}}}N_{G_3}(x)
=(1\pm \sqrt{\delta})\beta^{\ell_1}\cdot (\beta^2 r)^{\ell_2}n
=(1\pm \sqrt{\delta})d_W^{\ell_1}d_V^{\ell_2}n
\end{align*}
whenever $1\leq \ell_1+\ell_2\leq s$.
We apply Theorem~\ref{thm:keevash} to $G_3$ to obtain a decomposition of $G_3$ into copies of $J^+$ where the apex vertices are contained in $W$.
Observe that this yields the desired decomposition of $G$ into~$\cH$.
Indeed, for $H\in \cH'$, let $\cJ_H$ be the set of all copies of $J^+$ in $G_3$ whose apex vertex is~$v_H$; hence $|\cJ_H|=|I^H|/|V(J)|=\beta n/|V(J)|$.
We define a packing $\phi'$ of $\cH$ in $G$ as follows.
For all $H\in \cH\setminus  \cH'$, let $\phi'\mid_{V(H)}:=\phi\mid_{V(H)}$.
For all $H\in \cH'$, let $\phi'\mid_{V(H^-)}:=\phi\mid_{V(H^-)}$ and each component of $H-V(H^-)$ (which is isomorphic to $J$) is mapped to $C-v_H$ for some $C\in\cJ_H$ such that every $C-v_H$ is the image of exactly one component isomorphic to $J$ of $H-V(H^-)$.
\endproof

%\vspace{-.8cm}

\section*{Acknowledgement}
\noindent
The second author thanks Jaehoon Kim for stimulating discussions at early stages of the project.

\bibliographystyle{amsplain_v2.0customized}
\bibliographystyle{amsplain}
\bibliography{References_approx}

%\vfill
%
%\small
%\vskip2mm plus 1fill
%\noindent
%Version \today{}
%\bigbreak
%
%
%
%\noindent
%Stefan Glock
%{\tt <s.glock@bham.ac.uk>}\\
%Felix Joos
%{\tt <f.joos@bham.ac.uk>}\\
%School of Mathematics\\ 
%University of Birmingham\\
%United Kingdom

\end{document}